\documentclass[12pt]{article}

\usepackage{amsmath}
\usepackage{amssymb}
\usepackage{amsfonts}
\usepackage[colorlinks=true,linkcolor=RawSienna,citecolor=RawSienna,urlcolor=RawSienna,bookmarksopen=true,pdfstartview=FitB]{hyperref}
\usepackage[dvipsnames]{xcolor}
\usepackage[top=1.1in, bottom=1.1in, left=1in, right=1in]{geometry}
\usepackage[onehalfspacing]{setspace}
\usepackage{natbib}

\usepackage{algorithm}
\usepackage{algpseudocode}
\usepackage{graphicx}
\usepackage{caption}
\usepackage{subcaption}
\usepackage{enumerate}

\usepackage{bm}

\newtheorem{theorem}{Theorem}[section]
\newtheorem{corollary}{Corollary}[section]

\newtheorem{lemma}{Lemma}[section]
\newtheorem{proposition}{Proposition}[section]

\begin{document}

\title{\vspace{-1.2cm} Enhancing Signal Proportion Estimation Through Leveraging Arbitrary Covariance Structures}
\author{Jingtian Bai and Xinge Jessie Jeng\thanks{%
		Address for correspondence: Department of Statistics, North Carolina State University, 2311 Stinson Dr., Raleigh, NC 27695-8203, USA.  Email: \texttt{xjjeng@ncsu.edu}.}\\
Department of Statistics, North Carolina State University }
\date{}
\maketitle

\begin{abstract}

Accurately estimating the proportion of true signals among a large number of variables is crucial for enhancing the precision and reliability of scientific research. Traditional signal proportion estimators often assume independence among variables and specific signal sparsity conditions, limiting their applicability in real-world scenarios where such assumptions may not hold. This paper introduces a novel signal proportion estimator that leverages arbitrary covariance dependence information among variables, thereby improving performance across a wide range of sparsity levels and dependence structures. Building on previous work that provides lower confidence bounds for signal proportions, we extend this approach by incorporating the principal factor approximation procedure to account for variable dependence. 
Our theoretical insights offer a deeper understanding of how signal sparsity, signal intensity, and covariance dependence interact. By comparing the conditions for estimation consistency before and after dependence adjustment, we highlight the advantages of integrating dependence information across different contexts.
This theoretical foundation not only validates the effectiveness of the new estimator but also guides its practical application, ensuring reliable use in diverse scenarios. 
Through extensive simulations, we demonstrate that our method outperforms state-of-the-art estimators in both estimation accuracy and the detection of weaker signals that might otherwise go undetected.

\medskip

\textit{Keywords}: Adaptive estimation, Estimation consistency, Phase diagram, Lower bound estimators, Principal factor approximation, Sparse inference

\end{abstract}

%%%%%%%%%%%%%%%%%%%%% screening

\section{Introduction}

Estimating the proportion of signal variables among a large number of noise variables is a fundamental problem in many scientific fields, including genomics, neuroscience, and environmental science. In these applications, datasets often contain a mix of truly informative variables (signals) and irrelevant ones (noise). Accurate estimation of the signal proportion provides crucial insights into the underlying structure of the data and informs the choice and calibration of statistical models, ultimately leading to more robust scientific conclusions. In this paper, we consider the multivariate normal model
\begin{equation}\label{def:model}
(X_1,\cdots,X_m)^T\sim N(\bm{\mu},{\bm{\Sigma}}),
\end{equation}
where $\bm{\mu} = (\mu_1,\cdots,\mu_m)^T$ is a sparse mean vector, with $\mu_i = 0$ indicating a noise variable and $\mu_i \ne 0$ indicating a signal. The covariance matrix $\bm{\Sigma}$ is assumed to be an arbitrary correlation matrix. This flexible model has been widely adopted for its generality and relevance in high-dimensional data analysis.

Let $I_0 = \{i \in \{1,\cdots, m\}: \mu_i = 0\}$ and $I_1 = \{i\in \{1,\cdots, m\}: \mu_i \ne 0\}$ denote the index sets of noise and signal variables, respectively. Our objective is to estimate the signal proportion,
\[
\pi = |I_1|/m.
\]
This estimation problem poses several key challenges. First, the distributions of signal variables are only partially specified due to the unknown nonzero means. Second, the performance of estimators can vary substantially with the unknown signal proportion $\pi$ \citep{MR06}. Third, arbitrary dependencies in $\bm{\Sigma}$ can severely degrade the performance of methods designed under independence assumptions \citep{jeng2023estimating}.

A variety of signal proportion estimators have been proposed, particularly in the context of multiple testing. Foundational work in \cite{storey2002direct, GW04, MR06, JC07, jin2008proportion} demonstrated that accurate estimation of $\pi$ can significantly improve the performance multiple testing procedures. For a comprehensive review, \cite{chen2019uniformly} categorizes estimators by methodology, model assumptions, consistency properties, and practical limitations.

In recent years, the problem of signal proportion estimation has gained renewed attention in the context of high-dimensional false negative control (HD-FNC), which targets the discovery of weak signals. In HD-FNC procedures, accurate estimation of $\pi$ is a crucial step for effectively retaining signals that would otherwise be missed \cite{jeng2016rare, jeng2019variable, jeng2019efficient}. This has motivated the development of methods that are adaptive to both unknown signal sparsity levels and general dependence structures \cite{jeng2022weak, jeng2023estimating}.

Although signal proportion estimation is conceptually related to variable selection, the two problems are methodologically distinct. A variable selection procedure that consistently identifies the true set of signal variables can naturally yield an estimate of $\pi$ by counting the selected signals. However, achieving selection consistency typically requires strong assumptions, such as a minimum signal strength or strict sparsity conditions. In contrast, consistent estimation of the signal proportion can often be achieved under much weaker assumptions, without requiring perfect separation between signal and noise variables. This fundamental distinction is well illustrated by the phase diagram framework in high-dimensional inference; see, for example \citep{donoho2004higher, tony2011optimal,  arias2011global, ji2012ups, cai2017large, arias2017distribution, gao2021concentration, donoho2015special, ji2014rate, jin2017phase, chen2019two, chen2023testing}.

Among existing proportion estimation approaches, methods that provide lower confidence bounds for $\pi$ are especially useful when conservative estimates are needed. \cite{MR06} introduced lower $100(1-\alpha)\%$ confidence bounds using empirical processes of $p$-values under independence. Later work, including \cite{blanchard2020post}, \cite{katsevich2020simultaneous}, and \cite{jeng2023estimating}, extended these results to dependent settings, enhancing applicability. In particular, \cite{jeng2023estimating} developed an omnibus lower bound estimator that adapts to a wide range of dependence structures while preserving lower bound guarantees, and showed through theoretical and empirical results that it is more robust and efficient than existing alternatives.

Despite these advances, existing adaptive methods often do not fully leverage known dependence structures. In many real-world applications, prior knowledge about dependence is available and can be exploited for improved estimation. Examples include:

\begin{itemize}
	\item Gene expression studies, where known biological pathways inform gene-gene correlations;
	\item Genome-wide association studies (GWAS), where linkage disequilibrium patterns define genetic correlations;
	\item Functional MRI data, where brain connectivity maps inform spatial dependence;
	\item Climate models, where physical laws and geography dictate spatio-temporal dependencies.
\end{itemize}
In these scenarios, incorporating known dependence structures can significantly enhance the accuracy and reliability of signal proportion estimation.

In this paper, we propose a new estimator for $\pi$ that fully utilizes known dependence structures while maintaining consistency across a broad range of signal sparsity levels. To retain the lower bound property, we develop two complementary approaches: one that strictly guarantees a lower bound at the cost of higher computational burden, and another that is computationally efficient and asymptotically conservative under mild conditions. Both approaches integrate the principal factor approximation (PFA) of \cite{fan2012estimating} to enhance the signal-to-noise ratio by leveraging dependence information, resulting in substantial power gains.

We present a theoretical analysis that examines how signal sparsity, signal strength, and covariance dependence jointly affect the performance of signal proportion estimation. To facilitate this analysis, we extend the phase diagram framework, which was originally developed to characterize the fundamental limits of signal detection and classification in high-dimensional inference \citep{donoho2004higher, arias2011global, ji2012ups, cai2017large, chen2019two, gao2021concentration}, to the context of signal proportion estimation. This extension yields new insights into the regions where consistent estimation is achievable and where leveraging dependence information can provide significant advantages. Building on the new insight, we introduce a novel analytic tool that quantifies the potential benefit of dependence adjustment, offering practical guidance on when such methods are likely to outperform existing approaches.

Through extensive simulations, we compare the proposed estimators with state-of-the-art methods under diverse dependence scenarios—including those derived from real data and commonly used dependence models. Our results show that the proposed estimators outperform existing methods when dependence-adjustment is advantageous, while maintaining comparable performance otherwise.

The remainder of the paper is organized as follows. Section \ref{sec:family} reviews existing lower bound estimators under both independence and dependence. Section \ref{sec:true_factor} introduces our new dependence-adjusted estimator using the PFA framework. Section \ref{sec:estimated_factor} extends the method to estimated factors. Section \ref{sec:simulation} presents simulation results. Section \ref{sec:conclusion} concludes the paper. Detailed proofs are provided in the Supplementary Materials. 

\section{Some lower bound estimators}
\label{sec:family}

In this section, we review the family of lower bound estimators originally studied in \cite{MR06} under independence and further investigated in \cite{jeng2023estimating} under dependence.  A slightly simplified version of this family of estimators is presented here.
Given observations $x_1, \ldots, x_m$ of the test statistics $X_1, \ldots, X_m$, the family of estimators are $\{\hat{\pi}_\theta({\mathbf{x}}), \theta\in[0,1]\}$, where each member indexed by a parameter $\theta\in[0,1]$ is defined as 
\begin{equation} \label{def:pi_hat_theta}
\hat{\pi}_\theta({\mathbf{x}})=\sup\limits_{t>0}\dfrac{m^{-1}\sum\limits_{i=1}^{m}\mathbb{I}\{x_i> t\}-\bar{\Phi}(t)-{c}_{m,\theta}(\alpha; F_{\mathbf{x}}) [\bar{\Phi}(t)]^\theta}{1-\bar{\Phi}(t)}.
\end{equation}
In equation (\ref{def:pi_hat_theta}), 
$\bar{\Phi}(t)=1-\Phi(t)$ and ${c}_{m,\theta}(\alpha; F_{ \mathbf{x}})$ is a bounding sequence constructed based on $F_{ \mathbf{x}}$, the joint null distribution of $x_1, \ldots, x_m$. 
Specifically, under model (\ref{def:model}), $F_{ \mathbf{x}} = N(0, {\bm{\Sigma}})$, which is known a priori, and ${c}_{m,\theta}(\alpha; F_{ \mathbf{x}})$ is constructed in the following two steps: \\
(a) Let 
\begin{equation} \label{def:V_Z0_tilde}
V_{m,\theta}({\mathbf{x}}^0)=\sup\limits_{t>0}\dfrac{\left|m^{-1}\sum\limits_{i=1}^{m}\mathbb{I}\{x_i^0> t\}-\bar{\Phi}(t)\right|}{[\bar{\Phi}(t)]^\theta}, 
\end{equation}
where $ (x_1^0, \ldots, x_m^0)^T \sim F_{ \mathbf{x}} = N(\mathbf{0}, {\bm{\Sigma}})$. \\
(b) For a given $\alpha \in (0, 1)$, determine $c_{m,\theta}(\alpha; F_{ \mathbf{x}})$ such that the function $m c_{m,\theta}(\alpha; F_{ \mathbf{x}})$ is monotonically non-decreasing in $m$ and satisfies
\begin{equation} \label{def:c_m}
\text{P}(V_{m,\theta}({\mathbf{x}}^0)>c_{m,\theta}(\alpha; F_{ \mathbf{x}}))<\alpha.
\end{equation}

As shown in the above,  $c_{m,\theta}(\alpha; F_{ \mathbf{x}})$ serves as an upper bound for $V_{m,\theta}({\mathbf{x}}^0)$ at the confidence level of $1-\alpha$, where $V_{m,\theta}({\mathbf{x}}^0)$ depends on the joint null distribution of $(x_1^0, \ldots, x_m^0)$. Intuitively, $c_{m,\theta}(\alpha; F_{ \mathbf{x}})$ represents the typical range of $\left[m^{-1}\sum_{i=1}^{m}\mathbb{I}\{x_i> t\}-\bar{\Phi}(t)\right]/ [\bar{\Phi}(t)]^\theta$
when all variables are noise. Thus, if the observed value of $m^{-1}\sum_{i=1}^{m}\mathbb{I}\{x_i> t\}$ exceeds this typical range, it indicates evidence for the presence of signals. Furthermore, the property in (\ref{def:c_m}) implies that $\hat{\pi}_\theta({\mathbf{x}})$, as defined in (\ref{def:pi_hat_theta}), serves as a lower bound estimate for $\pi$. The detailed proofs for this property are provided in Lemma \ref{lemma:lowerbd_x} in the Supplementary Materials.%\hyperref[sec:appendix]{Appendix}. Section \ref{sec:appendix}. 

In practice, $c_{m,\theta}(\alpha; F_{ \mathbf{x}})$ can be numerically determined as the $(1-\alpha)$-quantile of $V_{m,\theta}({\mathbf{x}}^0)$. Additional details regarding the implementation are provided in Algorithm 2 in Section \ref{sec:algorithm}. 
This version of the $\hat{\pi}_\theta({\mathbf{x}})$  estimators
is designed for signals with one-sided effects, and minor modifications can be made to accommodate two-sided signal effects.

This family of estimators have been studied under independence in \cite{MR06}, where the most powerful estimator in the family was identified as $\hat \pi_{0.5}(\mathbf{x})$.  
\cite{jeng2023estimating} investigated the $\hat{\pi}_\theta({\mathbf{x}})$ family under dependence and  found that no single estimator within this family consistently outperforms others across different dependence scenarios. Motivated by this finding, \cite{jeng2023estimating} proposed an adaptive lower bound estimator that is better suited to arbitrary covariance dependence. The adaptive estimator is of the form
\begin{equation} \label{def:pi_hat_x_tilde}
\hat{\pi}({\mathbf{x}})=\max\{\hat{\pi}_{0.5}({\mathbf{x}}),~\hat{\pi}_1({\mathbf{x}})\},
\end{equation}
where $\hat{\pi}_{0.5}({\mathbf{x}})$ and $\hat{\pi}_1({\mathbf{x}})$ are constructed by (\ref{def:pi_hat_theta}) with $\theta=0.5$ and $1$, respectively. 

It has been shown theoretically that for a discretized version of the $\hat{\pi}_\theta({\mathbf{x}})$ family, referred to as $\hat{\pi}^*_\theta({\mathbf{x}})$, the discretized version of $\hat{\pi}({\mathbf{x}})$, denoted as $\hat{\pi}^*({\mathbf{x}})$, is more powerful than any individual member within the $\hat{\pi}^*_\theta({\mathbf{x}})$ family under arbitrary covariance dependence (Theorem 2.4 in \cite{jeng2023estimating}).
More precisely, the discretized $\hat{\pi}^*_\theta({\mathbf{x}})$ is defined by the same procedure as in (\ref{def:pi_hat_theta}) - (\ref{def:c_m}), with the operation ``$\sup\limits_{t>0}$" replaced by ``$\max\limits_{t\in\mathbb{T}}$", where $\mathbb{T}=[\sqrt{\log \log m},\sqrt{5\log m}]\cap\mathbb{N}$. Furthermore, the finite sample performance of $\hat{\pi}({\mathbf{x}})$ has been extensively studied and compared with several representative and state-of-the-art methods in \cite{jeng2023estimating}, demonstrating that $\hat{\pi}({\mathbf{x}})$ exhibits superior power and robustness across various dependence structures and sparsity levels.

\section{Enhancing Proportion Estimation with Dependence}
\label{sec:true_factor}

We aim to develop a new powerful and adaptive method that can consistently estimate a wide range of $\pi$ values while effectively leveraging the dependence structure in ${\bm{\Sigma}}$. The development is organized into three steps in this section: 1. Assess the joint effects of signal sparsity, signal intensity, and covariance dependence on existing methods. 2. Obtain dependence-adjusted statistics using principal factor approximation.
3. Derive our dependence-adjusted proportion estimator.

\subsection{Analyzing dependence-involved joint effects} \label{sec:pi_x}

To facilitate joint analysis on the effects of signal sparsity, signal intensity, and covariance dependence on the estimation problem, we employ the following calibrations of the model parameters: 

Let 
\begin{equation} \label{def:gamma}
\pi = m^{-\gamma}, \quad \gamma \in (0, 1);
\end{equation}
\begin{equation} \label{def:h}
\mu_i = \mu \cdot \mathbb{I}\{i\in I_1\}, \quad  \mu = \sqrt{2 h \log m}, \quad h>0;
\end{equation}
and 
\begin{equation} \label{def:eta}
{1\over m^2}\sum_{i=1}^{m}\sum_{j=1}^{m}|\Sigma_{ij}| = m^{-\eta}, \quad \eta \in [0, 1].
\end{equation}
The calibrations in (\ref{def:gamma}) -(\ref{def:h}) have been widely used in high-dimensional sparse inference to describe various levels of signal sparsity and intensity through $\gamma$ and $h$ (see, e.g,  \cite{donoho2004higher, tony2011optimal, ji2012ups,  arias2011global, cai2017large, jin2017phase}). The calibration in (\ref{def:eta}) was recently proposed to quantify the overall covariance dependence through parameter $\eta$ \citep{jeng2024weak}. Through (\ref{def:gamma})-(\ref{def:eta}), the effects of signal sparsity, signal intensity, and covariance dependence are represented by the three new parameters on comparable constant scales, which facilitates theoretical insights into their combined effects  on signal proportion estimation.

In this section, we extend the study of the existing estimator  $\hat{\pi}({\mathbf{x}})$ to our model setting, incorporating parameter calibrations from (\ref{def:gamma}) - (\ref{def:eta}). 
The following Proposition \ref{thm:pi_hat_x} explicate the joint effects of these model parameters on the consistency of the discretized $\hat{\pi}^*({\mathbf{x}})$. This result serves as a benchmark for evaluating the potential advantages of our new developments in dependence-adjusted proportion estimation as shown in Section \ref{sec:pi_z}. 
For simplicity of notation, the lowercase letters 
$x_i$, and $x_i^0$ are used to denote random variables, provided there is no ambiguity in meaning. Detailed proof of Proposition \ref{thm:pi_hat_x} is provided in Section \ref{sec:proof_prop3.1}. 

\begin{proposition} \label{thm:pi_hat_x}
	Consider model (\ref{def:model}) along with the calibrations in (\ref{def:gamma}) - (\ref{def:eta}). Given $\alpha>0$, 
	we have 
	\begin{equation} \label{eq:lowerbd_x}
	\text{P}(\hat{\pi}^*({\mathbf{x}}) < \pi) \ge 1-\alpha.
	\end{equation}
	Moreover, let $c^*_{m,1}(\alpha; F_{{\mathbf{x}}}) = C_{1,\alpha}\sqrt{\log m}$ for some constant $C_{1,\alpha}$ and $c^*_{m,0.5}(\alpha; F_{{\mathbf{x}}}) \\= C_{0.5,\alpha}\sqrt{m^{-\eta}\log m}$ for some constant $C_{0.5,\alpha}$. With $\alpha = \alpha_m \to 0$, we have $\text{P}((1-\epsilon)\pi < \hat{\pi}^*({\mathbf{x}}) < \pi)\rightarrow1$ for any constant $\epsilon>0$ if, and (almost) only if
	\begin{equation} \label{cond:h}
	h> [\gamma - ({\eta}-\gamma)_+]_+.
	\end{equation}
\end{proposition}  

Result in (\ref{eq:lowerbd_x}) shows that $\hat{\pi}^*({\mathbf{x}})$ is a lower bound estimator of the true $\pi$. The consistency of $\hat{\pi}^*({\mathbf{x}})$ can be achieved under condition (\ref{cond:h}), which reveals that  the sufficient and almost necessary signal strength, represented by $h$, increases with $\gamma$ and decreases with $\eta$. In other words, as signals become sparser ($\gamma$ increases) and/or covariance dependence strengthens ($\eta$ decreases), a higher signal intensity is required for consistency. It is worth noting that this relationship is not strictly linear but rather piecewise linear. Specifically, for a given level of dependence strength $\eta$, when signals are relatively denser such that $\gamma < \eta$, the condition is $h>[2\gamma - \eta]_+$. Conversely, when signals are relatively sparser with $\gamma \geq \eta$, the condition is $h>\gamma$. %We refer to the condition as almost necessary because the consistency of $\hat{\pi}^*({\mathbf{x}})$ fails to hold when $h<[\gamma - ({\eta}-\gamma)_+]_+$.

{We assume equal signal strength, as specified in (\ref{def:h}), to facilitate the derivation of the sufficient and almost necessary condition given in (\ref{cond:h}). In the more general setting where signal strengths vary, we define the minimum signal strength as $\mu_{\min} = \min\limits_{i\in I_1} \mu_i = \sqrt{2h\log m}$. Under this formulation, the sufficiency of condition (\ref{cond:h}) continues to hold.}

\begin{figure}[!h]
	\centering
	\includegraphics[width=1\textwidth]{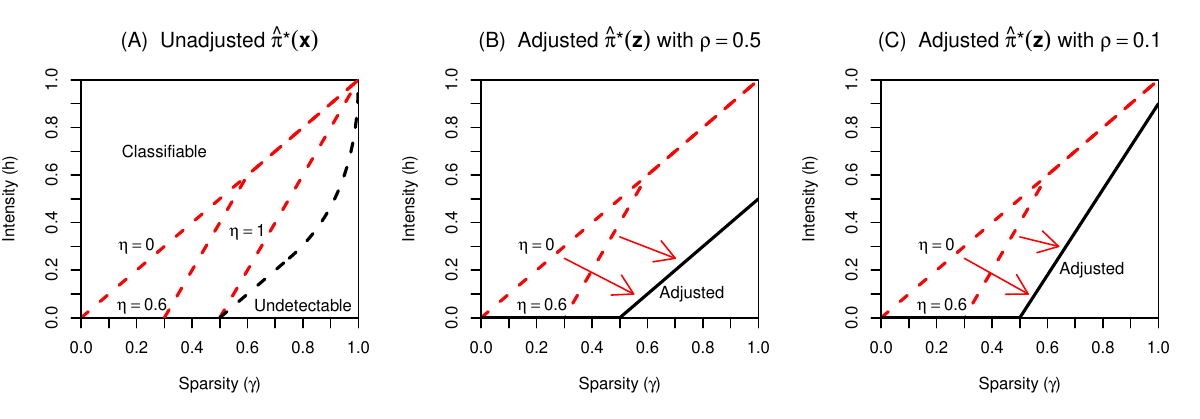}
	\caption{Estimable regions for the unadjusted estimator $\hat{\pi}^*({\mathbf{x}})$ (Panel (A)) and for the adjusted estimator $\hat{\pi}^*({\mathbf{z}})$ with a block diagonal correlation matrix (Panels (B) and (C)). The signal proportion can be consistently estimated by $\hat{\pi}^*({\mathbf{x}})$ if the model parameters fall within the region depicted by the dashed red line in Panel (A), and by $\hat{\pi}^*({\mathbf{z}})$ if they fall within the region depicted by the solid  black line in Panel (B) and (C). } \label{fig:diagram}
\end{figure}

Condition (\ref{cond:h}) can be illustrated in a two-dimensional diagram (Figure \ref{fig:diagram}(A)) as the dashed red lines that moves with the dependence parameter $\eta$. As the covariance dependence weakens ($\eta$ increases), the estimable region of $\hat{\pi}^*({\mathbf{x}})$ enlarges, making consistent estimation easier. 

Figure \ref{fig:diagram}(A) also illustrates some related results in high-dimensional inference for references. These results were derived under model (\ref{def:model}) with similar parameter calibrations in (\ref{def:gamma})-(\ref{def:h}) under independence ($\eta=1$). Specifically, the diagonal line represents the classification boundary separating classifiable and unclassificable regions \citep{donoho2004higher, cai2017large}, whereas the dotted lower curve represents the detection boundary separating detectable and undetectable regions under independence \citep{donoho2004higher, tony2011optimal}. One can see that under independence, signals in the region between these two lines are detectable but indistinguishable from the noise; however, their proportion is consistently estimable by $\hat \pi(\mathbf{x})$ if they fall above the dashed red line for $\eta=1$. More studies on related phase diagrams can be found in \cite{ji2012ups, ji2014rate, arias2011global, cai2017large, jin2017phase, chen2019two, gao2021concentration}, etc.

It is important to note that the arbitrary covariance structure in $\bm{\Sigma}$, represented through the parameter $\eta$, directly influences condition (\ref{cond:h}), making it more challenging for the estimator $\hat{\pi}^*({\mathbf{x}})$ to achieve consistency as the dependency strengthens. 
%This is because, as presented in Section \ref{sec:family}, the $\bm{\Sigma}$ matrix is used in the joint null distribution $F_{\mathbf{x}}= N(\bm{0}, \bm{\Sigma})$ to establish the bounding sequences, ${c}_{m,0.5}(\alpha; F_{ \mathbf{x}})$ and ${c}_{m,1}(\alpha; F_{ \mathbf{x}})$, which ensure the lower bound property of $\hat{\pi}^*({\mathbf{x}})$. 
This result, however, does not reflect the optimal scenario where dependence information could be leveraged to enhance the power of proportion estimation. For instance, in an extreme case where all signal variables exhibit high correlation with each other but remain independent from noise variables, stronger dependence should benefit the estimation of the signal proportion.

This consideration motivates us to develop a dependence-adjusted approach that can effectively harness the dependence information in $\bm{\Sigma}$ for proportion estimation. Such an approach would be successful under a more relaxed condition than (\ref{cond:h}), allowing for improved estimation in the presence of strong dependence among the variables.

\subsection{Dependence adjustment via principal factor approximation} \label{sec:PFA}

To achieve dependence-adjusted proportion estimation, we employ the principal factor approximation (PFA) approach developed in \cite{fan2012estimating}. Although originally designed for false discovery rate (FDR) control under dependence, we integrate this approach into the lower bound estimation procedure outlined in Section \ref{sec:family} to substantially enhance the estimation power. 

The PFA approach is carried out as follows. First, perform spectral decomposition on the correlation matrix ${\bm{\Sigma}}$ as 
\begin{equation}\label{def:lambda}
{\bm{\Sigma}}=\sum_{i=1}^{m}\lambda_i\bm{\gamma}_i\bm{\gamma}_i^T,
\end{equation}
where $\lambda_1\geq\cdots\geq\lambda_m$ are eigenvalues of ${\bm{\Sigma}}$ and $\bm{\gamma}_1,\cdots,\bm{\gamma}_m$ are corresponding orthogonal eigenvectors. 
For an integer $k (=k_m)$ that may increase with $m$, decompose the right-hand side of (\ref{def:lambda}) into two parts:
\[
{\bm{\Sigma}}=\sum_{i=1}^{k}\lambda_i\bm{\gamma}_i\bm{\gamma}_i^T+\sum_{i=k+1}^{m}\lambda_i\bm{\gamma}_i\bm{\gamma}_i^T=\mathbf{LL}^T+\mathbf{A},
\]
where $\mathbf{L}=(\sqrt{\lambda_1}\bm{\gamma}_1,\cdots,\sqrt{\lambda_k}\bm{\gamma}_k)$ is an $m\times k$ matrix,  $\mathbf{A}=\sum_{i=k+1}^{m}\lambda_i\bm{\gamma}_i\bm{\gamma}_i^T$ with $\|\mathbf{A}\|_F^2=\lambda_{k+1}^2+\cdots+\lambda_m^2$, and $\|\cdot\|_F$ stands for the Frobenius norm.
Then $X_1,\cdots,X_m$ can be rewritten as    
\begin{equation} \label{def:factor_model_1}
{\mathbf{X}}=\bm{\mu}+\mathbf{LW}+\mathbf{K}
\end{equation}    
or
\begin{equation}\label{def:factor_model_2}
X_i=\mu_i+\mathbf{b}_i^T\mathbf{W}+K_i, \qquad i=1,\cdots,m,
\end{equation}
where  $\mathbf{W}=(W_1,\cdots,W_k)^T\sim N(\bm{0}_k,\mathbf{I}_k)$ are the principal factors, the $i$th row of $\mathbf{L}$ is  $\mathbf{b}_i=(b_{i1},\cdots,b_{ik})^T$, $(b_{1j},\cdots,b_{mj})^T = \sqrt{\lambda}_j\bm{\gamma}_j$, and $\mathbf{K}=(K_1,\cdots,K_m)^T\sim N(\bm{0}_m,\mathbf{A})$ are the random errors, independent of the principal factors. In practice, the number of principal factors is determined as the smallest $k$ such that $\sqrt{\lambda_{k+1}^2+\cdots+\lambda_m^2}/(\lambda_1+\cdots+\lambda_m) <\epsilon$ for a prefixed small $\epsilon$. 

Define the dependence-adjusted statistics 
\begin{equation}\label{def:Z}
Z_i=a_i(X_i-\mathbf{b}_i^T\mathbf{W}), \qquad i=1,\cdots,m,
\end{equation}
where $a_i=(1-\|\mathbf{b}_i\|_2^2)^{-1/2}$. Then 
\begin{equation}\label{def:model.Z}
(Z_1,\cdots,Z_m)^T\sim N((a_1\mu_1,\cdots,a_m\mu_m)^T,\tilde{\bm{\Sigma}}),
\end{equation}
where $\tilde{\bm{\Sigma}}$ is the correlation matrix of $Z_1,\cdots, Z_m$, obtained from the covariance matrix $\mathbf{A}=\sum_{i=k+1}^{m}\lambda_i\bm{\gamma}_i\bm{\gamma}_i^T$.

Therefore, for $i \in I_0$, $Z_i \sim N(0, 1)$ as $a_i\mu_i=0$, whereas  for $i\in I_1$, $Z_i\sim N(a_i\mu_i, 1)$. Note that the signal proportion of $Z_1,\cdots, Z_m$ remains the same as $\pi = |I_1|/m$. On the other hand, the signal intensities of $Z_i$, $i \in I_1$, are amplified by $a_i (>1)$. In other words, accounting for the principal factors of $\bm{\Sigma}$ results in less variability in the remaining error terms $(K_1,\cdots,K_m)^T\sim N(\bm{0}_m,\mathbf{A})$, which contributes to the increased signal-to-noise ratio of $Z_1,\cdots, Z_m$.

It is important to note that $Z_1,\cdots, Z_m$ have a different correlation matrix $\tilde{\bm{\Sigma}}$. Although the remaining error terms carry less variability, as indicated by the smaller Frobenius norm of $\mathbf{A}$ compared to that of ${\bm{\Sigma}}$, the correlation matrix $\tilde{\bm{\Sigma}}$ derived from  $\mathbf{A}$ can be arbitrary. 
We apply the same calibration technique as in (\ref{def:eta}) to measure the degree of overall dependence in $\tilde{\bm{\Sigma}}$: let 
\begin{equation} \label{def:tilde_eta}
{1\over m^2}\sum_{i=1}^{m}\sum_{j=1}^{m}|\tilde \Sigma_{ij}| = m^{-\tilde \eta}, \quad \tilde \eta \in [0, 1].
\end{equation}
The parameter $\tilde \eta$ can be greater or smaller than the original parameter $\eta$.

\subsection{Dependence-adjusted proportion estimator} \label{sec:pi_z}

We implement the dependence-adjusted statistics $Z_1,\cdots, Z_m$ into the lower bound estimation procedure presented in Section \ref{sec:family} and obtain a new estimator $\hat \pi (\mathbf{z})$ as
\begin{equation} \label{def:pi_z}
\hat \pi (\mathbf{z}) = \max\{\hat{\pi}_{0.5}(\mathbf{z}),~\hat{\pi}_1(\mathbf{z})\},
\end{equation}
where $\hat{\pi}_\theta(\mathbf{z})$ is defined as in (\ref{def:pi_hat_theta}) with $(x_1, \ldots, x_m)$ replaced by $(z_1, \ldots, z_m)$, $c_{m,\theta}(\alpha; F_{\mathbf{x}})$ replaced by $c_{m,\theta}(\alpha/2; F_{\mathbf{z}})$, where  $F_{\mathbf{z}} = N(\bm{0}, \tilde{\bm{\Sigma}})$ is the adjusted joint null distribution of $(z_1, \cdots, z_m)$, and $(x^0_1, \cdots, x^0_m)$ replaced by $(z^0_1, \cdots, z^0_m) \sim F_{\mathbf{z}}$.

Consequently, we present the consistency of the discretized version $\hat{\pi}^*({\mathbf{z}})$. 
To that end, define 
\[
a_{\min} = \min_{1 \le i \le m}\{a_i\} \qquad \text{and} \qquad a_{\max} = \max_{1 \le i \le m}\{a_i\}.
\] 
For notation simplicity, the lowercase letters $z_i$ and $z^0_i$ are used to denote random variables, provided there is no ambiguity in meaning. The proof of the following Theorem \ref{thm:pi_hat_z} is provided in Section \ref{sec:proof_thm3.1}. 

\begin{theorem} \label{thm:pi_hat_z}
	Consider model (\ref{def:model}) and the dependence-adjusted statistics in (\ref{def:Z}), along with the calibrations in (\ref{def:gamma}) - (\ref{def:h}) and (\ref{def:tilde_eta}). Given $\alpha>0$, 
	we have 
	\begin{equation} \label{eq:lowerbd_z}
	\text{P}(\hat{\pi}^*({\mathbf{z}}) < \pi) \ge 1-\alpha.
	\end{equation}
	Moreover, let $c^*_{m,1}(\alpha/2; F_{{\mathbf{z}}}) = C'_{1,\alpha}\sqrt{\log m}$ for some constant $C'_{1,\alpha}$ and $c^*_{m,0.5}(\alpha/2; F_{{\mathbf{z}}}) = C'_{0.5,\alpha}\sqrt{m^{-\tilde \eta}\log m}$ for some constant $C'_{0.5,\alpha}$. With $\alpha = \alpha_m \to 0$, we have $\text{P}((1-\epsilon)\pi < \hat{\pi}^*({\mathbf{z}}) < \pi)\rightarrow1$ for any constant $\epsilon>0$ if
	\begin{equation} \label{cond:ah}
	a^2_{\min} \cdot h> [\gamma - ({\tilde \eta}-\gamma)_+]_+.
	\end{equation}
\end{theorem}  

Compared with condition (\ref{cond:h}) in Proposition \ref{thm:pi_hat_x}, condition (\ref{cond:ah}) differs in two important ways. First, the left-hand side is multiplied by $a^2_{\min} (>1)$, which substantially relaxes the condition. Second, $\tilde\eta$ replaces $\eta$, which can further relax the condition when $\tilde\eta > \eta$. This occurs when the dependence structure in the original data $X_1,\cdots, X_m$ is well captured by the factor model (\ref{def:factor_model_1}), thereby reducing the residual correlations in the error terms.

Because the quantities $\eta$, $\tilde\eta$, and $a_{\min}$ can be directly computed from the known correlation matrix $\bm{\Sigma}$, it is possible to assess whether condition (\ref{cond:ah}) is indeed looser than condition (\ref{cond:h}), and thus determine when the dependence-adjusted estimator is preferable to the unadjusted method. For example, if $\tilde\eta > \eta$, then (\ref{cond:ah}) is more relaxed than (\ref{cond:h}) for any $a_{\min} \ge 1$ and any unknown $\gamma$. Moreover, when $a_{\min}$ is substantially greater than 1, its squared effect in relaxing (\ref{cond:ah}) typically dominates the influence of $\tilde\eta$ and the unknown $\gamma$.

Such a practical guideline is highly desirable given its relevance for real-world applications. Below, we present an example that illustrates the insight gained from the key conditions characterizing the advantage of dependence adjustment. Additional examples under various covariance structures are provided in Section \ref{sec:sim_evaluation}.

\noindent {\bf An illustrative example}. 
Consider model (\ref{def:model}), where ${\bm{\Sigma}}$ is a block-diagonal matrix with each block having size $m^d \times m^d$, where $d \in [0,1]$. There are $m^{1-d}$ such blocks, and each block is given by ${\mathbf{B}} = \rho \mathbf{1}{m^d} \mathbf{1}{m^d}^T + (1 - \rho) \mathbf{I}_{m^d}$, with $\rho \in [0,1]$.
The following corollary compares the theoretical performance of $\hat{\pi}^*({\mathbf{x}})$ and $\hat{\pi}^*({\mathbf{z}})$ in this example. The proof is provided in detail in Section \ref{sec:proof_cor3.1}.  

\begin{corollary} \label{thm:compare}
	Consider model (\ref{def:model}), where ${\bm{\Sigma}}$ is a block-diagonal matrix with each block of size $m^d \times m^d, d \in (0,1]$, and the off-diagonal elements within each block are equal to $\rho \in (0, 1)$. Consider the calibrations in (\ref{def:gamma}) - (\ref{def:eta}) and (\ref{def:tilde_eta}). 
	
	For the unadjusted estimator $\hat{\pi}^*({\mathbf{x}})$ with $c^*_{m,1}(\alpha; F_{{\mathbf{x}}}) = C_{1,\alpha}\sqrt{\log m}$ and $c^*_{m,0.5}(\alpha; F_{{\mathbf{x}}}) = C_{0.5,\alpha}\sqrt{m^{-\eta}\log m}$, as $\alpha = \alpha_m \to 0$, $\text{P}((1-\epsilon)\pi < \hat{\pi}^*({\mathbf{x}}) < \pi)\rightarrow1$ for any $\epsilon>0$ if and (almost) only if
	\begin{equation} \label{cond:h_example}
	h> [\gamma - (1-d-\gamma)_+]_+. 
	\end{equation}
	
	For the dependence-adjusted estimator $\hat{\pi}^*({\mathbf{z}})$ with $k=m^{1-d}$, $c^*_{m,1}(\alpha/2; F_{{\mathbf{z}}}) = C'_{1,\alpha}\sqrt{\log m}$, and $c^*_{m,0.5}(\alpha/2; F_{{\mathbf{z}}}) = C'_{0.5,\alpha}\sqrt{m^{-\tilde \eta}\log m}$, as $\alpha = \alpha_m \to 0$,  $\text{P}((1-\epsilon)\pi < \hat{\pi}^*({\mathbf{z}}) < \pi)\rightarrow1$ for any $\epsilon>0$ if
	\begin{equation} \label{cond:ah_example}
	h> (1-\rho)[2\gamma - 1]_+.
	\end{equation}
\end{corollary}  

In this example, for the given $\bm{\Sigma}$, we have $\eta = 1-d-o(1)$. By setting $k$ equal to the number of blocks ($k=m^{1-d}$), we obtain $a_{\min} \ge 1/\sqrt{1-\rho}$ and $\tilde \eta = 1-o(1)$, 
which are used to derive condition (\ref{cond:ah_example}). Comparing conditions (\ref{cond:h_example}) and (\ref{cond:ah_example}), the latter is more relaxed because $(1-\rho) < 1$ and $2\gamma -1 \le \gamma - (1-d-\gamma)_+$ for any $\rho \in (0, 1)$, $d \in (0, 1]$ and $\gamma \in (0,1)$. This highlights the advantage of the dependence-adjusted estimator.  
The comparison is illustrated in Figure \ref{fig:diagram}(B) and   \ref{fig:diagram}(C), with $\rho = 0.5$ and $0.1$, respectively, and  $d=0.4$ or 1 ($\eta=0.6$ or $0$). The dashed red lines represent condition (\ref{cond:h_example}) for $\hat{\pi}^*({\mathbf{x}})$, while the solid black lines represent condition (\ref{cond:ah_example}) for $\hat{\pi}^*({\mathbf{z}})$. Both panels show that $\hat{\pi}^*({\mathbf{z}})$ achieves improved estimation by enlarging the estimable regions compared to $\hat{\pi}^*({\mathbf{x}})$. Moreover, $\hat{\pi}^*({\mathbf{z}})$ benefits more from stronger overall dependence ($\rho = 0.5$) than from weaker dependence ($\rho = 0.1$), reinforcing the idea that incorporating dependence information can substantially improve proportion estimation, and that our theoretical analysis offers a practical way to assess when this advantage exists by examining the observed values of $\eta$, $\tilde\eta$ and $a_{\min}$. 

\section{Dependence-adjusted estimator with estimated principal factors} \label{sec:estimated_factor}

In real data applications, the principal factors $\mathbf{W} = (W_1, \ldots, W_k)$ of the observed $\bm{\Sigma}$ are realized but not directly observable. 
Consequently, the dependence-adjusted proportion estimator $\hat \pi(\mathbf{z})$ cannot be directly constructed. Instead, we apply $\hat \pi(\widehat{\mathbf{z}})$ in real practice with 
\begin{equation} \label{def:z_hat}
\hat z_i = a_i(x_i - \mathbf{b}^T_i \widehat{\mathbf{w}}),
\end{equation}
where $\widehat{\mathbf{w}} = (\hat w_1, \ldots, \hat w_k)$ are the estimated principal factors of $\bm{\Sigma}$. We employ the robust $L_1$-regression method in \cite{fan2012estimating} to obtain $\widehat{\mathbf{w}}$ as follows.

Given the observed values $x_1, \ldots, x_m$, obtain
\begin{equation}\label{def:est.W}
\widehat{\mathbf{w}}=\arg\min\limits_{\bm{\beta}\in\mathbb{R}^k}\sum_{i=1}^m |x_{i}-\mathbf{b}_{i}^T\bm{\beta}|,
\end{equation}
where $\mathbf{b}_i$ is defined as in (\ref{def:factor_model_2}). The following lemma is adapted from Theorem 3 in \cite{fan2012estimating} with relaxed conditions. Its proof is presented in Section \ref{sec:proof_lemma4.1}. {Note that $k(=k_m)$ may increase with $m$.}

\begin{lemma} \label{lemma:w_hat-w}
	Consider model (\ref{def:model}) and its principal factor approximation by (\ref{def:factor_model_2}), along with estimator (\ref{def:est.W}) and the calibrations in (\ref{def:gamma}) and (\ref{def:tilde_eta}). Under the following conditions \\
	(C1)~ $C\sqrt{m^{1-\min\{\tilde \eta,~2\gamma\}}}\leq a_{\min}\sqrt{k}\leq a_{\max}\sqrt{k}=o(\sqrt{m})$ for some constant $C$.\\
	(C2)~ For a constant $d>0$,
	\begin{equation*}
	\lim\limits_{m\rightarrow\infty}\sup\limits_{\|\mathbf{u}\|_2=1}m^{-1}\sum_{i=1}^{m}\mathbb{I}{\{|\mathbf{b}_i^T\mathbf{u}|\leq d\}}=0.
	\end{equation*}
	We have 	
	\begin{equation}\label{hatW.1}
	\|\widehat{\mathbf{w}} - \mathbf{w}\|_2 = O_p\left(\sqrt{k/m}\right).
	\end{equation}
\end{lemma}    

Since implementing $\widehat{\mathbf{w}}$ in (\ref{def:z_hat}) introduces additional uncertainty into the proportion estimation procedure, we develop two distinct approaches to maintain the lower bound property and examine their estimation consistency. The first approach, referred to as the lower bound approach, strictly preserves the lower bound property but is computationally more intensive. The second approach, termed the approximation approach, is computationally faster but only retains the lower bound property asymptotically and under additional conditions.

\subsection{The lower bound approach} \label{sec:L_approach}

In this approach, given the $\widehat{\mathbf{w}}$ obtained in (\ref{def:est.W}), the proportion estimator is constructed as 
\begin{equation} \label{def:pi_L}
\hat \pi_L (\widehat{\mathbf{z}}) = \max\{\hat{\pi}_{L, 0.5}(\widehat{\mathbf{z}}),~\hat{\pi}_{L,1}(\widehat{\mathbf{z}})\},
\end{equation}
where $\hat{\pi}_{L,\theta}(\widehat{\mathbf{z}})$ is defined as in (\ref{def:pi_hat_theta}) with $(x_1, \cdots, x_m)$ replaced by $(\hat z_1, \cdots, \hat z_m)$ from (\ref{def:z_hat}),  $c_{m,\theta}(\alpha; F_{\mathbf{x}})$ replaced by $c_{m,\theta}(\alpha/2; F_{\widehat{\mathbf{z}}})$, where  $F_{\widehat{\mathbf{z}}}$ is the joint null distribution of $(\hat z_1, \ldots, \hat z_m)$, and $(x^0_1, \cdots, x^0_m)^T$ replaced by $(\hat z^0_1, \cdots, \hat z^0_m)^T \sim F_{\widehat{\mathbf{z}}}$. Because the $F_{\widehat{\mathbf{z}}}$ distribution cannot be explicitly defined, we generate $\hat z^0_i$ through $\hat z^0_i = a_i(x_i^0-\mathbf{b}_i^T \widehat{\mathbf{w}}^0)$, where $(x^0_1, \cdots, x^0_m)^T \sim N(0, \bm{\Sigma})$, and $\widehat{\mathbf{w}}^0$ is obtained by solving (\ref{def:est.W}) with $x_i$ replaced by $x^0_i$. 

To investigate the theoretical properties of $\hat \pi_L (\widehat{\mathbf{z}})$, we examine its discritized version, $\hat{\pi}^*_L(\widehat{\mathbf{z}})$, and present the following Theorem \ref{thm:pi_hat_L} on the consistency of $\hat{\pi}^*_L(\widehat{\mathbf{z}})$. To facilitate the theoretical analysis of $\hat{\pi}^*_L(\widehat{\mathbf{z}})$, we replace ``$\max\limits_{t\in\mathbb{T}}$" with ``$\max\limits_{t\in\mathbb{T}'}$", where $\mathbb{T}'=[\sqrt{2[\gamma- (\tilde \eta-\gamma)_+]_+\log m},\sqrt{5\log m}]\cap\mathbb{N}$ is a more restricted range than $\mathbb{T}$. For notation simplicity, the lowercase letters $\hat z_i$ and $\hat z^0_i$ are used to denote random variables, provided there is no ambiguity in meaning. The proof of Theorem \ref{thm:pi_hat_L} is provided in Section \ref{sec:proof_thm4.1}. 

\begin{theorem} \label{thm:pi_hat_L}
	Consider model (\ref{def:model}) and the dependence-adjusted statistics in (\ref{def:z_hat}) with the estimated principal factors in (\ref{def:est.W}), along with the calibrations in (\ref{def:gamma}) - (\ref{def:h}) and (\ref{def:tilde_eta}). Given $\alpha>0$, we have 
	\begin{equation} \label{eq:lowerbd_L}
	\text{P}(\hat\pi_L^*(\widehat{\mathbf{z}}) < \pi) \ge 1-\alpha.
	\end{equation} 
	Moreover, let $c^*_{m,1}(\alpha/2; F_{\widehat{\mathbf{z}}}) = C''_{1,\alpha}\sqrt{\log m}$ for some constant $C''_{1,\alpha}$ and $c^*_{m,0.5}(\alpha/2; F_{\widehat{\mathbf{z}}}) = C''_{0.5,\alpha}\sqrt{m^{-\tilde \eta}\log m}$ for some constant $C''_{0.5,\alpha}$.	With $\alpha = \alpha_m \to 0$, under condition (\ref{cond:ah}) and the conditions of Lemma \ref{lemma:w_hat-w} with (C1) strengthened by\\
	(C1')~ $C_1\sqrt{m^{1-\tilde \eta}}\leq a_{\min}\sqrt{k}\leq a_{\max}\sqrt{k} \leq C_2\sqrt{m^{1-2(\tilde \eta - \gamma)_+}/\log m}$ for some constant $C_1$ and $C_2$,  \\
	we have 
	\begin{equation} \label{eq:upper_L}
	\text{P}((1-\epsilon)\pi < \hat\pi_L^*(\widehat{\mathbf{z}}) < \pi)\rightarrow 1
	\end{equation}
	for any constant $\epsilon>0$. 
\end{theorem}  
The result in (\ref{eq:lowerbd_L}) demonstrates that $\hat\pi_L^*(\widehat{\mathbf{z}})$ maintains the lower bound property without additional conditions. However, additional conditions (C1') and (C2) are required for the consistency of $\hat\pi_L^*(\widehat{\mathbf{z}})$, as compared to Theorem \ref{thm:pi_hat_z}.

\subsection{The approximation approach} \label{sec:A_approach}

The second approach simplifies the estimation procedure by streamlining the generation of the critical sequence. Specifically, the estimator from this approach is denoted as  
\begin{equation} \label{def:pi_A}
\hat \pi_A (\widehat{\mathbf{z}}) = \max\{\hat{\pi}_{A, 0.5}(\widehat{\mathbf{z}}),~\hat{\pi}_{A,1}(\widehat{\mathbf{z}})\},
\end{equation}
where $\hat{\pi}_{A,\theta}(\widehat{\mathbf{z}})$ is defined as in (\ref{def:pi_hat_theta}) with $(x_1, \cdots, x_m)$ replaced by $(\hat z_1, \cdots, \hat z_m)$ from (\ref{def:z_hat}),  $c_{m,\theta}(\alpha; F_{\mathbf{x}})$ replaced by $c_{m,\theta}(\alpha/2; F_{\mathbf{z}})$, and $(x^0_1, \cdots, x^0_m)^T$ replaced by $(z^0_1, \cdots, z^0_m)^T \sim F_{\mathbf{z}}$. 

Since the distribution $F_{\mathbf{z}} = N(0, \tilde{\bm{\Sigma}})$ can be directly specified, with the correlation matrix $\tilde{\bm{\Sigma}}$ derived from the covariance matrix $\mathbf{A}=\sum_{i=k+1}^{m}\lambda_i\bm{\gamma}_i\bm{\gamma}_i^T$, generating  $c_{m,\theta}(\alpha/2; F_{\mathbf{z}})$ is simpler than generating  $c_{m,\theta}(\alpha/2; F_{\widehat{\mathbf{z}}})$ in the lower bound approach. 
However, this simplification comes with a trade-off: the lower bound property is preserved only asymptotically and under additional conditions as shown in the following theorem. 

Let the discretized version be denoted as $\hat{\pi}^*_A(\widehat{\mathbf{z}})$. To facilitate the theoretical analysis of $\hat{\pi}^*_A(\widehat{\mathbf{z}})$, we replace ``$\max\limits_{t\in\mathbb{T}}$" with ``$\max\limits_{t\in\mathbb{T}'}$". Detailed proof of the following theorem is presented in Section \ref{sec:proof_thm4.2}. 

\begin{theorem} \label{thm:pi_hat_A}
	Consider model (\ref{def:model}) and the dependence-adjusted statistics in (\ref{def:z_hat}) with the estimated principal factors in (\ref{def:est.W}), along with the calibrations in (\ref{def:gamma}) - (\ref{def:h}) and (\ref{def:tilde_eta}). Given $\alpha>0$,  under the conditions of Lemma \ref{lemma:w_hat-w} with (C1) strengthened by \\
	(C1'')~ $C\sqrt{m^{1-\tilde \eta}}\leq a_{\min}\sqrt{k}\leq a_{\max}\sqrt{k} = o(\sqrt{m^{1-2(\tilde \eta - \gamma)_+}/\log m})$ for some constant $C$, \\
	we have
	\begin{equation} \label{eq:lowerbd_A}
	\text{P}\left(\hat\pi_A^*(\widehat{\mathbf{z}}\right) < \pi(1+\delta)) \ge 1-\alpha +o(1)
	\end{equation} 
	for any constant $\delta>0$. Moreover, let $c^*_{m,1}(\alpha/2; F_{{\mathbf{z}}}) = C'_{1,\alpha}\sqrt{\log m}$ and $c^*_{m,0.5}(\alpha/2; F_{{\mathbf{z}}}) = C'_{0.5,\alpha}\sqrt{m^{-\tilde \eta}\log m}$ as in Theorem \ref{thm:pi_hat_z}. With $\alpha = \alpha_m \to 0$, under condition (\ref{cond:ah}), we have 
	\begin{equation} \label{eq:upperbd_A}
	\text{P}((1-\epsilon)\pi < \hat\pi_A^*(\widehat{\mathbf{z}}) < \pi)\rightarrow 1
	\end{equation}
	for any constant $\epsilon>0$.
\end{theorem}  

\subsection{Algorithm and Implementation} \label{sec:algorithm}

In this section, we provide algorithms to obtain the dependence-adjusted lower bound estimator $\hat{\pi}_L(\widehat{\mathbf{z}})$ and the approximated version $\hat{\pi}_A(\widehat{\mathbf{z}})$. Specifically, we provide algorithms to compute the dependence-adjusted statistics, to construct the bounding sequences, and to estimate the signal proportion.

\begin{algorithm}[!h]
	\small
	\caption{Compute dependence-adjusted statistics}
	\label{Alg:pfa}
	\begin{enumerate}
		\item \textbf{Input}: $m$ observed statistics: $X_i$ for $i=1,\cdots,m$; the eigenvalues $\lambda_1\geq\cdots\geq\lambda_m$ and eigenvectors $\bm{\gamma}_1,\cdots,\bm{\gamma}_m$ of the $m\times m$ correlation matrix ${\bm{\Sigma}}$; and the number of principal factors $k$.
		% \item Obtain eigenvalues $\lambda_1\geq\cdots\geq\lambda_m$ and eigenvectors $\bm{\gamma}_1,\cdots,\bm{\gamma}_m$ of ${\bm{\Sigma}}$ using singular value decomposition.
		\item Compute 
		\[
		\mathbf{L}=(\sqrt{\lambda_1}\bm{\gamma}_1,\cdots,\sqrt{\lambda_k}\bm{\gamma}_k)
		\] 
		and obtain $\mathbf{b}_i$ as the $i$th row of $\mathbf{L}$ for $i=1,\cdots,m$.
		\item Compute
		\[ \mathbf{A}=\sum_{i=k+1}^{m}\lambda_i\bm{\gamma}_i\bm{\gamma}_i^T
		\] 
		\item Estimate principal factors
		\[
		\widehat{\mathbf{w}}=\arg\min\limits_{\bm{\beta}\in\mathbb{R}^k}\sum_{i=1}^{m}|X_i-\mathbf{b}_i^T\bm{\beta}|.
		\]
		\item Compute the dependence-adjusted statistics $\widehat{Z}_i$ for $i=1,\cdots,m$:
		\[
		a_i=(1-\|\mathbf{b}_i\|_2^2)^{-1/2};\quad\widehat{Z}_i=a_i(X_i-\mathbf{b}_i^T\widehat{\mathbf{w}}).
		\]
		\item \textbf{Output}: the dependence-adjusted statistics
		$\widehat{Z}_i$ for $i=1,\cdots,m$, and the correlation matrix $\tilde{\bm{\Sigma}}$ obtained from the covariance matrix $\mathbf{A}$.
	\end{enumerate}
\end{algorithm}

\begin{algorithm}[!h]
	\small
	\caption{Construct the bounding sequences for two-sided signals}
	\label{Alg:boundseq}
	\begin{enumerate}
		\item  \textbf{Input}: $N$ sets of statistics from the joint null distribution:  $X_{r,i}^0$ for $r=1,\cdots,N$ and $i=1,\cdots,m$; and the lower bound level $\alpha$. 		
		\item  \textbf{For }$r=1,\cdots,N$ \textbf{do}:
		\begin{enumerate}
			\item Rank the $r$-th set of statistics: $|X_{r,(1)}^0|>|X_{r,(2)}^0|>\cdots>|X_{r,(m)}^0|$.			
			\item Compute
			\[
			V_{r,0.5}=\max\limits_{1<i< m}\dfrac{|i/m-2\bar{\Phi}(|X_{r,(i)}^0|)|}{\sqrt{2\bar{\Phi}(|X_{r,(i)}^0|)}};\quad V_{r,1}=\max\limits_{1< i<m}\dfrac{|i/m-2\bar{\Phi}(|X_{r,(i)}^0|)|}{2\bar{\Phi}(|X_{r,(i)}^0|)}.
			\]			
		\end{enumerate}
		\item \textbf{Output}: the bounding sequences:
		\begin{eqnarray}
		c_{m,0.5}&=&\text{the }(1-\alpha)\text{-th quantile of }V_{r,0.5}, \quad r=1,\cdots, N. \nonumber \\
		c_{m,1}&=&\text{the }(1-\alpha)\text{-th quantile of }V_{r,1}, \quad r=1,\cdots, N. \nonumber
		\end{eqnarray}	    
	\end{enumerate}
\end{algorithm}

\begin{algorithm}[!h]
	\small
	\caption{Estimate the signal proportion based on the bounding sequences}
	\label{Alg:signprop}
	\begin{enumerate}
		\item  \textbf{Input}: $m$ observed statistics: $X_i$ for $i=1,\cdots,m$; and bounding sequences $c_{m,0.5}$ and $c_{m,1}$.
		\item  Rank the statistics: $|X_{(1)}|>|X_{(2)}|>\cdots>|X_{(m)}|$.
		\item Compute
		\[
		\hat{\pi}_{0.5}=\max\limits_{1<i<m}\dfrac{i/m-2\bar{\Phi}(|X_{(i)}|)-c_{m,0.5}\sqrt{2\bar{\Phi}(|X_{(i)}|)}}{1-2\bar{\Phi}(|X_{(i)}|)};\]
		\[ \hat{\pi}_{1}=\max\limits_{1<i<m}\dfrac{i/m-2\bar{\Phi}(|X_{(i)}|)-c_{m,1}2\bar{\Phi}(|X_{(i)}|)}{1-2\bar{\Phi}(|X_{(i)}|)}.
		\]
		\item \textbf{Output}: the signal proportion estimate: $\hat{\pi}=\max\{\hat{\pi}_{0.5},\hat{\pi}_{1}\}$.	  
	\end{enumerate}
\end{algorithm}

\begin{algorithm}[!h]
	\small
	\caption{Dependence-Adjusted Proportion Estimator - Lower Bound Approach}
	\label{Alg:lowerbound}
	\begin{enumerate}
		\item \textbf{Input:} $m$ observed statistics $X_i$ for $i=1,\cdots,m$; $m\times m$ correlation matrix ${\bm{\Sigma}}$; number of principal factors $k$; number of replication sets $N$; lower bound level $\alpha$.
		\item Compute eigenvalues $\lambda_1\geq\cdots\geq\lambda_m$ and eigenvectors $\bm{\gamma}_1,\cdots,\bm{\gamma}_m$ of ${\bm{\Sigma}}$ using singular value decomposition.
		\item Input $X_i$, $\lambda_i$, $\bm{\gamma}_i$ for $i=1,\cdots,m$, and $k$, run \textbf{Algorithm \ref{Alg:pfa}}, and obtain the dependence-adjusted statistics $\widehat{Z}_i$ for $i=1,\cdots,m$.
		\item Randomly generate $N$ independent sets of statistics $(X_{r,1}^0,\cdots,X_{r,m}^0)^T \sim N((0, \cdots,0)^T,{\bm{\Sigma}})$ for $r=1,\cdots,N$.
		\item \textbf{For} $r=1,\cdots,N$: Input $X_{r,i}^0$, $\lambda_i$, $\bm{\gamma}_i$ for $i=1,\cdots,m$, and $k$, run \textbf{Algorithm \ref{Alg:pfa}}, and obtain the dependence-adjusted statistics $\widehat{Z}_{r,i}^0$ for $i=1,\cdots,m$.
		\item Input $\widehat{Z}_{r,i}^0$ for $r=1,\cdots,N$ and $i=1,\cdots,m$, and $\alpha/2$, run \textbf{Algorithm \ref{Alg:boundseq}}, and obtain the bounding sequences $\hat{c}_{m,0.5}$ and $\hat{c}_{m,1}$.
		\item Input $\widehat{Z}_i$ for $i=1,\cdots,m$, and $\hat{c}_{m,0.5}$ and $\hat{c}_{m,1}$, run \textbf{Algorithm \ref{Alg:signprop}}, and obtain the signal proportion estimate $\hat{\pi}_L$.
		\item \textbf{Output:} The dependence-adjusted lower bound estimator $\hat{\pi}_L$.
	\end{enumerate}
\end{algorithm}

\begin{algorithm}[!h]
	\small
	\caption{Dependence-Adjusted Proportion Estimator - Approximation Approach}
	\label{Alg:approxlowerbound}
	\begin{enumerate}
		\item  \textbf{Input}: $m$ observed statistics: $X_i$ for $i=1,\cdots,m$; $m\times m$ correlation matrix ${\bm{\Sigma}}$; number of principal factors $k$; number of replication sets $N$; lower bound level $\alpha$.
		\item Compute eigenvalues $\lambda_1\geq\cdots\geq\lambda_m$ and eigenvectors $\bm{\gamma}_1,\cdots,\bm{\gamma}_m$ of ${\bm{\Sigma}}$ using singular value decomposition.
		\item Input $X_i$, $\lambda_i$, $\bm{\gamma}_i$ for $i=1,\cdots,m$, and $k$, run \textbf{Algorithm \ref{Alg:pfa}}, and obtain the dependence-adjusted correlation matrix $\tilde{\bm{\Sigma}}$ and the dependence-adjusted statistics
		$\widehat{Z}_i$ for $i=1,\cdots,m$.
		\item Randomly generate $N$ independent sets of statistics: $(Z_{r,1}^0, \ldots, Z_{r,m}^0)^T \sim N((0,\cdots,0)^T, \tilde{\bm{\Sigma}})$ for $r=1,\cdots,N$.
		\item Input $Z_{r,i}^0$ for $r=1,\cdots,N$ and $i=1,\cdots,m$, and $\alpha/2$, run \textbf{Algorithm \ref{Alg:boundseq}}, and obtain the bounding sequences $c_{m,0.5}$ and $c_{m,1}$.
		\item Input $\widehat{Z}_i$ for $i=1,\cdots,m$, and $c_{m,0.5}$ and $c_{m,1}$, run \textbf{Algorithm \ref{Alg:signprop}}, and output the signal proportion estimate $\hat{\pi}_A$.
		\item \textbf{Output}: The dependence-adjusted lower bound approximated estimator $\hat{\pi}_A$.
	\end{enumerate}
\end{algorithm}

In Algorithm $\ref{Alg:pfa}$, the dependence-adjusted statistics are computed using the PFA approach, detailed in Section \ref{sec:PFA}. Algorithms \ref{Alg:boundseq} and \ref{Alg:signprop} together provide the lower bound estimator  discussed in Section \ref{sec:family}. 
The dependence-adjusted estimator using the lower bound approach, $\hat{\pi}_L(\widehat{\mathbf{z}})$, is presented in Algorithm \ref{Alg:lowerbound}. Similarly, the dependence-adjusted estimator using the approximation approach, $\hat{\pi}_A(\widehat{\mathbf{z}})$, is shown in Algorithm \ref{Alg:approxlowerbound}.

Compared to Algorithm \ref{Alg:lowerbound}, Algorithm \ref{Alg:approxlowerbound} is computationally fast because it requires estimating $\mathbf{w}$ using robust $L_1$-regression only once, rather than $N+1$ times. Here, $N$ represents the number of replication sets of statistics generated from the joint null distribution, which is typically set to a large, predetermined value. In our simulation studies, we set $N=1000$.

The number of principal factors $k$ should be selected with care. A larger $k$ increases the $a_i$ values, which can have two opposing effects on estimation. On the one hand, larger $a_i$ values amplify the signals, as each $\mu_i$ is scaled by $a_i$ for $i \in I_1$, making them easier to detect. On the other hand, since the principal factors need to be estimated, the discrepancy between the true $Z_i$ and the estimated $\widehat{Z}_i$, given by $a_i \mathbf{b}_i^T(\widehat{\mathbf{W}} - \mathbf{W})$ for $i = 1, \ldots, m$, will also be larger with increasing $a_i$ values.
\cite{fan2012estimating} suggested choosing the smallest $k$ such that the random errors $(K_1, \ldots, K_m)^T$ are weakly dependent, using the criterion
\[
m^{-1}\|\mathbf{A}\|_F=m^{-1}\sqrt{\lambda_{k+1}^2+\cdots+\lambda_m^2}<\epsilon,
\]
where $\epsilon$ is a predetermined small value such as 0.01 or 0.05. In our simulation studies, we set $\epsilon = 1/\sqrt{m}$ to reflect the effect of dimensionality. Determining the optimal number of principal factors remains an open problem that requires further research.

%and we will show examples of how the choice of $k$ may affect our performance in the simulation studies.

\section{Simulation} \label{sec:simulation}

\subsection{Simulation setup}

In our simulation studies, test statistics are generated as $(X_1,\cdots,X_m)^T\sim N(\bm{\mu},{\bm{\Sigma}})$. Given a signal proportion $\pi$, the set of indices for signal variables, $I_1$, is randomly sampled from $\{1,\cdots,m\}$ with size $|I_1| = \pi m$. We set $\pi = 0.02$ or 0.1, representing relative sparse or dense signal variables, respectively. 

For the correlation matrix $\bm{\Sigma}$, we consider six different dependence structures drawn from real high-dimensional datasets and popular models in the literature.

\begin{itemize}
	\item \textbf{(i) Gene Network: } $\bm{\Sigma}$ is the sample correlation matrix from a real dataset containing the expression levels of $4088$ genes from 71 individuals in a riboflavin production study. The dataset is available at  \url{https://www.annualreviews.org/content/journals/10.1146/annurev-statistics-022513-115545}. Here, $m=4088$. 
	\item \textbf{(ii) SNP LD:} $\bm{\Sigma}$ is the sample correlation matrix from a real dataset containing genotype data of chromosome 21 from 90 individuals in the International HapMap project. Routine LD pruning results in 8657 SNPs, so $m=8657$. This dataset has also been analyzed in \cite{jeng2023estimating}.
	\item \textbf{(iii) Factor Model: }$\bm{\Sigma}$ is the correlation matrix of $\mathbf{V}$, where $\mathbf{V}=\tau\mathbf{h}\mathbf{h}^T+\mathbf{I}_m$ with $\tau\in(0,1)$ and $\mathbf{h}\sim N(\mathbf{0},\mathbf{I}_m)$. We set $m=2000$ and $\tau=0.5$.
	\item \textbf{(iv) Block Model:} $\bm{\Sigma}$ is an $m\times m$ correlation matrix with 20 diagonal blocks, where block sizes are randomly sampled from 10 to 100. Within each block, the off-diagonal entries have a value of $\rho$. We set $m=2,000$ and $\rho=0.5$.
	\item \textbf{(v) Autocorrelation}: $\Sigma_{ij} = 0.2^{|i-j|}$, $i=1, \cdots, m$ and $j = 1, \cdots, m$, with $m=2000$. 
	\item \textbf{(vi) Small Blocks}: $\bm{\Sigma}$ consists of small diagonal blocks, each of size 4, with off-diagonal entries having a value of $\rho$. We set $m=2,000$ and $\rho=0.5$.
\end{itemize}

\subsection{Analytical Evaluation of Dependence Adjustment Benefits} \label{sec:sim_evaluation}

We begin by assessing the potential advantage of the proposed dependence-adjusted estimator over existing non-adjusted approaches by comparing the stringency of condition (\ref{cond:ah}) for $\hat{\pi}^*(\mathbf{z})$ with that of condition (\ref{cond:h}) for $\hat{\pi}^*(\mathbf{x})$  across the six different dependence structures. 
This assessment is conducted prior to applying any of the estimation approaches, 
serving as a guideline for determining whether the proposed dependence-adjusted method should be used in specific scenarios.
Such a guideline is highly desirable due to its practical relevance.

To facilitate the comparison, denote the theoretical mean boundary for $\hat{\pi}^*(\mathbf{x})$ as
\begin{equation}\label{def:mu_*}
\mu_*=\sqrt{2h_*\log m},\quad \text{where } h_*= [\gamma-(\eta-\gamma)_+]_+
\end{equation}
and the theoretical mean boundary for $\hat{\pi}^*(\mathbf{z})$ as
\begin{equation}\label{def:mu_*}
\tilde\mu_* = \sqrt{2\tilde h_*\log m},\quad \text{where } \tilde h_*= [\gamma-(\tilde \eta-\gamma)_+]_+ / a^2_{\min}. 
\end{equation}
For the six dependence structures, we calculate the values of $\eta$, $\tilde{\eta}$, $a_{\min}$, $a_{\max}$, and $k$ from the given $\bm{\Sigma}$, and obtain the value of $\gamma$ from the given $\pi$. We then compare the values of $\mu_*$ and $\tilde{\mu}_*$. The results are presented in Table \ref{table:evaluation}.

\begin{table}[!h]
	\centering
	\caption{Parameter values and the theoretical mean boundaries before and after dependence adjustment. The key indicators for potential power gain are highlighted.} \label{table:evaluation}
	\begin{tabular}{llllllllll}
		\hline
		Dependence & \multicolumn{2}{l}{Signal Sparsity} & \multicolumn{2}{l}{Before Adjustment} & \multicolumn{5}{l}{After Adjustment}\\
		Structure & $\pi$ & $\gamma$ & $\eta$ & $\mu_*$ & $\tilde{\eta}$ & $a_{\min}$ & $a_{\max}$ & $k$ & $\tilde{\mu}_*$ \\ 
		\hline
		(i) Gene Network & 0.02 & 0.47 & {\bf 0.13} & 2.80 & {\bf 0.24} & {\bf 1.32} & 8.58 & 28 & 2.12 \\ 
		& 0.10 & 0.28 & {\bf 0.13} & 2.15 & {\bf 0.24} & {\bf 1.32} & 8.58 & 28 & 1.63 \\ 
		\hline
		(ii) SNP LD & 0.02 & 0.43 & 0.27 & 2.80 & 0.08 & {\bf 2.78} & 214 & 86 & 1.01 \\ 
		& 0.10 & 0.25 & 0.27 & 2.08 & 0.08 & {\bf 2.78} & 214 & 86 & 0.77 \\ 
		\hline 
		(iii) Factor Model & 0.02 & 0.51 & {\bf 0.22} & 2.80 & {\bf 0.92} & 1.00 & 2.87 & 1 & 1.29 \\ 
		& 0.10 & 0.30 & {\bf 0.22} & 2.15 & {\bf 0.92} & 1.00 & 2.87 & 1 & 0.00 \\ 
		\hline 
		(iv) Block Model & 0.02 & 0.51 & {\bf 0.60} & 2.56 & {\bf 0.91} & 1.00 & 1.43 & 16 & 1.34 \\ 
		& 0.10 & 0.30 & {\bf 0.60} & 0.31 & {\bf 0.91} & 1.00 & 1.43 & 16 & 0.00 \\ 
		\hline 
		(v) Autocorrelation & 0.02 & 0.51 & 0.95 & 1.12 & 0.78 & 1.00 & 1.04 & 75 & 1.93 \\ 
		& 0.10 & 0.30 & 0.95 & 0.00 & 0.78 & 1.00 & 1.04 & 75 & 0.00 \\ 
		\hline
		(vi) Small Blocks & 0.02 & 0.51 & 0.88 & 1.51 & 0.89 & 1.00 & 1.63 & 241 & 1.44 \\ 
		& 0.10 & 0.30 & 0.88 & 0.00 & 0.89 & 1.00 & 1.63 & 241 & 0.00 \\ 
		\hline
	\end{tabular}
\end{table}

As we can see, the overall covariance strength, presented through $\eta$, is strongest in the (i) Gene Network case, relatively strong in the (ii) SNP LD and (iii) Factor Model cases, moderately weak in the (iv) Block Model case, and weakest in the (v) Autocorrelation and (vi) Small Blocks cases. 

After dependence adjustment, $a_{\min} > 1$ in cases (i) and (ii), indicating increased signal intensity. The overall covariance strength is much reduced, with $\tilde{\eta} > \eta$, in cases (i), (iii), and (iv), leading to decreased variability in the data. Consequently, the theoretical mean boundaries are much lowered ($\tilde{\mu}_* < \mu_*$) across cases (i) through (iv), suggesting substantial power gains from the dependence-adjusted approach in these scenarios. Even without knowing the exact value of $\gamma$ in real applications, the decrease in $\tilde{\mu}_*$ is implied by observing an increase in $\tilde\eta$.

On the other hand, in the (v) Autocorrelation and (vi) Small Blocks cases, $\eta$ values are close to 1, indicating very weak overall dependence. After dependence adjustment, $\tilde\eta$ does not significantly increase and may even decrease, with the adjusted mean boundary $\tilde \mu_*$ showing little change or even increasing. These indicators do not suggest a substantial power gain for the dependence-adjusted estimator in these scenarios. It is important to note that since $\mu_*$ and $\tilde{\mu}_*$ are defined on the scale of $\sqrt{\log m}$, zero values of $\mu_*$ and $\tilde{\mu}_*$ in Table \ref{table:evaluation} correspond to values on the order of $o(\sqrt{\log m})$.

\subsection{Method comparison}

We compare our dependence-adjusted proportion estimators using both the lower bound approach ($\hat\pi_L(\widehat{\mathbf{z}})$) and the approximation approach ($\hat\pi_A(\widehat{\mathbf{z}})$) against the existing adaptive estimator ($\hat\pi({\mathbf{x}})$). As discussed in Section \ref{sec:L_approach} - \ref{sec:A_approach}, $\hat\pi_L(\widehat{\mathbf{z}})$ and $\hat\pi({\mathbf{x}})$ share the same lower bound property, while $\hat\pi_A(\widehat{\mathbf{z}})$ only asymptotically holds the lower bound property. For all three methods, the lower bound level is set at $\alpha=0.1$.

Other existing estimators, including both lower bound and non-lower bound types, have been extensively evaluated alongside $\hat\pi({\mathbf{x}})$ in \cite{jeng2023estimating}, where $\hat\pi({\mathbf{x}})$ showed advantages across various dependence structures and signal sparsity levels in finite samples. To avoid redundancy, we include one widely used non-lower bound estimator in this comparison:
\[
\hat{\pi}_{\lambda}=1-\min\left\{1,\quad\dfrac{\sum\limits_{i=1}^{m}\mathbb{I}\{p_i>\lambda\}}{m(1-\lambda)}\right\},
\]
where $p_1,\cdots,p_m$ are the $p$-values of $x_1, \cdots, x_m$, and $\lambda$ is a tuning parameter, set to the lower 40\% sample quantile of $p_1,\cdots,p_m$ as used in the PFA software package. This estimator, originally developed in \cite{storey2002direct} and \cite{storey2004strong}, is widely used in the literature but does not possess the lower bound property.

Two signal sparsity levels are considered: $\pi = 0.02$ and $\pi = 0.1$, representing moderately dense and moderately sparse cases, respectively. While traditional estimators like $\hat\pi_{\lambda}$ typically focus on relatively dense signal proportions, studies on the lower bound family have explored a broader range of signal sparsity ($\gamma \in (0, 1)$). 
The finite-sample performances of $\hat\pi_L(\widehat{\mathbf{z}})$, $\hat\pi_A(\widehat{\mathbf{z}})$, $\hat\pi({\mathbf{x}})$, and $\hat \pi_{\lambda}$ are evaluated with $\mu_i=1, 2$, or $3$ for all $i\in I_1$. These non-zero mean values are chosen based on the range of $\mu_*$ and $\tilde \mu_*$ shown in Table \ref{table:evaluation}. 

The estimation results are presented as the ratio $\hat\pi / \pi$, where $\hat\pi$ represents the estimated value from each method. A ratio closer to 1 indicates more accurate estimation. Additionally, for a lower bound estimator, the ratio is generally expected to be less than 1.

\subsection{Simulation results}

\begin{figure}[!h]
	\centering
	\includegraphics[width=1\textwidth]{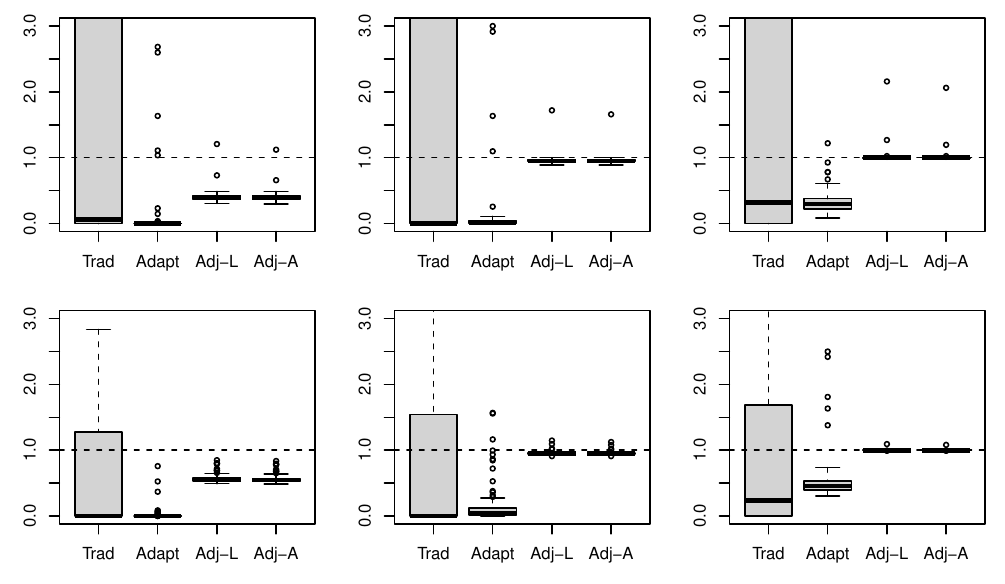}
	\caption{Boxplots of $\hat{\pi}/\pi$ for the (i) Gene Network case across four methods: Traditional ($\hat\pi_{\lambda}$), Adaptive ($\hat \pi(\mathbf{x})$), Adjust-L ($\hat\pi_L(\widehat{\mathbf{z}})$), and Adjust-A ($\hat\pi_A(\widehat{\mathbf{z}})$). The top row has $\pi=0.02$, while the bottom row has $\pi=0.1$. The three columns correspond to signal intensity levels 1, 2, and 3, respectively.}\label{fig:box_gene}
\end{figure}

We generate boxplots of the $\hat{\pi}/\pi$ results for different methods across 100 replications. %, focusing on the dependence structures in cases (i) - (iv), which exhibit relatively strong overall dependence.
Figure \ref{fig:box_gene} presents the results for the (i) Gene Network case. The top row corresponds to $\pi=0.02$, while the bottom row corresponds to $\pi=0.1$. The three columns, from left to right, represent signal intensity levels 1, 2, and 3, respectively. In this example, the traditional estimator $\hat \pi_{\lambda}$ and the adaptive estimator $\hat \pi(\mathbf{x})$ have comparable median values of $\hat{\pi}/\pi$, but $\hat\pi_\lambda$ shows higher volatility and is more prone to overestimation, particularly in the sparser scenario with $\pi=0.02$. In contrast, the two dependence-adjusted approaches, $\hat\pi_L(\widehat{\mathbf{z}})$ and $\hat\pi_A(\widehat{\mathbf{z}})$, significantly outperform both $\hat \pi_{\lambda}$ and $\hat \pi(\mathbf{x})$, with similar performance between the two adjusted methods.

\begin{figure}[!htbp]
	\centering
	\includegraphics[width=1\textwidth]{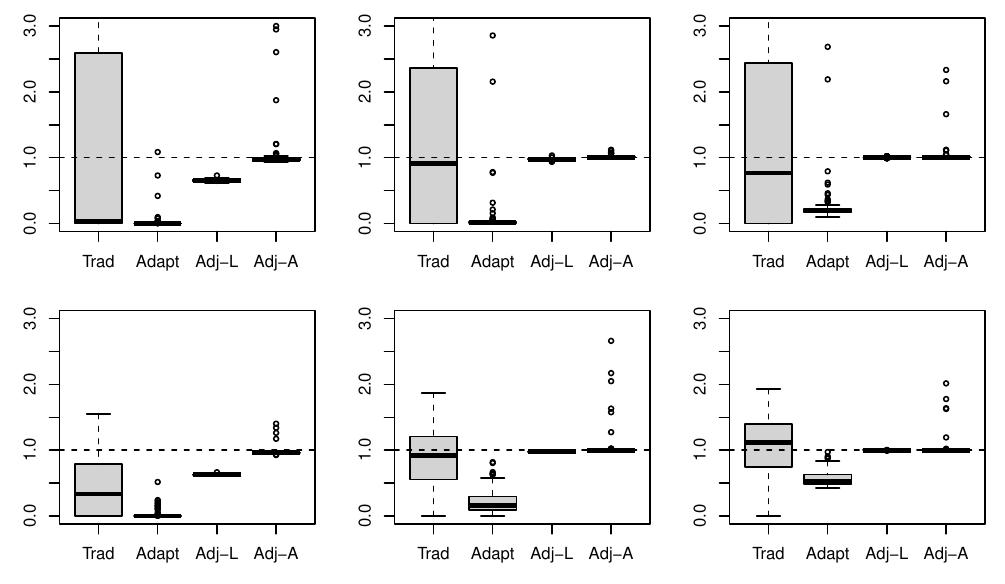}
	\caption{Boxplots of $\hat{\pi}/\pi$ for the (ii) SNP LD case across four methods. The notations and row/column assignments follow those in Figure \ref{fig:box_gene}.}\label{fig:box_snp}
\end{figure}

\begin{figure}[!htbp]
	\centering
	\includegraphics[width=1\textwidth]{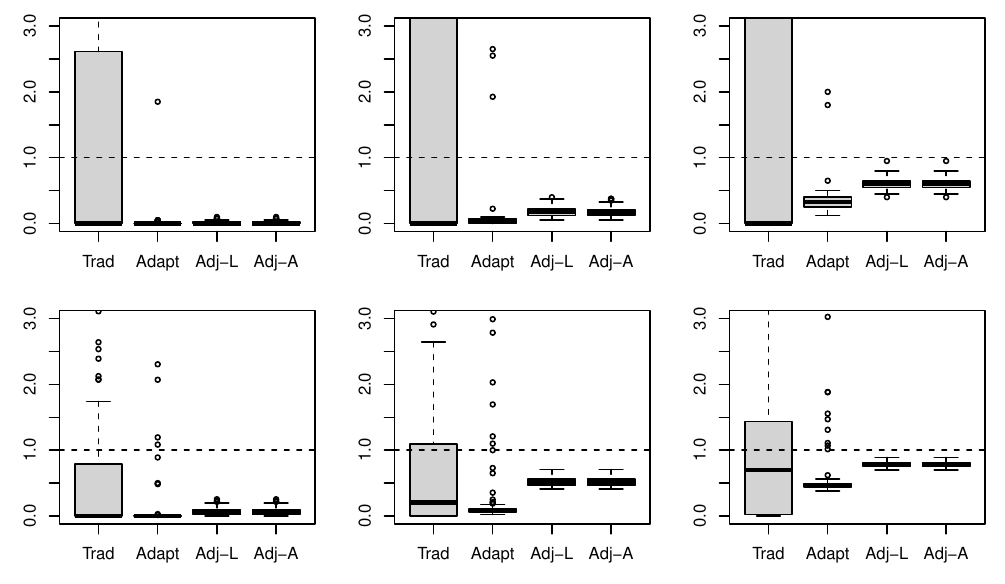}
	\caption{Boxplots of $\hat{\pi}/\pi$ for the (iii) Factor Model case across four methods. The notations and row/column assignments follow those in Figure \ref{fig:box_gene}.}\label{fig:box_factor}
\end{figure}

\begin{figure}[!htbp]
	\centering
	\includegraphics[width=1\textwidth]{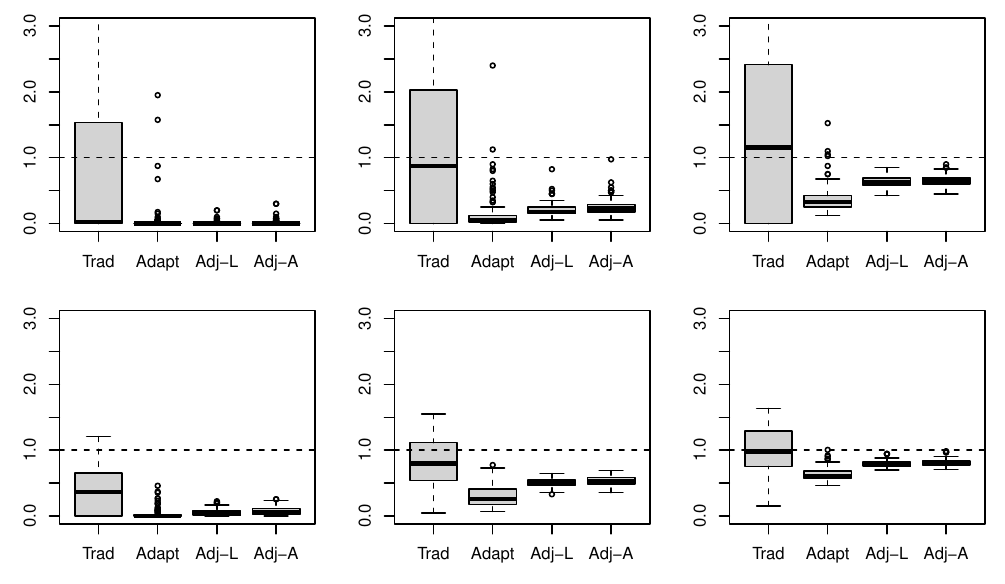}
	\caption{Boxplots of $\hat{\pi}/\pi$ for the (iv) Block Model case across four methods. The notations and row/column assignments follow those in Figure \ref{fig:box_gene}.}\label{fig:box_block}
\end{figure}

\begin{figure}[!htbp]
	\centering
	\includegraphics[width=1\textwidth]{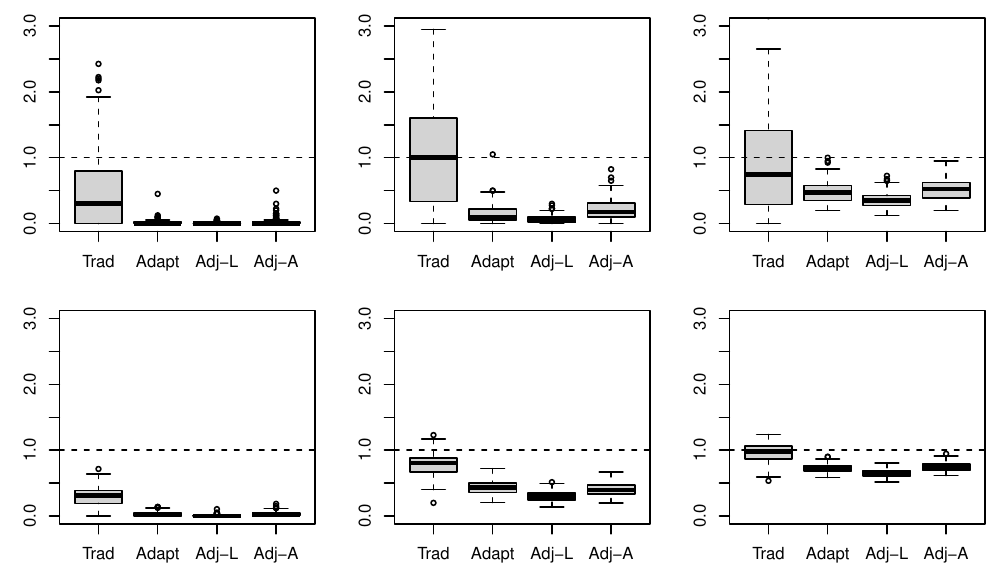}
	\caption{Boxplots of $\hat{\pi}/\pi$ for the (v) Autocorrelation case across four methods. The notations and row/column assignments follow those in Figure \ref{fig:box_gene}.}\label{fig:box_auto}
\end{figure}

\begin{figure}[!htbp]
	\centering
	\includegraphics[width=1\textwidth]{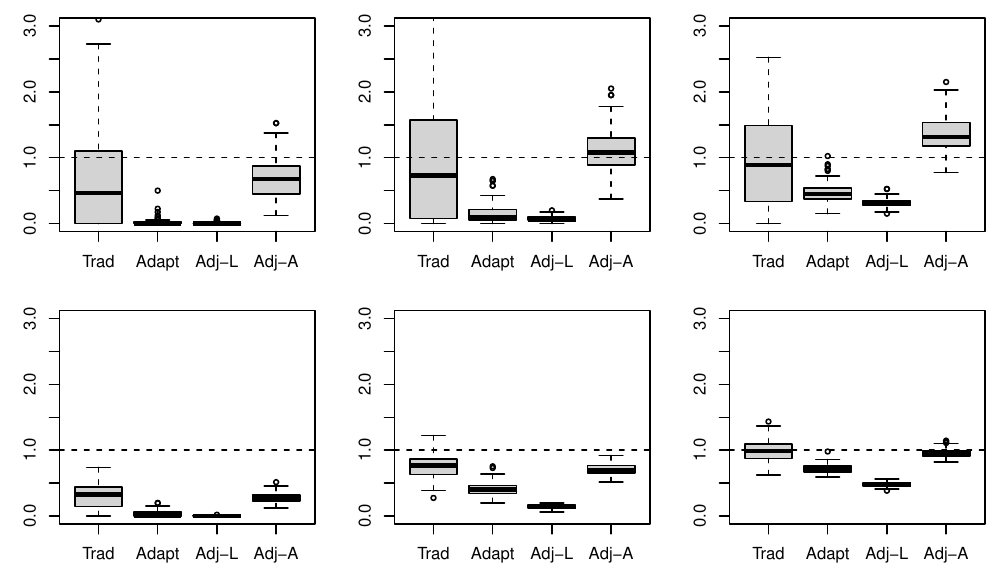}
	\caption{Boxplots of $\hat{\pi}/\pi$ for the (vi) Small Blocks case across four methods. The notations and row/column assignments follow those in Figure \ref{fig:box_gene}.}\label{fig:box_smallblock}
\end{figure}

Similar results are observed for the (ii) SNP LD case in Figure \ref{fig:box_snp}, where the advantages of $\hat\pi_L(\widehat{\mathbf{z}})$ and $\hat\pi_A(\widehat{\mathbf{z}})$ over $\hat \pi_{\lambda}$ and $\hat \pi(\mathbf{x})$ are clear. In this example, the difference between $\hat\pi_L(\widehat{\mathbf{z}})$ and $\hat\pi_A(\widehat{\mathbf{z}})$ is more pronounced in the first column, where the signal intensity is lowest. While $\hat\pi_A(\widehat{\mathbf{z}})$ exhibits better median values, it has more instances of overestimation. This can be explained by the high $a_{\max}$ value in this case, as shown in Table \ref{table:evaluation}, which challenges condition (C1'') in Theorem \ref{thm:pi_hat_A}, thereby impacting the asymptotic lower bound property of $\hat\pi_A(\widehat{\mathbf{z}})$.

Figures \ref{fig:box_factor} and \ref{fig:box_block} display the results for the (iii) Factor Model and (iv) Block Model cases, respectively. 
In these examples, $\hat\pi_\lambda$ continues to exhibit high volatility and a significant portion of overestimated results, while the other methods -- $\hat \pi(\mathbf{x})$, $\hat\pi_L(\widehat{\mathbf{z}})$, and $\hat\pi_A(\widehat{\mathbf{z}})$ -- demonstrate consistent lower bound property. Among these three, the dependence-adjusted estimators $\hat\pi_L(\widehat{\mathbf{z}})$ and $\hat\pi_A(\widehat{\mathbf{z}})$ continue to outperform $\hat \pi(\mathbf{x})$.

Figures {\ref{fig:box_auto}} and \ref{fig:box_smallblock} correspond to the (v) Autocorrelation and (vi) Small Blocks cases, respectively, where power gains for the dependence-adjusted approaches are not anticipated, as discussed in Section \ref{sec:sim_evaluation}. This expectation is confirmed by the results in Figures {\ref{fig:box_auto}} and \ref{fig:box_smallblock}. 
In Figures {\ref{fig:box_auto}}, 
$\hat\pi_L(\widehat{\mathbf{z}})$ performs slightly worse than $\hat\pi(\widehat{\mathbf{x}})$, while $\hat\pi_A(\widehat{\mathbf{z}})$ performs comparably to $\hat\pi_L(\widehat{\mathbf{z}})$. Moreover, in Figure \ref{fig:box_smallblock}, $\hat\pi_A(\widehat{\mathbf{z}})$ has an excessive amount of overestimated results. This can be attributed to the high $k$ value in this case, as shown in Table \ref{table:evaluation}, which challenges condition (C1'') in Theorem \ref{thm:pi_hat_A}, thereby affecting the asymptotic lower bound property of $\hat\pi_A(\widehat{\mathbf{z}})$.

\section{Conclusion} \label{sec:conclusion}

This study makes a significant contribution to the field of signal proportion estimation by introducing an innovative estimator that effectively accounts for the dependencies among variables. Unlike traditional methods that often assume independence, our estimator leverages the underlying dependence structures, making it more applicable to real-world data where such assumptions are frequently violated. 

Our approach features a lower bound strategy to prevent the overestimation of the true signal proportion, which is crucial in fields where incorrect signal identification can have significant consequences. This ensures that the estimates remain cautious and reliable. Additionally, by integrating a principal factor approximation procedure, our estimator enhances accuracy across various dependence scenarios, including those with high correlations among signal variables. The estimator's versatility in handling different sparsity levels and dependence structures makes it a robust tool for a wide range of scientific research.

Theoretical insights presented in this paper deepen the understanding of how signal sparsity, signal intensity, and covariance dependence interact.
By comparing the conditions for estimation consistency before and after dependence adjustment, we demonstrate the benefits of incorporating dependence information in various settings. This theoretical framework not only validates the effectiveness of the new estimator but also guides its practical application, ensuring its reliable use in diverse scenarios. The extensive simulation studies and comparisons with existing methods further confirm the superior performance of our estimator, particularly in scenarios where dependence structures are strong and complex.

Overall, the novel signal proportion estimator introduced in this work represents a meaningful advancement in statistical methodology. We anticipate that this estimator will significantly benefit the scientific community by enabling more accurate and reliable analyses, ultimately leading to more informed and impactful discoveries.

A natural extension is to relax the assumption that the correlation matrix $\bm{\Sigma}$ is known. In many occasions, it is estimated either from the data used to form the test statistics or from an auxiliary base dataset. Studying how estimation error affects signal proportion estimation would broaden the practical applicability of our method.

	\section{Appendix: Detailed proofs for the theoretical results}\label{sec:appendix}%% if no title is needed, leave empty \section*{}.
	
	In this section, we provide detailed proofs for the theoretical results presented throughout the paper. The notation ``$C$" is used to represent a generic constant, the value of which may vary depending on the context.
	
	For simplicity of notation, the lowercase letters 
	$x_i$, $x_i^0$,  $z_i$, $z^0_i$, $\hat z_i$, and $\hat z^0_i$ are used to denote random variables, provided there is no ambiguity in meaning.

	\subsection{Preliminary lemmas}
	
	\begin{lemma} \label{lemma:lowerbd_x}
		Consider model (\ref{def:model}) and the family of lower bound estimators in (\ref{def:pi_hat_theta}) with bounding sequences $c_{m,\theta}(\alpha; F_{{\mathbf{x}}})$ satisfying the following properties: \\
		(i) $mc_{m,\theta}(\alpha; F_{\mathbf{x}})>(m-s)c_{m-s, \theta}(\alpha; F_{\mathbf{x}})$ and \\
		(ii) $\text{P}(V_{m,\theta}(\mathbf{x}^0)>c_{m,\theta}(\alpha; F_{\mathbf{x}}))<\alpha$ for all $m$. \\
		Then, for any given $\theta \in [0, 1]$, we have 
		\[
		\text{P}(\hat \pi_\theta ({\mathbf{x}})<\pi)\geq1-\alpha.
		\]
	\end{lemma}
	
	{\bf Proof:} Note that
	\begin{eqnarray}\label{lowerbound.1}
	\hat{\pi}_\theta(\mathbf{x})-\pi&=&\sup\limits_{t>0}\left\{\dfrac{m^{-1}\sum\limits_{i=1}^{m}\mathbb{I}\{x_i> t\}-\bar{\Phi}(t)-c_{m,\theta}(\alpha; F_{\mathbf{x}})[\bar{\Phi}(t)]^\theta}{1-\bar{\Phi}(t)}-\pi\right\}\nonumber\\
	&=&\sup\limits_{t>0}\dfrac{m^{-1}\sum\limits_{i\in I_0}^{}\mathbb{I}\{x_i> t\}+m^{-1}\sum\limits_{i\in I_1}^{}\mathbb{I}\{x_i> t\}-\pi-(1-\pi)\bar{\Phi}(t)-c_{m,\theta}(\alpha; F_{\mathbf{x}})[\bar{\Phi}(t)]^\theta}{1-\bar{\Phi}(t)}\nonumber\\
	&\leq&\sup\limits_{t>0}\dfrac{(1-\pi)(m-s)^{-1}\sum\limits_{i\in I_0}^{}\mathbb{I}\{x_i> t\}+m^{-1}s-\pi-(1-\pi)\bar{\Phi}(t)-c_{m,\theta}(\alpha; F_{\mathbf{x}})[\bar{\Phi}(t)]^\theta}{1-\bar{\Phi}(t)}\nonumber\\
	&=&\sup\limits_{t>0}\dfrac{(1-\pi)\left\{(m-s)^{-1}\sum\limits_{i\in I_0}^{}\mathbb{I}\{x_i> t\}-\bar{\Phi}(t)\right\}-c_{m,\theta}(\alpha; F_{\mathbf{x}})[\bar{\Phi}(t)]^\theta}{1-\bar{\Phi}(t)}\nonumber\\
	&\leq&\sup\limits_{t>0}\dfrac{(1-\pi)\left\{(m-s)^{-1}\sum\limits_{i\in I_0}^{}\mathbb{I}\{x_i> t\}-\bar{\Phi}(t)-c_{m-s,\theta}(\alpha; F_{\mathbf{x}})[\bar{\Phi}(t)]^\theta\right\}}{1-\bar{\Phi}(t)},\nonumber\\
	\end{eqnarray}
	where the last inequality comes from property (i). This implies that
	\begin{eqnarray}\label{lowerbound.2}
	&&\text{P}(\hat{\pi}_\theta(\mathbf{x})>\pi)\nonumber\\
	&\leq&\text{P}\left(\sup\limits_{t>0}\dfrac{(1-\pi)\left\{(m-s)^{-1}\sum\limits_{i\in I_0}^{}\mathbb{I}\{x_i> t\}-\bar{\Phi}(t)-c_{m-s,\theta}(\alpha; F_{\mathbf{x}})[\bar{\Phi}(t)]^\theta\right\}}{1-\bar{\Phi}(t)}>0\right)\nonumber\\
	&\leq&\text{P}\left(\sup\limits_{t>0}\left\{\dfrac{(m-s)^{-1}\sum\limits_{i\in I_0}^{}\mathbb{I}\{x_i> t\}-\bar{\Phi}(t)}{[\bar{\Phi}(t)]^\theta}-c_{m-s,\theta}(\alpha; F_{\mathbf{x}})\right\}>0\right)\nonumber\\
	&\leq&\text{P}(V_{m-s,\theta}(\mathbf{x}^0)>c_{m-s,\theta}(\alpha; F_{\mathbf{x}}))<\alpha,\nonumber\\
	\end{eqnarray}
	where the last step comes from property (ii). Thus the conclusion follows.
	
	\begin{lemma} \label{lemma:lowerbd_z}
		Consider model (\ref{def:model}), the dependence adjusted statistics in (\ref{def:Z}), and the estimator in (\ref{def:pi_z}) with bounding sequences $c_{m,\theta}(\alpha/2; F_{{\mathbf{z}}})$ satisfying the following properties: \\
		(i) $mc_{m,\theta}(\alpha/2; F_{\mathbf{z}})>(m-s)c_{m-s, \theta}(\alpha/2; F_{\mathbf{z}})$ and \\
		(ii) $\text{P}(V_{m,\theta}(\mathbf{z}^0)>c_{m,\theta}(\alpha/2; F_{\mathbf{z}}))<\alpha/2$ for all $m$. \\
		Then, we have 
		\[
		\text{P}(\hat \pi ({\mathbf{z}})<\pi)\geq1-\alpha.
		\]
	\end{lemma}
	
	{\bf Proof:} Here for any given $\theta\in[0,1]$, we have
	\begin{equation} \label{def:pi_hat_theta_z}
	\hat{\pi}_\theta(\mathbf{z})=\sup\limits_{t>0}\dfrac{m^{-1}\sum\limits_{i=1}^{m}\mathbb{I}\{z_i> t\}-\bar{\Phi}(t)-c_{m,\theta}(\alpha/2; F_{\mathbf{z}}) [\bar{\Phi}(t)]^\theta}{1-\bar{\Phi}(t)}
	\end{equation}
	and 
	\begin{equation} \label{def:V_Z0}
	V_{m,\theta}(\mathbf{z}^0)=\sup\limits_{t>0}\dfrac{\left|m^{-1}\sum\limits_{i=1}^{m}\mathbb{I}\{z_i^0> t\}-\bar{\Phi}(t)\right|}{[\bar{\Phi}(t)]^\theta}. 
	\end{equation}
	With a similar argument to the proof of Lemma \ref{lemma:lowerbd_x}, replacing $x_i$ with $z_i$, $c_{m,\theta}(\alpha; F_{\mathbf{x}})$ with $c_{m,\theta}(\alpha/2; F_{\mathbf{z}})$, and $V_{m,\theta}(\mathbf{x}^0)$ with $V_{m,\theta}(\mathbf{z}^0)$, we can show that $\text{P}(\hat{\pi}_\theta(\mathbf{z})<\pi)\geq1-\alpha/2$.
	
	Taking $\theta=0.5$, we have $\text{P}(\hat{\pi}_{0.5}(\mathbf{z})<\pi)\geq 1-\alpha/2$. Taking $\theta=1$, we have $\text{P}(\hat{\pi}_{1}(\mathbf{z})<\pi)\geq 1-\alpha/2$. Then by (\ref{def:pi_z}), we have $\text{P}(\hat{\pi}(\mathbf{z}) < \pi) \ge (1-\alpha/2)^2\ge 1-\alpha$.

	\begin{lemma} \label{lemma:lowerbd_z_hat}
		Consider model (\ref{def:model}), the dependence adjusted statistics with estimated principal factors in (\ref{def:z_hat}), and the estimator in (\ref{def:pi_L}) with bounding sequences $c_{m,\theta}(\alpha/2; F_{{\widehat{\mathbf{z}}}})$ satisfying the following properties: \\
		(i) $mc_{m,\theta}(\alpha/2; F_{{\widehat{\mathbf{z}}}})>(m-s)c_{m-s, \theta}(\alpha/2; F_{{\widehat{\mathbf{z}}}})$ and \\
		(ii) $\text{P}(V_{m,\theta}({\widehat{\mathbf{z}}}^0)>c_{m,\theta}(\alpha/2; F_{{\widehat{\mathbf{z}}}}))<\alpha/2$ for all $m$. \\
		Then, we have 
		\[
		\text{P}(\hat \pi_L ({\widehat{\mathbf{z}}})<\pi)\geq1-\alpha.
		\]
	\end{lemma}
	
	The proof of Lemma \ref{lemma:lowerbd_z_hat} closely follows that of Lemma \ref{lemma:lowerbd_z} and is therefore omitted.
	
	\subsection{Proof of Proposition \ref{thm:pi_hat_x}} \label{sec:proof_prop3.1}
	
	Recall that the discretized $\hat \pi^*(\mathbf{x})$ is defined as 
	\begin{equation} \label{def:pi_hat*_z_tilde}
	\hat{\pi}^*(\mathbf{x})=\max\{\hat{\pi}^*_{0.5}(\mathbf{x}),~\hat{\pi}^*_1(\mathbf{x})\},
	\end{equation} 
	where 
	\begin{equation} \label{def:hat_pi*_theta_tilde}
	\hat{\pi}^*_\theta(\mathbf{x})=\max\limits_{t\in\mathbb{T}}\dfrac{m^{-1}\sum\limits_{i=1}^{m}\mathbb{I}\{x_i> t\}-\bar{\Phi}(t)-c^*_{m,\theta}(\alpha; F_{\mathbf{x}})[\bar{\Phi}(t)]^\theta}{1-\bar{\Phi}(t)}
	\end{equation}
	for a given $\theta \in [0,1]$, with $\mathbb{T}=[\sqrt{\log \log m},\sqrt{5\log m}]\cap\mathbb{N}$. The bounding sequence $c^*_{m,\theta}(\alpha; F_{\mathbf{x}})$ is constructed such that $\text{P}(V^*_{m,\theta}(\mathbf{x}^0) > c^*_{m,\theta}(\alpha; F_{\mathbf{x}})) < \alpha$, where 
	\begin{equation}
	V^*_{m,\theta}(\mathbf{x}^0)=\max\limits_{t\in\mathbb{T}}\dfrac{\left|m^{-1}\sum\limits_{i=1}^{m}\mathbb{I}\{x_i^0> t\}-\bar{\Phi}(t)\right|}{[\bar{\Phi}(t)]^\theta}. \nonumber
	\end{equation}
	
	It is easy to see that Lemma \ref{lemma:lowerbd_x}, in its discretized form, implies the follows: 
	\begin{equation} \label{4.1}
	\text{P}(\hat{\pi}^*_1(\mathbf{x}) < \pi) \ge 1-\alpha 
	\end{equation}
	and 
	\begin{equation} \label{4.2}
	\text{P}(\hat{\pi}^*_{0.5}(\mathbf{x}) < \pi) \ge 1-\alpha.  
	\end{equation}
	
	Before arriving at the result in (\ref{eq:lowerbd_x}), we first provide the following lemmas that specify the bounding sequences  for $\theta=1$ and 0.5, respectively, and establish the upper bound results of $\hat{\pi}^*_1(\mathbf{x})$ and $\hat{\pi}^*_{0.5}(\mathbf{x})$. Their proofs are provided in Sections \ref{sec:proof_7.4}, \ref{sec:proof_7.5}, and \ref{sec:proof_7.6}. 
	Following this, we prove the result in (\ref{cond:h}) and conclude by demonstrating the result in (\ref{eq:lowerbd_x}).
	
	\begin{lemma} \label{lemma:boundingSeq}
		Consider model (\ref{def:model}) along with the calibrations in (\ref{def:gamma}) - (\ref{def:eta}). Let $c^*_{m,1}(\alpha; F_{{\mathbf{x}}}) = C_{1,\alpha}\sqrt{\log m}$ for some constant $C_{1,\alpha}$ and $c^*_{m,0.5}(\alpha; F_{{\mathbf{x}}}) = C_{0.5,\alpha}\sqrt{m^{-\eta}\log m}$ for some constant $C_{0.5,\alpha}$. Then both  $c^*_{m,1}(\alpha; F_{{\mathbf{x}}})$ and $c^*_{m,0.5}(\alpha; F_{{\mathbf{x}}})$ satisfy properties (i) and (ii) in  Lemma \ref{lemma:lowerbd_x} in its discritized form. 
	\end{lemma}
	
	\begin{lemma} \label{lemma:suf_nec_1}
		Under the same conditions as in Lemma \ref{lemma:boundingSeq}, with $\alpha = \alpha_m\to 0$, we have $\text{P}(\hat{\pi}^*_{1}({\mathbf{x}}) > (1-\epsilon)\pi)\to 1$ for any constant $\epsilon>0$ if and (almost) only if $h > \gamma$. 
	\end{lemma}  
	
	\begin{lemma} \label{lemma:suf_nec_05}
		Under the same conditions as in Lemma \ref{lemma:boundingSeq}, with $\alpha = \alpha_m\to 0$, we have
		$\text{P}(\hat{\pi}^*_{0.5}({\mathbf{x}}) > (1-\epsilon)\pi)\to 1$ for any constant $\epsilon>0$ if and (almost) only if
		$h>[2\gamma-\eta]_+$. 
	\end{lemma}    
	Combining Lemmas \ref{lemma:boundingSeq} - \ref{lemma:suf_nec_05}, and the definition of $\hat{\pi}^*({\mathbf{x}})$ in (\ref{def:pi_hat*_z_tilde}) gives 
	$\text{P}((1-\epsilon)\pi < \hat{\pi}^*({\mathbf{x}}) < \pi )\to 1$ for any constant $\epsilon>0$ if and almost only if
	$h > \min\{\gamma,~ [2\gamma-\eta]_+\}$, which is equivalent to the condition in (\ref{cond:h}). 
	
	Next, we show (\ref{eq:lowerbd_x}). Lemmas \ref{lemma:suf_nec_1} - \ref{lemma:suf_nec_05} also imply that 
	\[
	\text{P}(\hat{\pi}^*_{1}({\mathbf{x}}) \ge \hat{\pi}^*_{0.5}({\mathbf{x}})) \to 1 \quad \text{when} \quad \gamma \le [2\gamma-\eta]_+
	\] 
	and 
	\[
	\text{P}(\hat{\pi}^*_{1}({\mathbf{x}}) < \hat{\pi}^*_{0.5}({\mathbf{x}})) \to 1 \quad \text{when} \quad \gamma > [2\gamma-\eta]_+
	\] 
	These, combined with the definition of $\hat{\pi}^*({\mathbf{x}})$ in (\ref{def:pi_hat*_z_tilde}) and the results in (\ref{4.1})- (\ref{4.2}),  gives (\ref{eq:lowerbd_x}). 
	
	\subsection{Proof of Lemma \ref{lemma:boundingSeq}} \label{sec:proof_7.4}
	
	Firstly, we show that the bounding sequence $c_{m,1}^*(\alpha; F_{\mathbf{x}}) =  C_{1,\alpha}\sqrt{\log m}$ satisfies (i) and (ii). Property (i) is trivial. 
	
	As for property (ii), let
	\[
	B(t)=\dfrac{\left|m^{-1}\sum\limits_{i=1}^{m}\mathbb{I}\{x_i^0> t\}-\bar{\Phi}(t)\right|}{\bar{\Phi}(t)}
	\]
	such that $V^*_{m,1}(\mathbf{x}^0)=\max\limits_{t\in\mathbb{T}}B(t)$.  By Markov’s inequality,
	\begin{eqnarray} \label{theta1_z.3}
	&&\text{P}(V^*_{m,1}(\mathbf{x}^0)>c_{m,1}^*(\alpha; F_{\mathbf{x}}))\nonumber\\
	&\le & (c_{m,1}^*(\alpha; F_{\mathbf{x}}))^{-1}\text{E}\left[\max\limits_{t\in\mathbb{T}}B(t)\right]=(c_{m,1}^*(\alpha; F_{\mathbf{x}}))^{-1}\int_{0}^{\infty}\text{P}\left(\max\limits_{t\in\mathbb{T}}B(t)>x\right)dx \nonumber \\ 
	&\le & (c_{m,1}^*(\alpha; F_{\mathbf{x}}))^{-1}\int_{0}^{\infty}\sum_{t\in\mathbb{T}}\text{P}(B(t)>x)dx=(c_{m,1}^*(\alpha; F_{\mathbf{x}}))^{-1}\sum_{t\in\mathbb{T}}\text{E}[B(t)]  \nonumber \\
	&\le & C\sqrt{\log m}(c_{m,1}^*(\alpha; F_{\mathbf{x}}))^{-1}\max\limits_{t\in\mathbb{T}}\text{E}[B(t)], 
	\end{eqnarray}
	where the last inequality holds due to $\mathbb{T}=[\sqrt{\log\log m},\sqrt{5\log m}]\cap\mathbb{N}$.
	
	Since
	\[
	\text{E}[B(t)]\leq \dfrac{m^{-1}\sum\limits_{i=1}^{m}\text{P}(x_i^0> t)+\bar{\Phi}(t)}{\bar{\Phi}(t)}=\dfrac{m^{-1}\sum\limits_{i=1}^{m}\bar{\Phi}(t)+\bar{\Phi}(t)}{\bar{\Phi}(t)}=2,
	\]
	(\ref{theta1_z.3}) gives
	\[
	\text{P}(V^*_{m,1}(\mathbf{x}^0)>c_{m,1}^*(\alpha; F_{\mathbf{x}}))\leq C\sqrt{\log m}(c_{m,1}^*(\alpha; F_{\mathbf{x}}))^{-1}.
	\]
	Thus, there exists $c_{m,1}^*(\alpha; F_{\mathbf{x}})=C_{1,\alpha}\sqrt{\log m}$ for a large enough constant $C_{1,\alpha}$ such that $\text{P}(V^*_{m,1}(\mathbf{x}^0)>c_{m,1}^*(\alpha; F_{\mathbf{x}}))<\alpha$. 
	
	Next we show that the bounding sequence $c^*_{m,0.5}(\alpha; F_{\mathbf{x}}) = C_{0.5,\alpha}\sqrt{m^{-\eta}\log m}$ satisfies (i) and (ii). Property (i) is trivial. 
	
	As for property (ii), let
	\[
	A(t)=\dfrac{\left|m^{-1}\sum\limits_{i=1}^{m}\mathbb{I}\{x_i^0> t\}-\bar{\Phi}(t)\right|}{[\bar{\Phi}(t)]^{1/2}}
	\]
	such that $V_{m,0.5}^*(\mathbf{x}^0)=\max\limits_{t\in\mathbb{T}}A(t)$. By Markov’s inequality,
	\begin{eqnarray} \label{theta0.5_z.3}
	&&\text{P}(V_{m,0.5}^*(\mathbf{x}^0)>c^*_{m,0.5}(\alpha; F_{\mathbf{x}}))\nonumber\\
	& \le & (c^*_{m,0.5}(\alpha; F_{\mathbf{x}}))^{-2}\text{E}\left[\left\{\max\limits_{t\in\mathbb{T}}A(t)\right\}^2\right]=(c^*_{m,0.5}(\alpha; F_{\mathbf{x}}))^{-2}\text{E}\left[\max\limits_{t\in\mathbb{T}}\{A(t)\}^2\right]\nonumber\\
	&\le& (c^*_{m,0.5}(\alpha; F_{\mathbf{x}}))^{-2}\text{E}\left[\sum_{t\in\mathbb{T}}\{A(t)\}^2\right] = (c^*_{m,0.5}(\alpha; F_{\mathbf{x}}))^{-2}\sum_{t\in\mathbb{T}}\text{E}\left[\{A(t)\}^2\right] \nonumber \\ 
	& \le & C\sqrt{\log m}(c^*_{m,0.5}(\alpha; F_{\mathbf{x}}))^{-2}\max\limits_{t\in\mathbb{T}}\text{E}\left[\{A(t)\}^2\right], 
	\end{eqnarray}
	where the last inequality holds due to $\mathbb{T}=[\sqrt{\log\log m},\sqrt{5\log m}]\cap\mathbb{N}$. 
	
	Here we have for $t\in\mathbb{T}$,
	\begin{eqnarray} \label{theta0.5_z.4}
	&&\text{E}\left[\{A(t)\}^2\right]= [\bar{\Phi}(t)]^{-1}\text{E}\left[\left\{m^{-1}\sum\limits_{i=1}^{m}\mathbb{I}\{x_i^0> t\}-\bar{\Phi}(t)\right\}^2\right]\nonumber\\
	&=&[\bar{\Phi}(t)]^{-1}\text{Var}\left(m^{-1}\sum\limits_{i=1}^{m}\mathbb{I}\{x_i^0> t\}\right)= m^{-2}[\bar{\Phi}(t)]^{-1}\text{Var}\left(\sum\limits_{i=1}^{m}\mathbb{I}\{x_i^0> t\}\right),
	\end{eqnarray}
	where the following lemma is useful to show the order of $\text{E}[\{A(t)\}^2]$.
		\begin{lemma} \label{lemma:var_sum_w}
			For $ (x_1^0, \ldots, x_m^0)^T \sim F_{ \mathbf{x}} = N(\mathbf{0}, {\bm{\Sigma}})$, we have
			\begin{equation}\label{var_sum_w}
			\text{Var}\left(\sum\limits_{i=1}^{m}\mathbb{I}\{x_i^0> t\}\right)\leq Cm^{2-\eta}e^{-t^2/2}.
			\end{equation}
		\end{lemma}
		
		Combining (\ref{theta0.5_z.4}) and (\ref{var_sum_w}), we obtain
	\begin{equation} \label{theta0.5_z.8}
	\text{E}\left[\{A(t)\}^2\right]= m^{-2}[\bar{\Phi}(t)]^{-1}\text{Var}\left(\sum\limits_{i=1}^{m}\mathbb{I}\{x_i^0> t\}\right)\leq C m^{-\eta}[\bar{\Phi}(t)]^{-1}e^{-t^2/2}.
	\end{equation}
	Using Mill's ratio, we have
	\begin{equation} \label{theta0.5_z.9}
	[\bar{\Phi}(t)]^{-1}e^{-t^2/2} \leq C[\phi(t)]^{-1}e^{-t^2/2}\left(\dfrac{1}{t}-\dfrac{1}{t^3}\right)^{-1}\leq C(t+o(t))\leq C\sqrt{\log m},
	\end{equation}
	where the last step comes from $\mathbb{T}=[\sqrt{\log\log m},\sqrt{5\log m}]\cap\mathbb{N}$.
	Combining (\ref{theta0.5_z.3}), (\ref{theta0.5_z.8}) and (\ref{theta0.5_z.9}) gives
	\[
	\text{P}(V_{m,0.5}^*(\mathbf{x}^0)>c^*_{m,0.5}(\alpha; F_{\mathbf{x}}))\leq C m^{-\eta}\log m(c^*_{m,0.5}(\alpha; F_{\mathbf{x}}))^{-2}.
	\]
	Thus, there exists $c^*_{m,0.5}(\alpha; F_{\mathbf{x}})=C_{0.5,\alpha}\sqrt{ m^{-\eta}\log m}$ for a large enough constant $C_{0.5,\alpha}$ such that $\text{P}(V_{m,0.5}^*(\mathbf{x}^0)>c^*_{m,0.5}(\alpha; F_{\mathbf{x}}))<\alpha$. 
	
	\subsubsection{Proof of Lemma \ref{lemma:var_sum_w}}
		Here
		\begin{equation}\label{theta0.5_z.5}
		\text{Var}\left(\sum_{i=1}^m\mathbb{I}\{x_i^0> t\}\right)=\sum_{i=1}^m\text{Var}\left(\mathbb{I}\{x_i^0> t\}\right)+\sum_{i\neq j}\text{Cov}\left(\mathbb{I}\{x_i^0> t\},\mathbb{I}\{x_j^0> t\}\right).
		\end{equation} 
		As for the first term, we have
		\begin{equation}\label{theta0.5_z.6}
		\sum_{i=1}^{m}\text{Var}\left(\mathbb{I}\{x_i^0> t\}\right)=\sum_{i=1}^{m}\text{P}(x_i^0> t)(1-\text{P}(x_i^0> t))\leq\sum_{i=1}^{m}\text{P}(x_i^0> t)=m\bar{\Phi}(t)\leq Cme^{-t^2/2},
		\end{equation}
		where the last inequality holds due to Mill's ratio. As for the second term, for $i\neq j=1,\cdots,m$,
		\begin{eqnarray}
		\text{Cov}\left(\mathbb{I}\{x_i^0> t\},\mathbb{I}\{x_j^0> t\}\right)&=&\text{Cov}\left(\mathbb{I}\{x_i^0\leq t\},\mathbb{I}\{x_j^0\leq t\}\right)\nonumber\\
		&=&\text{P}(x_i^0\leq t,x_j^0\leq t)-\text{P}(x_i^0\leq t)\text{P}(x_j^0\leq t).\nonumber
		\end{eqnarray}
		Since $x_i^0$ and $x_j^0$ are standard normal random variables with $\text{Cov}(x_i^0,x_j^0)=\Sigma_{ij}$, we can use Corollary 2.1 in \cite{li2002normal} to bound the right hand side to get
		\begin{equation}
		\text{Cov}\left(\mathbb{I}\{x_i^0> t\},\mathbb{I}\{x_j^0> t\}\right)\leq C|\Sigma_{ij}|e^{-t^2/2}.\nonumber
		\end{equation}
		Thus
		\begin{equation}\label{theta0.5_z.7}
		\sum_{i\neq j}\text{Cov}\left(\mathbb{I}\{x_i^0> t\},\mathbb{I}\{x_j^0> t\}\right)\leq C\sum_{i\neq j}|\Sigma_{ij}|e^{-t^2/2}.
		\end{equation}
		Combining (\ref{theta0.5_z.5}), (\ref{theta0.5_z.6}), and (\ref{theta0.5_z.7}) gives (\ref{var_sum_w}).
	
	\subsection{Proof of Lemma \ref{lemma:suf_nec_1}} \label{sec:proof_7.5}
	
	Firstly, we show that $\text{P}(\hat{\pi}^*_1(\mathbf{x})>(1-\epsilon)\pi)\rightarrow1$ for any constant $\epsilon>0$ if $h>\gamma$. Let
	\begin{equation}\label{tau1}
	\tau_m=\dfrac{\mu+\underline{\mu}_1}{2}\rightarrow\infty,
	\end{equation}	
	where
	\begin{equation}\label{mu1.bar}
	\underline{\mu}_1=\sqrt{2\gamma\log m}.
	\end{equation}
	Then, under the condition $h > \gamma$, we have
	\begin{equation} \label{3.3}
	\mu-\tau_m=\dfrac{\mu-\underline{\mu}_1}{2}\rightarrow\infty
	\end{equation}
	and 
	\begin{equation} \label{3.4}
	\tau_m-\underline{\mu}_1=\dfrac{\mu-\underline{\mu}_1}{2}\rightarrow\infty,
	\end{equation}
	and (\ref{3.4}) implies 
	\begin{equation} \label{3.6}
	\bar{\Phi}(\tau_m)/\bar{\Phi}(\underline{\mu}_1)\rightarrow0. 
	\end{equation}
	
	Now, by Mill's ratio, we have
	\[
	\bar{\Phi}(\underline{\mu}_1)\leq\dfrac{\phi(\underline{\mu}_1)}{\underline{\mu}_1}=\dfrac{Ce^{-\gamma\log m}}{\sqrt{\log m}}=\dfrac{C\pi}{\sqrt{\log m}}=\dfrac{C\pi}{c^*_{m,1}(\alpha; F_{\mathbf{x}})},
	\]
	where the last step is by the given value of $c_{m,1}^*(\alpha; F_{\mathbf{x}})$. Then, by (\ref{3.6}), it follows that
	\begin{equation} \label{3.1}
	\pi^{-1}c^*_{m,1}(\alpha; F_{\mathbf{x}})\bar{\Phi}(\tau_m)\rightarrow0.
	\end{equation}
	
	By definition of $\hat{\pi}^*_1(\mathbf{x})$, we have for $t>0$,
	\begin{eqnarray}
	\hat{\pi}^*_1(\mathbf{x})&>&m^{-1}\sum\limits_{i=1}^{m}\mathbb{I}\{x_i> t\}-\bar{\Phi}(t)-c^*_{m,1}(\alpha; F_{\mathbf{x}})\bar{\Phi}(t)\nonumber\\
	&=&\dfrac{1-\pi}{m-s}\sum\limits_{i\in I_0}^{}\mathbb{I}\{x_i> t\}+\dfrac{\pi}{s}\sum\limits_{i\in I_1}^{}\mathbb{I}\{x_i> t\}-\bar{\Phi}(t)-c^*_{m,1}(\alpha; F_{\mathbf{x}})\bar{\Phi}(t),\nonumber
	\end{eqnarray}
	which implies
	\begin{eqnarray}\label{theta1_z.4}
	\dfrac{\hat{\pi}^*_1(\mathbf{x})}{\pi}-1&>&-\dfrac{1}{\pi}c^*_{m,1}(\alpha; F_{\mathbf{x}})\bar{\Phi}(t)-\bar{\Phi}(t)\nonumber\\
	&+&\dfrac{1-\pi}{\pi}\left\{(m-s)^{-1}\sum\limits_{i\in I_0}^{}\mathbb{I}\{x_i> t\}-\bar{\Phi}(t)\right\}\nonumber\\
	&+&s^{-1}\sum\limits_{i\in I_1}^{}\mathbb{I}\{x_i> t\}-1.
	\end{eqnarray}
	Then it suffices to show that the right hand side of (\ref{theta1_z.4}) is of $o_p(1)$ at $t=\tau_m$.
	
	The first term
	\begin{equation}\label{theta1_z.5}
	A_1=-\dfrac{1}{\pi}c^*_{m,1}(\alpha; F_{\mathbf{x}})\bar{\Phi}(\tau_m)-\bar{\Phi}(\tau_m)=o(1)
	\end{equation}
	holds directly by (\ref{3.1}) and $\tau_m\rightarrow\infty$  as $m\rightarrow\infty$. 
	
	As for the second term, note that 
	\begin{eqnarray}\label{theta1_z.6}
	A_2&=&\dfrac{1-\pi}{\pi}\left\{(m-s)^{-1}\sum\limits_{i\in I_0}^{}\mathbb{I}\{x_i> \tau_m\}-\bar{\Phi}(\tau_m)\right\}\nonumber\\
	&=&\dfrac{1-\pi}{\pi}\left\{m^{-1}\sum\limits_{i=1}^{m}\mathbb{I}\{x_i^0> \tau_m\}-\bar{\Phi}(\tau_m)\right\}\nonumber\\
	&+&\dfrac{1-\pi}{\pi}\left\{(m-s)^{-1}\sum\limits_{i\in I_0}^{}\mathbb{I}\{x_i> \tau_m\}-m^{-1}\sum\limits_{i=1}^{m}\mathbb{I}\{x_i^0> \tau_m\}\right\}.
	\end{eqnarray}
	On the one hand,
	\begin{eqnarray}\label{theta1_z.7}
	&&\text{P}\left(\dfrac{\left|m^{-1}\sum\limits_{i=1}^{m}\mathbb{I}\{x_i^0> \tau_m\}-\bar{\Phi}(\tau_m)\right|}{\pi}>\dfrac{c^*_{m,1}(\alpha; F_{\mathbf{x}})\bar{\Phi}(\tau_m)}{\pi}\right)\nonumber\\
	&=&\text{P}\left(\dfrac{\left|m^{-1}\sum\limits_{i=1}^{m}\mathbb{I}\{x_i^0> \tau_m\}-\bar{\Phi}(\tau_m)\right|}{\bar{\Phi}(\tau_m)}>c^*_{m,1}(\alpha; F_{\mathbf{x}})\right)\nonumber\\
	&\le&\text{P}(V^*_{m,1}(\mathbf{x}^0)>c^*_{m,1}(\alpha; F_{\mathbf{x}}))<\alpha = \alpha_m\to 0,
	\end{eqnarray}
	where the last two inequalities hold by definition of $V^*_{m,1}(\mathbf{x}^0)$ and property (ii) of $c^*_{m,1}(\alpha; F_{\mathbf{x}})$, respectively. Now, by (\ref{3.1}) and (\ref{theta1_z.7}), we have
	\begin{equation}\label{theta1_z.8}
	\dfrac{1}{\pi}\left|m^{-1}\sum\limits_{i=1}^{m}\mathbb{I}\{x_i^0> \tau_m\}-\bar{\Phi}(\tau_m)\right|=o_p(1).
	\end{equation}
	On the other hand,
	\begin{eqnarray}\label{theta1_z.9}
	&&\dfrac{1}{\pi}\left|(m-s)^{-1}\sum\limits_{i\in I_0}^{}\mathbb{I}\{x_i^0> \tau_m\}-m^{-1}\sum\limits_{i=1}^{m}\mathbb{I}\{x_i^0> \tau_m\}\right|\nonumber\\
	&=&\dfrac{1}{\pi}\left|\pi(m-s)^{-1}\sum\limits_{i\in I_0}^{}\mathbb{I}\{x_i^0> \tau_m\}-m^{-1}\sum\limits_{i\in I_1}\mathbb{I}\{x_i^0> \tau_m\}\right|\nonumber\\
	&\le&(m-s)^{-1}\sum\limits_{i\in I_0}^{}\mathbb{I}\{x_i^0> \tau_m\}+s^{-1}\sum\limits_{i\in I_1}\mathbb{I}\{x_i^0> \tau_m\}
	\end{eqnarray}
	Using Markov’s inequality, we have for $a>0$,
	\begin{equation}
	\text{P}\left(s^{-1}\sum\limits_{i\in I_1}\mathbb{I}\{x_i^0> \tau_m\}>a\right)\leq\dfrac{1}{as}\sum\limits_{i\in I_1}\text{P}(x_i^0> \tau_m)=a^{-1}\bar{\Phi}(\tau_m)=o(1),\nonumber
	\end{equation}
	where the last step holds since  $\tau_m\rightarrow\infty$. This implies that
	\[
	s^{-1}\sum\limits_{i\in I_1}\mathbb{I}\{x_i^0> \tau_m\}=o_p(1).
	\] 
	With similar arguments, we can show that
	\[
	(m-s)^{-1}\sum\limits_{i\in I_0}^{}\mathbb{I}\{x_i^0> \tau_m\}=o_p(1).
	\]
	Combining the above with (\ref{theta1_z.9}) implies that
	\[
	\dfrac{1}{\pi}\left|(m-s)^{-1}\sum\limits_{i\in I_0}^{}\mathbb{I}\{x_i^0> \tau_m\}-m^{-1}\sum\limits_{i=1}^{m}\mathbb{I}\{x_i^0> \tau_m\}\right|=o_p(1).
	\]
	Since the joint distribution of $x_i, i\in I_0$, is the same as the joint distribution of $x_i^0, i\in I_0$, we have 
	\begin{equation}\label{theta1_z.10}
	\dfrac{1}{\pi}\left|(m-s)^{-1}\sum\limits_{i\in I_0}^{}\mathbb{I}\{x_i> \tau_m\}-m^{-1}\sum\limits_{i=1}^{m}\mathbb{I}\{x_i^0> \tau_m\}\right|=o_p(1).
	\end{equation}
	Combining (\ref{theta1_z.6}), (\ref{theta1_z.8}), and (\ref{theta1_z.10}) gives
	\begin{equation}\label{theta1_z.11}
	A_2=\dfrac{1-\pi}{\pi}\left\{(m-s)^{-1}\sum\limits_{i\in I_0}^{}\mathbb{I}\{x_i> \tau_m\}-\bar{\Phi}(\tau_m)\right\}=o_p(1).
	\end{equation}
	
	As for the last term,
	using Markov's inequality, for $a>0$,
	\begin{eqnarray}\label{theta1_z.12}
	\text{P}(|A_3|>a)&=&\text{P}\left(1-s^{-1}\sum\limits_{i\in I_1}^{}\mathbb{I}\{x_i> \tau_m\}>a\right)\nonumber\\
	&\leq&a^{-1}\left(1-s^{-1}\sum\limits_{i\in I_1}^{}\text{P}(x_i> \tau_m)\right)
	\leq a^{-1}\left(1-\min\limits_{i\in I_1}^{}\text{P}(x_i> \tau_m)\right)\nonumber\\
	&=&a^{-1}\left(1-\min\limits_{i\in I_1}^{}\text{P}(N(0,1)> \tau_m-\mu_i)\right)\nonumber\\
	&\leq& a^{-1}\left(1-\text{P}(N(0,1)> \tau_m-\mu)\right)=o(1),
	\end{eqnarray}
	where the last step holds due to (\ref{3.3}). Thus (\ref{theta1_z.12}) implies
	\begin{equation}\label{theta1_z.13}
	A_3=s^{-1}\sum\limits_{i\in I_1}^{}\mathbb{I}\{x_i> \tau_m\}-1=o_p(1).
	\end{equation}
	
	Combining (\ref{theta1_z.4}), (\ref{theta1_z.5}), (\ref{theta1_z.11}) and (\ref{theta1_z.13}), we conclude that $\text{P}(\hat{\pi}^*_1(\mathbf{x})>(1-\epsilon)\pi)\rightarrow1$ for any constant $\epsilon>0$.
	
	Next we show that $\text{P}(\hat{\pi}^*_1(\mathbf{x})>(1-\epsilon)\pi)\rightarrow1$ for any constant $\epsilon>0$ almost only if $h>\gamma$. It is sufficient to show that when $h<\gamma$, 
	\begin{equation} \label{3.5}
	\hat{\pi}^*_1(\mathbf{x})/ \pi = o_p(1).
	\end{equation} 
	To facilitate the proof, we introduce the following lemma, which follows directly from Lemma A.3 in \cite{jeng2023estimating}. 
	
	\begin{lemma} \label{lemma:nec_theta}
		Consider model (\ref{def:model}) along with the calibrations in (\ref{def:gamma}) - (\ref{def:eta}). If 
		\[
		\bar{\Phi}^{-1}\left(\{\pi/c^*_{m,\theta}(\alpha; F_{\mathbf{x}})\}^{1/\theta}\right)-\mu\to\infty,
		\]
		then $\hat{\pi}^*_\theta({\mathbf{x}})/\pi = o_p(1)$.
	\end{lemma} 
	
	Now, when $h<\gamma$, we have 
	\[
	\underline{\mu}_1-\mu\to\infty,
	\] 
	where $\underline{\mu}_1$ is defined in (\ref{mu1.bar}). 
	On the other hand, by Mill's ratio,
	\[
	\bar{\Phi}(\underline{\mu}_1)\geq\phi(\underline{\mu}_1)\left(\dfrac{1}{\underline{\mu}_1}-\dfrac{1}{\underline{\mu}^3_1}\right)=\dfrac{\phi(\underline{\mu}_1)}{\underline{\mu}_1+o(\underline{\mu}_1)}=\dfrac{Ce^{-\gamma\log m}}{\sqrt{\log m}}=\dfrac{C\pi}{\sqrt{\log m}}=\dfrac{C\pi}{c^*_{m,1}(\alpha; F_{\mathbf{x}})},
	\]
	which implies
	\[
	\underline{\mu}_1\leq\bar{\Phi}^{-1}\left(\dfrac{C\pi}{c^*_{m,1}(\alpha; F_{\mathbf{x}})}\right).
	\]
	Therefore, we have 
	\[
	\bar{\Phi}^{-1}\left(\pi/c^*_{m,1}(\alpha; F_{\mathbf{x}})\right)-\mu\to\infty.
	\]
	Then (\ref{3.5}) follows from Lemma \ref{lemma:nec_theta} with $\theta=1$. This concludes the proof of Lemma \ref{lemma:suf_nec_1}.

	\subsection{Proof of Lemma \ref{lemma:suf_nec_05}} \label{sec:proof_7.6}
	
	Firstly, we show that $\text{P}(\hat{\pi}^*_{0.5}(\mathbf{x})>(1-\epsilon)\pi)\rightarrow1$ for any constant $\epsilon>0$ if $h>[2\gamma-\eta]_+$. In this section, let
	\begin{equation}\label{tau2}
	\tau_m=\dfrac{\mu+\underline{\mu}_2}{2}\rightarrow\infty.
	\end{equation}
	where
	\begin{equation}\label{mu2.bar}
	\underline{\mu}_2=\sqrt{2[2\gamma-\eta]_+\log m+2\log\log m},
	\end{equation}
	Then, under the condition $h>[2\gamma-\eta]_+$, we have
	\[
	\mu-\tau_m=\dfrac{\mu-\underline{\mu}_2}{2}\rightarrow\infty
	\]
	and
	\[
	\tau_m-\underline{\mu}_2=\dfrac{\mu-\underline{\mu}_2}{2}\rightarrow\infty.
	\]
	The latter implies
	\begin{equation} \label{4.3}
	\bar{\Phi}(\tau_m)/\bar{\Phi}(\underline{\mu}_2)\rightarrow 0. 
	\end{equation}
	
	Now, given the value of the bounding sequence $c_{m,0.5}^*(\alpha; F_{\mathbf{x}})$, we have, 
	\[
	\bar{\Phi}(\underline{\mu}_2)\leq\dfrac{\phi(\underline{\mu}_2)}{\underline{\mu}_2}\leq\dfrac{Ce^{-2\gamma\log m +\eta\log m - \log\log m/2}}{\sqrt{\log m + o(\log m)}}=\dfrac{C\pi^2}{m^{-\eta}\log m}=\dfrac{C\pi^2}{\{c^*_{m,0.5}(\alpha; F_{\mathbf{x}})\}^2}
	\]
	for $2\gamma-\eta>0$, and 
	\[
	\bar{\Phi}(\underline{\mu}_2)\leq\dfrac{\phi(\underline{\mu}_2)}{\underline{\mu}_2}\leq\dfrac{Ce^{-\log\log m}}{\sqrt{\log\log m }}\leq\dfrac{C}{\log m}\leq\dfrac{C\pi^2}{m^{-\eta}\log m}=\dfrac{C\pi^2}{\{c^*_{m,0.5}(\alpha; F_{\mathbf{x}})\}^2}
	\]
	for $2\gamma-\eta\leq 0$. These, combined with (\ref{4.3}), give 
	\begin{equation} \label{4.4}
	\pi^{-1}c^*_{m,0.5}(\alpha; F_{\mathbf{x}})\sqrt{\bar{\Phi}(\tau_m)}\rightarrow0.
	\end{equation}
	
	By definition of $\hat{\pi}^*_{0.5}(\mathbf{x})$, we have for $t>0$,
	\begin{eqnarray}
	\hat{\pi}^*_{0.5}(\mathbf{x})&>&m^{-1}\sum\limits_{i=1}^{m}\mathbb{I}\{x_i> t\}-\bar{\Phi}(t)-c^*_{m,0.5}(\alpha; F_{\mathbf{x}})\sqrt{\bar{\Phi}(t)}\nonumber\\
	&=&\dfrac{1-\pi}{m-s}\sum\limits_{i\in I_0}^{}\mathbb{I}\{x_i> t\}+\dfrac{\pi}{s}\sum\limits_{i\in I_1}^{}\mathbb{I}\{x_i> t\}-\bar{\Phi}(t)-c^*_{m,0.5}(\alpha; F_{\mathbf{x}})\sqrt{\bar{\Phi}(t)},\nonumber
	\end{eqnarray}
	which implies
	\begin{eqnarray}\label{theta0.5_z.10}
	\dfrac{\hat{\pi}^*_{0.5}(\mathbf{x})}{\pi}-1&>&-\dfrac{1}{\pi}c^*_{m,0.5}(\alpha; F_{\mathbf{x}})\sqrt{\bar{\Phi}(t)}-\bar{\Phi}(t)\nonumber\\
	&+&\dfrac{1-\pi}{\pi}\left\{(m-s)^{-1}\sum\limits_{i\in I_0}^{}\mathbb{I}\{x_i> t\}-\bar{\Phi}(t)\right\}\nonumber\\
	&+&s^{-1}\sum\limits_{i\in I_1}^{}\mathbb{I}\{x_i> t\}-1.
	\end{eqnarray}
	Similarly to the arguments that lead to $\text{P}(\hat{\pi}^*_1(\mathbf{x})>(1-\epsilon)\pi)\rightarrow1$ in the proof of Lemma \ref{lemma:suf_nec_1}, and by using the result in (\ref{4.4}), we can show that the right hand side of (\ref{theta0.5_z.10}) is $o_p(1)$ when evaluated at $t=\tau_m$. 
	Therefore, $\text{P}(\hat{\pi}^*_{0.5}(\mathbf{x})>(1-\epsilon)\pi)\rightarrow1$ for any constant $\epsilon>0$.
	
	Next we show that $\text{P}(\hat{\pi}^*_{0.5}(\mathbf{x})>(1-\epsilon)\pi)\rightarrow1$ for any constant $\epsilon>0$ almost only if $h>[2\gamma-\eta]_+$. 
	It is sufficient to show that when $h<[2\gamma-\eta]_+$, 
	\begin{equation} \label{4.6} 
	\hat{\pi}^*_{0.5}(\mathbf{x})/ \pi = o_p(1).
	\end{equation} 
	By Mill's ratio and $h<[2\gamma-\eta]_+$,
	\begin{eqnarray}
	\dfrac{\{\pi/c^*_{m,0.5}(\alpha; F_{\mathbf{x}})\}^2}{\bar{\Phi}(\mu)}&\leq&\dfrac{\{\pi/c^*_{m,0.5}(\alpha; F_{\mathbf{x}})\}^2}{\phi(\mu)}\left(\dfrac{1}{\mu}-\dfrac{1}{\mu^3}\right)^{-1}=\dfrac{\pi^2}{\phi(\mu)}\dfrac{\mu+o(\mu)}{\{c^*_{m,0.5}(\alpha; F_{\mathbf{x}})\}^2}\nonumber\\
	&=&\dfrac{m^{-2\gamma}}{m^{-h}}\dfrac{\sqrt{2h\log m}+o(\sqrt{\log m})}{C^2_{0.5,\alpha}m^{-\eta}\log m}=O\left(\dfrac{m^{h-(2\gamma-\eta)}}{\sqrt{\log m}}\right)=o(1),\nonumber
	\end{eqnarray}
	which implies 
	\[
	\bar{\Phi}^{-1}\left(\{\pi/c^*_{m,0.5}(\alpha; F_{\mathbf{x}})\}^2\right)-\mu\to\infty.
	\]
	Then (\ref{4.6}) follows from Lemma \ref{lemma:nec_theta} with $\theta=0.5$. This concludes the proof of Lemma \ref{lemma:suf_nec_05}.

	\subsection{Proof of Theorem \ref{thm:pi_hat_z}} \label{sec:proof_thm3.1}
	
	First, it is easy to see that (\ref{eq:lowerbd_z}) follows Lemma \ref{lemma:lowerbd_z} in discritized form. 
	
	As for (\ref{cond:ah}), 
	recall that after the dependence adjustment, we have $(Z_1,\cdots,Z_m)^T\sim N((a_1\mu_1,\cdots,a_m\mu_m)^T,\tilde{\bm{\Sigma}})$ as in (\ref{def:model.Z}). Under the calibration in (\ref{def:h}), for  $i = \{1, \cdots, m\}$,
	\[
	a_i \mu_i \ge a_{\min} \cdot \mu = \sqrt{2(a^2_{\min} \cdot h) \log m}. 
	\]
	With $c^*_{m,1}(\alpha/2; F_{{\mathbf{z}}}) = C'_{1,\alpha}\sqrt{\log m}$ and $c^*_{m,0.5}(\alpha/2; F_{{\mathbf{z}}}) = C'_{0.5,\alpha}\sqrt{m^{-\tilde \eta}\log m}$, similar reasoning as in the proof of Lemma \ref{lemma:boundingSeq} shows that both  $c^*_{m,1}(\alpha; F_{{\mathbf{x}}})$ and $c^*_{m,0.5}(\alpha; F_{{\mathbf{x}}})$ satisfy properties (i) and (ii) in  Lemma \ref{lemma:lowerbd_z} in its discritized form.
	Furthermore, similar reasoning as in the proof of (\ref{cond:h}) in Proposition \ref{thm:pi_hat_x} leads to (\ref{cond:ah}). The repeated details are omitted to save space.

	\subsection{Proof of Corollary \ref{thm:compare}} \label{sec:proof_cor3.1}
	
	Given $\bm{\Sigma}$, we have 
	\[
	\|\bm{\Sigma}\|_1/m^2=\rho(m^d)^2(m^{1-d})/m^2+(1-\rho)m/m^2=m^{d-1}\{\rho+(1-\rho)m^{-d}\}. 
	\]
	Thus, $\eta=1-d-o(1)$, and Proposition \ref{thm:pi_hat_x} implies (\ref{cond:h_example}).

	To prove the result in (\ref{cond:ah_example}), by Theorem \ref{thm:pi_hat_z}, it is sufficient to show 
	\begin{equation} \label{1.1}
	a_{\min}\geq1/\sqrt{1-\rho}
	\end{equation} 
	and 
	\begin{equation} \label{1.2}
	\tilde{\eta}=1-o(1). 
	\end{equation} 
	Here, the eigenvalues of $\bm{\Sigma}$ are the roots of $\det(\bm{\Sigma}-\lambda\mathbf{I}_m)=0$. Since $\bm{\Sigma}-\lambda\mathbf{I}_m$ is also block-diagonal, we have
	\[
	\det(\bm{\Sigma}-\lambda\mathbf{I}_m)=\{\det(\mathbf{B}-\lambda\mathbf{I}_{m^d})\}^{m^{1-d}}=\{(m^d\rho+1-\rho-\lambda)(1-\rho-\lambda)^{m^d-1}\}^{m^{1-d}},
	\]
	where the last step follows due to
	\begin{eqnarray}
	\det(\mathbf{B}-\lambda\mathbf{I}_{m^d})&=&\det(\rho\mathbf{1}_{m^d}\mathbf{1}_{m^d}^T+(1-\rho-\lambda)\mathbf{I}_{m^d})\nonumber\\
	&=&\{1+\rho\mathbf{1}_{m^d}^T((1-\rho-\lambda)\mathbf{I}_{m^d})^{-1}\mathbf{1}_{m^d}\}\det((1-\rho-\lambda)\mathbf{I}_{m^d})\nonumber\\
	&=&(m^d\rho+1-\rho-\lambda)(1-\rho-\lambda)^{m^d-1}.\nonumber
	\end{eqnarray} 
	As we can see, there are $m^{1-d}$ eigenvalues that are greater than 1, i.e., $\lambda_1=\cdots=\lambda_{m^{1-d}}=m^d\rho-\rho+1$. Other eigenvalues are $1-\rho$, which are smaller than 1. As for the $b$-th eigenvalue $\lambda_b=m^d\rho-\rho+1$, the corresponding eigenvector $\bm{\gamma}_b$ has $1/\sqrt{m^d}$ in entries correspond to block $b$ and $0$ in other entries, where $b=1,\cdots,m^{1-d}$. This is because for $\bm{\gamma}=\mathbf{1}_{m^d}/\sqrt{m^d}$, s.t. $\|\bm{\gamma}\|_2=1$, we have
	\[
	\{\mathbf{B}-(m^d\rho-\rho+1)\mathbf{I}_{m^d}\}\bm{\gamma}=(\rho\mathbf{1}_{m^d}\mathbf{1}_{m^d}^T-m^d\rho\mathbf{I}_{m^d})\mathbf{1}_{m^d}/\sqrt{m^d}=\mathbf{0}.
	\]
	By choosing $k=m^{1-d}$, we have
	\[
	\mathbf{A}=\bm{\Sigma}-\sum_{j=1}^{k}\lambda_j\bm{\gamma}_j\bm{\gamma}_j^T,
	\]
	which is also block diagonal with $k$ blocks $\mathbf{A}_1,\cdots,\mathbf{A}_k$, where
	\[
	\mathbf{A}_b=\mathbf{B}-(m^d\rho-\rho+1)(\mathbf{1}_{m^d}/\sqrt{m^d})(\mathbf{1}_{m^d}/\sqrt{m^d})^T=\dfrac{\rho-1}{m^d}\mathbf{1}_{m^d}\mathbf{1}_{m^d}^T+(1-\rho)\mathbf{I}_{m^d}.
	\]
	Here for $i=1,\cdots,m$ and $j=1,\cdots,k$, we have $b_{ij}=\sqrt{\lambda_j}\bm{\gamma}_{j,i}$. Note that $i$ only belongs to one block, which is denoted by $b(i)\in\{1,\cdots,k\}$, then
	\[
	\|\mathbf{b}_i\|_2^2 =\sum_{j=1}^{k}\lambda_j(\bm{\gamma}_{j,i})^2=\lambda_{b(i)}(\bm{\gamma}_{b(i),i})^2=(m^d\rho-\rho+1)(1/\sqrt{m^d})^2=\dfrac{m^d\rho-\rho+1}{m^d},
	\]
	and
	\[
	a_i=\dfrac{1}{\sqrt{1-\|\mathbf{b}_i\|_2^2}}=\dfrac{1}{\sqrt{(1-\rho)(1-1/m^d)}}.
	\]
	Therefore, (\ref{1.1}) follows. 
	
	Next, we can standardize $\mathbf{A}$ to get $\tilde{\bm{\Sigma}}$, which is also block diagonal with $k$ blocks $\tilde{\bm{\Sigma}}_1,\cdots,\tilde{\bm{\Sigma}}_k$, where
	\[
	\tilde{\bm{\Sigma}}_b=-\dfrac{1}{m^d-1}\mathbf{1}_{m^d}\mathbf{1}_{m^d}^T+\dfrac{m^d}{m^d-1}\mathbf{I}_{m^d}.
	\]
	So
	\[
	\|\tilde{\bm{\Sigma}}\|_1/m^2=\sum_{b=1}^{k}\|\tilde{\bm{\Sigma}}_b\|_1/m^2=2/m.
	\]
	Then, (\ref{1.2}) follows. 
	
	\subsection{Proof of Theorem \ref{thm:pi_hat_L}} \label{sec:proof_thm4.1}
	
	Recall that the discretized $\hat \pi^*_L(\widehat{\mathbf{z}})$ is defined as 
	\begin{equation} \label{def:pi_hat*_L}
	\hat\pi^*_L(\widehat{\mathbf{z}})=\max\{\hat{\pi}^*_{L,0.5}(\widehat{\mathbf{z}}),~\hat{\pi}^*_{L,1}(\widehat{\mathbf{z}})\},
	\end{equation} 
	where 
	\begin{equation} \label{def:hat_pi*_theta_tilde}
	\hat{\pi}^*_{L,\theta}(\widehat{\mathbf{z}})=\max\limits_{t\in\mathbb{T}'}\dfrac{m^{-1}\sum\limits_{i=1}^{m}\mathbb{I}\{\hat z_i> t\}-\bar{\Phi}(t)-c^*_{m,\theta}(\alpha/2; F_{\widehat{\mathbf{z}}})[\bar{\Phi}(t)]^\theta}{1-\bar{\Phi}(t)}
	\end{equation}
	for a given $\theta \in [0,1]$, with $\mathbb{T}'=[\sqrt{2[\gamma- (\tilde \eta-\gamma)_+]_+\log m},\sqrt{5\log m}]\cap\mathbb{N}$. Note that in (C1'), the order of $\sqrt{m^{1-\tilde \eta}}$ is less than that of $\sqrt{m^{1-2(\tilde \eta - \gamma)_+}}$, which indicates $2\gamma>\tilde \eta$, and $[\gamma- (\tilde \eta-\gamma)_+]_+>0$. The bounding sequence $c^*_{m,\theta}(\alpha/2; F_{\widehat{\mathbf{z}}})$ is constructed such that $\text{P}(V^*_{m,\theta}(\widehat{\mathbf{z}}^0) > c^*_{m,\theta}(\alpha/2; F_{\widehat{\mathbf{z}}})) < \alpha/2$ with
	\begin{equation}
	V^*_{m,\theta}(\widehat{\mathbf{z}}^0) = \max\limits_{t\in\mathbb{T}'}\dfrac{\left|m^{-1}\sum\limits_{i=1}^{m}\mathbb{I}\{\hat z_i^0> t\}-\bar{\Phi}(t)\right|}{[\bar{\Phi}(t)]^\theta}. \nonumber
	\end{equation}
	It is easy to see that (\ref{eq:lowerbd_L}) follows Lemma \ref{lemma:lowerbd_z_hat} in discretized form.
	
	We consider (\ref{eq:upper_L}) next. Since
	\[
	\text{P}((1-\epsilon)\pi < \hat{\pi}^*_{L}(\widehat{\mathbf{z}})<\pi) \ge \text{P}((1-\epsilon)\pi < \hat{\pi}^*_{L}(\widehat{\mathbf{z}})<\pi \mid \|\widehat{\mathbf{w}} - \mathbf{w}\|_2\leq C\sqrt{k/m}) \cdot \text{P}(\|\widehat{\mathbf{w}} - \mathbf{w}\|_2\leq C\sqrt{k/m})
	\]
	and, by Lemma \ref{lemma:w_hat-w},
	\[
	\text{P}(\|\widehat{\mathbf{w}} - \mathbf{w}\|_2\leq C\sqrt{k/m}) \to 1, 
	\] 
	Then it is sufficient to show
	\begin{equation} \label{5.0}
	\text{P}((1-\epsilon)\pi < \hat{\pi}^*_{L}(\widehat{\mathbf{z}})<\pi \mid \|\widehat{\mathbf{w}} - \mathbf{w}\|_2\leq C\sqrt{k/m}) \to 1. 
	\end{equation}
	
	The following lemmas specify  $c^*_{m,1}(\alpha/2; F_{\widehat{\mathbf{z}}})$ and  $c^*_{m,0.5}(\alpha/2; F_{\widehat{\mathbf{z}}})$ and demonstrate the upper bound results of $\hat{\pi}^*_{L,1}(\widehat{\mathbf{z}})$ and $\hat{\pi}^*_{L,0.5}(\widehat{\mathbf{z}})$ given the event $\|\widehat{\mathbf{w}} - \mathbf{w}\|_2\leq C\sqrt{k/m}$. 
	
	\begin{lemma} \label{lemma:boundingSeq_L}
		Consider model (\ref{def:model}) and the dependence-adjusted statistics in (\ref{def:z_hat}) with the estimated principal factors in (\ref{def:est.W}), along with the calibrations in (\ref{def:gamma}) - (\ref{def:h}) and (\ref{def:tilde_eta}). 
		Let $c^*_{m,1}(\alpha/2; F_{\widehat{\mathbf{z}}}) = C''_{1,\alpha}\sqrt{\log m}$ for some constant $C''_{1,\alpha}$ and $c^*_{m,0.5}(\alpha/2; F_{\widehat{\mathbf{z}}}) = C''_{0.5,\alpha}\sqrt{m^{-\tilde \eta}\log m}$ for some constant $C''_{0.5,\alpha}$. 
		Given that $\|\widehat{\mathbf{w}} - \mathbf{w}\|_2\leq C\sqrt{k/m}$ and condition (C1'), both  $c^*_{m,1}(\alpha/2; F_{\widehat{\mathbf{z}}})$ and $c^*_{m,0.5}(\alpha/2; F_{\widehat{\mathbf{z}}})$ satisfy properties (i) and (ii) in  Lemma \ref{lemma:lowerbd_z_hat} in their discritized form. 
	\end{lemma}
	
	\begin{lemma} \label{lemma:suf_1_L}
		Assume the same conditions as in Lemma \ref{lemma:boundingSeq_L},
		we have $\text{P}(\hat{\pi}^*_{L,1}(\widehat{\mathbf{z}}) > (1-\epsilon)\pi)\to 1$ for any constant $\epsilon>0$ if $a^2_{\min} \cdot h > \gamma$. 
	\end{lemma}  
	
	\begin{lemma} \label{lemma:suf_05_L}
		Assume the same conditions as in Lemma \ref{lemma:boundingSeq_L}, 
		we have
		$\text{P}(\hat{\pi}^*_{L,0.5}(\widehat{\mathbf{z}}) > (1-\epsilon)\pi)\to 1$ for any constant $\epsilon>0$ if 
		$a^2_{\min} \cdot h>[2\gamma-\tilde \eta]_+$. 
	\end{lemma}    
	Combining Lemmas \ref{lemma:boundingSeq_L} - \ref{lemma:suf_05_L}, and the definition of $\hat\pi^*_L(\widehat{\mathbf{z}})$ in (\ref{def:pi_hat*_L}) gives (\ref{5.0}) if 
	$a^2_{\min} \cdot h > \min\{\gamma,~ [2\gamma-\tilde \eta]_+\}$, which is equivalent to the condition in (\ref{cond:ah}). Therefore, (\ref{eq:upper_L}) follows. 
	
	\subsection{Proof of Lemma \ref{lemma:boundingSeq_L}} \label{sec:proof_7.8}
	
	First, given that $\|\widehat{\mathbf{w}} - \mathbf{w}\|_2\leq C\sqrt{k/m}$, we have
	\begin{equation} \label{keyLemma.1}
	|\mathbf{b}_i^T(\widehat{\mathbf{w}} - \mathbf{w})| \le \sqrt{\sum_{h=1}^{k} b_{ih}^2} \sqrt{\sum_{h=1}^k |\hat w_h - w_h|^2} \le \|\widehat{\mathbf{w}} - \mathbf{w}\|_2 \le C\sqrt{k/m}, \quad i=1, \cdots, m,
	\end{equation}
	almost surely, where the first inequality is by Cauchy-Schwarz inequality and the second inequality is by the fact that $\sum_{h=1}^{k} b_{ih}^2 \le 1$ for $i=1, \cdots, m$. 
	
	Now, we show that given $\|\widehat{\mathbf{w}} - \mathbf{w}\|_2\leq C\sqrt{k/m}$, the bounding sequence $c_{m,1}^*(\alpha/2; F_{\widehat{\mathbf{z}}})$ satisfies (i) and (ii). 
	
	Property (i) is trivial. 
	
	As for property (ii), let
	\[
	\widehat B(t)=\dfrac{\left|m^{-1}\sum\limits_{i=1}^{m}\mathbb{I}\{\hat z_i^0> t\}-\bar{\Phi}(t)\right|}{\bar{\Phi}(t)}
	\]
	such that $V^*_{m,1}(\widehat{\mathbf{z}}^0)=\max\limits_{t\in\mathbb{T}'}\widehat B(t)$. Argument similar to (\ref{theta1_z.3}) in proof of Lemma \ref{lemma:boundingSeq} gives
	\begin{equation} \label{theta1_z_hat.3}
	\text{P}(V^*_{m,1}(\widehat{\mathbf{z}}^0)>c_{m,1}^*(\alpha/2; F_{\widehat{\mathbf{z}}})) \leq C\sqrt{\log m}(c_{m,1}^*(\alpha/2; F_{\widehat{\mathbf{z}}}))^{-1}\max\limits_{t\in\mathbb{T}'}\text{E}[\widehat B(t)].
	\end{equation}
	
	Note that $\hat{z}_i^0=z_i^0-a_i\mathbf{b}_i^T(\widehat{\mathbf{w}}^0-\mathbf{w}^0)$, so (\ref{keyLemma.1}) gives 
	\begin{eqnarray} \label{theta1_z_hat.4}
	\text{P}(\hat z_i^0> t)&=&\text{P}(z_i^0>t+a_i\mathbf{b}_i^T(\widehat{\mathbf{w}}^0-\mathbf{w}^0))\leq\text{P}(z_i^0>t-a_{\max}|\mathbf{b}_i^T(\widehat{\mathbf{w}}^0-\mathbf{w}^0)|)\nonumber\\
	&\leq&\text{P}(z_i^0>t-a_{\max}C\sqrt{k/m})=\bar{\Phi}(t-a_{\max}C\sqrt{k/m}).
	\end{eqnarray}
	Thus
	\begin{eqnarray}\label{theta1_z_hat.5}
	\text{E}[\widehat B(t)]& \le & [\bar{\Phi}(t)]^{-1}\left(m^{-1}\sum\limits_{i=1}^{m}\text{P}(\hat z_i^0> t)+\bar{\Phi}(t)\right) \nonumber \\ 
	& \le & [\bar{\Phi}(t)]^{-1}\left(m^{-1}\sum\limits_{i=1}^{m}\bar{\Phi}(t-a_{\max}C\sqrt{k/m})+\bar{\Phi}(t)\right) \nonumber \\ 
	& = & \dfrac{\bar{\Phi}(t-a_{\max}C\sqrt{k/m})}{\bar{\Phi}(t)}+1. 
	\end{eqnarray}
	Using Taylor's expansion and Mill's ratio, we have
	\begin{eqnarray}\label{theta1_z_hat.6}
	\dfrac{\bar{\Phi}(t-a_{\max}C\sqrt{k/m})}{\bar{\Phi}(t)}&\leq&\dfrac{\bar{\Phi}(t)+\phi(t-a_{\max}C\sqrt{k/m})a_{\max}C\sqrt{k/m}}{\bar{\Phi}(t)}\nonumber\\
	&\leq& 1+\dfrac{\phi(t-a_{\max}C\sqrt{k/m})a_{\max}C\sqrt{k/m}}{\phi(t)}\left(\dfrac{1}{t}-\dfrac{1}{t^3}\right)^{-1}\nonumber\\
	&=&1+Ce^{a_{\max}Ct\sqrt{k/m}}a_{\max}\sqrt{k/m}(t+o(t))=O(1),
	\end{eqnarray}
	where the last step comes from $a_{\max}\sqrt{k}\leq C\sqrt{m/\log m}$ and $t\in\mathbb{T}'$.
	
	Combining (\ref{theta1_z_hat.3}), (\ref{theta1_z_hat.5}), and (\ref{theta1_z_hat.6}) gives
	\[
	\text{P}(V^*_{m,1}(\widehat{\mathbf{z}}^0)>c_{m,1}^*(\alpha/2; F_{\widehat{\mathbf{z}}})) \leq C\sqrt{\log m}(c_{m,1}^*(\alpha/2; F_{\widehat{\mathbf{z}}}))^{-1}.
	\]
	Thus,  $c_{m,1}^*(\alpha/2; F_{\widehat{\mathbf{z}}})=C''_{1,\alpha}\sqrt{\log m}$, with a large enough constant $C''_{1,\alpha}$, satisfies property (ii) in Lemma \ref{lemma:boundingSeq_L}. 
	
	Next, we show that, given $\|\widehat{\mathbf{w}} - \mathbf{w}\|_2\leq C\sqrt{k/m}$, the bounding sequence $c_{m,0.5}^*(\alpha/2; F_{\widehat{\mathbf{z}}})$ satisfying properties (i) and (ii). 
	
	Property (i) is trivial. 
	
	As for property (ii), let
	\[
	\widehat A(t)=\dfrac{\left|m^{-1}\sum\limits_{i=1}^{m}\mathbb{I}\{\hat z_i^0> t\}-\bar{\Phi}(t)\right|}{[\bar{\Phi}(t)]^{1/2}}
	\]
	such that $V_{m,0.5}^*(\widehat{\mathbf{z}}^0)=\max\limits_{t\in\mathbb{T}'}\widehat A(t)$. By Markov’s inequality, arguments similar to (\ref{theta0.5_z.3}) in proof of Lemma \ref{lemma:boundingSeq} gives
	\begin{equation}\label{theta0.5_z_hat.2}
	\text{P}(V_{m,0.5}^*(\widehat{\mathbf{z}}^0)>c_{m,0.5}^*(\alpha/2; F_{\widehat{\mathbf{z}}}))\leq C\sqrt{\log m}(c_{m,0.5}^*(\alpha/2; F_{\widehat{\mathbf{z}}}))^{-2}\max\limits_{t\in\mathbb{T}'}\text{E}\left[\left\{\widehat A(t)\right\}^2\right].
	\end{equation}
	Here we have for $t\in\mathbb{T}'$,
	\begin{eqnarray}\label{theta0.5_z_hat.3}
	\text{E}\left[\left\{\widehat A(t)\right\}^2\right]&=& [\bar{\Phi}(t)]^{-1}\left\{\text{Var}\left(m^{-1}\sum\limits_{i=1}^{m}\mathbb{I}\{\hat z_i^0> t\}-\bar{\Phi}(t)\right)+\left(\text{E}\left[m^{-1}\sum\limits_{i=1}^{m}\mathbb{I}\{\hat z_i^0> t\}-\bar{\Phi}(t)\right]\right)^2\right\}\nonumber\\&=& m^{-2}[\bar{\Phi}(t)]^{-1}\text{Var}\left(\sum\limits_{i=1}^{m}\mathbb{I}\{\hat z_i^0> t\}\right)+[\bar{\Phi}(t)]^{-1}\left\{m^{-1}\sum\limits_{i=1}^{m}\text{P}(\hat z_i^0> t)-\bar{\Phi}(t)\right\}^2.
	\end{eqnarray}
	
	The following lemma is useful to show the order of the first term in (\ref{theta0.5_z_hat.3}).
		\begin{lemma} \label{lemma:var_sum_w_hat}
			Assume the same conditions as in Lemma \ref{lemma:boundingSeq_L}. For $(\hat z^0_1, \cdots, \hat z^0_m)^T \sim F_{\widehat{\mathbf{z}}}$ as described in \ref{sec:L_approach}, we have
			\begin{equation}\label{var_sum_w_hat}
			\text{Var}\left(\sum\limits_{i=1}^{m}\mathbb{I}\{\hat z_i^0> t\}\right)\leq Cm^{2-\tilde{\eta}}e^{-(t-C/\sqrt{\log m})^2/2}.
			\end{equation}
		\end{lemma}		
		From (\ref{var_sum_w_hat}), we obtain
	\begin{equation}\label{theta0.5_z_hat.7}
	 m^{-2}[\bar{\Phi}(t)]^{-1}\text{Var}\left(\sum\limits_{i=1}^{m}\mathbb{I}\{\hat z_i^0> t\}\right)\leq Cm^{-\tilde{\eta}}[\bar{\Phi}(t)]^{-1}e^{-(t-C/\sqrt{\log m})^2/2}.
	\end{equation}
Using Mill's ratio, we have
\begin{eqnarray}\label{theta0.5_z_hat.8}
[\bar{\Phi}(t)]^{-1}e^{-(t-C/\sqrt{\log m})^2/2} &\le&C[\phi(t)]^{-1}e^{-(t-C/\sqrt{\log m})^2/2}\left(\dfrac{1}{t}-\dfrac{1}{t^3}\right)^{-1}\nonumber\\
&\le&Ce^{t^2/2}e^{-(t-C/\sqrt{\log m})^2/2}(t+o(t))=Ce^{Ct/\sqrt{\log m}}(t+o(t))\nonumber\\
&\leq&C\sqrt{\log m},
\end{eqnarray}
where the last step comes from $t\in\mathbb{T}'$. Combining (\ref{theta0.5_z_hat.7}) and (\ref{theta0.5_z_hat.8}) gives
\begin{equation}\label{theta0.5_z_hat.left}
m^{-2}[\bar{\Phi}(t)]^{-1}\text{Var}\left(\sum\limits_{i=1}^{m}\mathbb{I}\{\hat z_i^0> t\}\right)\leq Cm^{-\tilde{\eta}}\sqrt{\log m}.
\end{equation}

	As for the second term in (\ref{theta0.5_z_hat.3}), we have from (\ref{theta1_z_hat.4}) and by Taylor's expansion,
	\begin{eqnarray}
	m^{-1}\sum\limits_{i=1}^{m}\text{P}(\hat z_i^0> t)-\bar{\Phi}(t)&=&\bar{\Phi}(t-a_{\max}C\sqrt{k/m})-\bar{\Phi}(t)\nonumber\\
	&\leq&\phi(t-a_{\max}C\sqrt{k/m})a_{\max}C\sqrt{k/m}\leq Ce^{-t^2/2}e^{Cta_{\max}\sqrt{k/m}}a_{\max}\sqrt{k/m}\nonumber\\
	&\leq&Ce^{-t^2/2}a_{\max}\sqrt{k/m},\nonumber
	\end{eqnarray}
	where the last step comes from $a_{\max}\sqrt{k}\leq C\sqrt{m/\log m}$ and $t\in\mathbb{T}'$. Using Mill's ratio, we have
	\begin{eqnarray}\label{theta0.5_z_hat.right}
	&&[\bar{\Phi}(t)]^{-1}\left\{m^{-1}\sum\limits_{i=1}^{m}\text{P}(\hat z_i^0> t)-\bar{\Phi}(t)\right\}^2\nonumber\\
	&\leq&[\phi(t)]^{-1}\left\{m^{-1}\sum\limits_{i=1}^{m}\text{P}(\hat z_i^0> t)-\bar{\Phi}(t)\right\}^2\left(\dfrac{1}{t}-\dfrac{1}{t^3}\right)^{-1}\nonumber\\
	&\leq&Ce^{t^2/2}\left(e^{-t^2/2}a_{\max}\sqrt{k/m}\right)^2(t+o(t))=Cm^{-1}e^{-t^2/2}(a_{\max}\sqrt{k})^2(t+o(t))\nonumber\\
	&\leq& Cm^{-\tilde{\eta}}\sqrt{\log m},
	\end{eqnarray}
    where the last step holds since condition (C1') gives $(a_{\max}\sqrt{k})^2 \leq Cm^{1-2(\tilde \eta - \gamma)_+}/\log m$, $t\in\mathbb{T}'$ gives $e^{-t^2/2}\leq m^{(\tilde \eta-\gamma)_+-\gamma}$, and $e^{-t^2/2}(a_{\max}\sqrt{k})^2\leq Cm^{1-\gamma-(\tilde \eta - \gamma)_+}/\log m\leq Cm^{1-\tilde \eta}/\log m$.
    
    Combining (\ref{theta0.5_z_hat.2}), (\ref{theta0.5_z_hat.3}), (\ref{theta0.5_z_hat.left}), and (\ref{theta0.5_z_hat.right}) gives		
	\[
	\text{P}(V^*_{m,0.5}(\widehat{\mathbf{z}}^0)>c_{m,0.5}^*(\alpha/2; F_{\widehat{\mathbf{z}}})) \leq Cm^{-\tilde{\eta}}\log m(c_{m,0.5}^*(\alpha/2; F_{\widehat{\mathbf{z}}}))^{-2}.
	\]
	Thus, the bounding sequence $c_{m,0.5}^*(\alpha/2; F_{\widehat{\mathbf{z}}})=C''_{0.5,\alpha}\sqrt{m^{-\tilde{\eta}}\log m}$ with a large enough constant $C''_{0.5,\alpha}$ satisfies property (ii) in Lemma \ref{lemma:boundingSeq_L}. 
	
	\subsubsection{Proof of Lemma \ref{lemma:var_sum_w_hat}}
		Here
		\begin{equation}\label{theta0.5_z_hat.4}
		\text{Var}\left(\sum_{i=1}^m\mathbb{I}\{\hat{z}_i^0> t\}\right)=\sum_{i=1}^m\text{Var}\left(\mathbb{I}\{\hat{z}_i^0> t\}\right)+\sum_{i\neq j}\text{Cov}\left(\mathbb{I}\{\hat{z}_i^0> t\},\mathbb{I}\{\hat{z}_j^0> t\}\right).
		\end{equation} 
		As for the first term, we have
		\begin{eqnarray}\label{theta0.5_z_hat.5}
		\sum_{i=1}^{m}\text{Var}\left(\mathbb{I}\{\hat{z}_i^0> t\}\right)&=&\sum_{i=1}^{m}\text{P}(\hat{z}_i^0> t)(1-\text{P}(\hat{z}_i^0> t))\leq\sum_{i=1}^{m}\text{P}(\hat{z}_i^0> t)\nonumber\\
		&\leq& m\bar{\Phi}(t-a_{\max}C\sqrt{k/m})\leq \dfrac{m\phi(t-a_{\max}C\sqrt{k/m})}{t-a_{\max}C\sqrt{k/m}}\nonumber\\
		&\leq& Cme^{-(t-C/\sqrt{\log m})^2/2},
		\end{eqnarray}
		where the last three inequalities hold due to (\ref{theta1_z_hat.4}), Mill's ratio, $a_{\max}\sqrt{k}\leq C\sqrt{m/\log m}$ and $t\in\mathbb{T}'$, respectively. As for the second term, for $i\neq j \in \{1,\cdots,m\}$,
		\begin{eqnarray}\label{theta0.5_z_hat.6.1}
		\left|\text{Cov}\left(\mathbb{I}\{\hat{z}_i^0> t\},\mathbb{I}\{\hat{z}_j^0> t\}\right)\right|&=&\left|\text{Cov}\left(\mathbb{I}\{\hat{z}_i^0\leq t\},\mathbb{I}\{\hat{z}_j^0\leq t\}\right)\right|\nonumber\\
		&=&\left|\text{P}(\hat{z}_i^0\leq t,\hat{z}_j^0\leq t)-\text{P}(\hat{z}_i^0\leq t)\text{P}(\hat{z}_j^0\leq t)\right|\nonumber\\
		&=&\Big|\text{P}(z_i^0\leq t+a_i\mathbf{b}_i^T(\widehat{\mathbf{w}}^0-\mathbf{w}^0),z_j^0\leq t+a_j\mathbf{b}_j^T(\widehat{\mathbf{w}}^0-\mathbf{w}^0))\nonumber\\
		&&-\text{P}(z_i^0\leq t+a_i\mathbf{b}_i^T(\widehat{\mathbf{w}}^0-\mathbf{w}^0))\text{P}(z_j^0\leq t+a_j\mathbf{b}_j^T(\widehat{\mathbf{w}}^0-\mathbf{w}^0))\Big|\nonumber\\
		&\leq&\Big|\text{P}(z_i^0\leq t-a_{\max}C\sqrt{k/m},z_j^0\leq t-a_{\max}C\sqrt{k/m})\nonumber\\
		&&-\text{P}(z_i^0\leq t-a_{\max}C\sqrt{k/m})\text{P}(z_j^0\leq t-a_{\max}C\sqrt{k/m})\Big|\nonumber\\
		&\leq&\Big|\text{P}(z_i^0\leq t-C/\sqrt{\log m},z_j^0\leq t-C/\sqrt{\log m})\nonumber\\
		&&-\text{P}(z_i^0\leq t-C/\sqrt{\log m})\text{P}(z_j^0\leq t-C/\sqrt{\log m})\Big|\nonumber\\
		&\leq&C|\tilde\Sigma_{ij}|e^{-(t-C/\sqrt{\log m})^2/2},
		\end{eqnarray}
		where we use Corollary 2.1 in \cite{li2002normal} in the last step, which is similar to proof of Lemma \ref{lemma:boundingSeq}, and the first two inequalities hold by (\ref{keyLemma.1}) and $a_{\max}\sqrt{k}\leq C\sqrt{m/\log m}$, respectively. 
		
		Note that in the first two inequalities in deriving (\ref{theta0.5_z_hat.6.1}), we utilize the monotonicity of
		\[
		f_{ij}(x,y)=\text{P}(z_i^0\leq x, z_j^0\leq y)-\text{P}(z_i^0\leq x)\text{P}(z_j^0\leq y),~x,y>0.
		\] 
		Specifically, $|f_{ij}(x,y)|$ decreases with $x$ and $y$ when $(x,y)\in[t\pm C/\sqrt{\log m}]^2$. Recall that $z_i^0\sim N(0,1)$, $z_j^0\sim N(0,1)$, and $z_i^0$ and $z_j^0$ jointly follow bivariate normal with correlation $\tilde\Sigma_{ij}$. Set
		\[
		z_{ij}^0=\dfrac{z_j^0-\tilde\Sigma_{ij} z_i^0}{\sqrt{1-\tilde\Sigma_{ij}^2}},
		\]
		then $z_{ij}^0\sim N(0,1)$, and $z_i^0$ and $z_{ij}^0$ are independent. We have
		\begin{eqnarray}
		f_{ij}(x,y)&=&\int_{-\infty}^{x}f_{z_i^0}(t)\text{P}\left(z_{ij}^0\leq\dfrac{y-\tilde\Sigma_{ij} z_i^0}{\sqrt{1-\tilde\Sigma_{ij}^2}}\Bigg|z_i^0=t\right)dt-\text{P}(z_i^0\leq x)\text{P}(z_j^0\leq y)\nonumber\\
		&=&\int_{-\infty}^{x}\phi(t)\Phi\left(\dfrac{y-\tilde\Sigma_{ij} t}{\sqrt{1-\tilde\Sigma_{ij}^2}}\right)dt-\Phi(x)\Phi(y).\nonumber
		\end{eqnarray}
		Then we have
		\[
		\dfrac{\partial}{\partial x}f_{ij}(x,y)=\phi(x)\left\{\Phi\left(\dfrac{y-\tilde\Sigma_{ij} x}{\sqrt{1-\tilde\Sigma_{ij}^2}}\right)-\Phi(y)\right\}.
		\]
		If $\tilde\Sigma_{ij}=0$, we have $f_{ij}(x,y)\equiv0$. If $\tilde\Sigma_{ij}<0$, we have $f_{ij}(x,y)<0$, and
		\[
		\dfrac{y-\tilde\Sigma_{ij} x}{\sqrt{1-\tilde\Sigma_{ij}^2}}>y\Rightarrow\dfrac{\partial}{\partial x}f_{ij}(x,y)>0,
		\]
		so $|f_{ij}(x,y)|$ decreases with $x$. If $\tilde\Sigma_{ij}>0$, we have $f_{ij}(x,y)>0$, and
		\[
		\dfrac{\partial}{\partial x}f_{ij}(x,y)<0~\text{when}~x>\dfrac{\tilde\Sigma_{ij}}{1+\sqrt{1-\tilde\Sigma_{ij}^2}}y;
		\]
		\[
		\dfrac{\partial}{\partial x}f_{ij}(x,y)>0~\text{when}~x<\dfrac{\tilde\Sigma_{ij}}{1+\sqrt{1-\tilde\Sigma_{ij}^2}}y.
		\]
		Here we need to consider only
		\[
		(x,y)\in[t-C/\sqrt{\log m},t+C/\sqrt{\log m}]\times[t-C/\sqrt{\log m},t+C/\sqrt{\log m}],
		\]
		where $t\in\mathbb{T}'$. If $0<\tilde\Sigma_{ij}<1$, as $m\to\infty$ we have
		\[
		x>\dfrac{\tilde\Sigma_{ij}}{1+\sqrt{1-\tilde\Sigma_{ij}^2}}y\Rightarrow\dfrac{\partial}{\partial x}f_{ij}(x,y)<0,
		\]
		so $|f_{ij}(x,y)|$ decreases with $x$. Thus $|f_{ij}(x,y)|$ decreases with $x$ for $-1\leq\tilde\Sigma_{ij}<1$. Similar arguments can show that $|f_{ij}(x,y)|$ decreases with $y$. Thus we can bound
		\[
		|f_{ij}(x,y)|\leq |f(t-C/\sqrt{\log m},t-C/\sqrt{\log m})|
		\]
		for $(x,y)\in[t\pm C/\sqrt{\log m}]^2$.
		Thus (\ref{theta0.5_z_hat.6.1}) gives
		\begin{equation}\label{theta0.5_z_hat.6}
		\sum_{i\neq j}\text{Cov}\left(\mathbb{I}\{\hat{z}_i^0> t\},\mathbb{I}\{\hat{z}_j^0> t\}\right)\leq C\sum_{i\neq j}|\tilde\Sigma_{ij}|e^{-(t-C/\sqrt{\log m})^2/2}.
		\end{equation}
		Combining (\ref{theta0.5_z_hat.4}), (\ref{theta0.5_z_hat.5}), and (\ref{theta0.5_z_hat.6}) gives (\ref{var_sum_w_hat}).
	
	\subsection{Proof of Lemma \ref{lemma:suf_1_L}} \label{sec:proof_7.9}
	
	Using similar reasoning as in the proof of Lemma \ref{lemma:suf_nec_1}, we have 
	\[
	\dfrac{\hat{\pi}^*_1(\mathbf{x})}{\pi}-1 > A_1 + A_2 + A_3, 
	\]
	where $A_1$, $A_2$, and $A_3$ are defined similarly as in Section \ref{sec:proof_7.5} with the following substitutions: replace $x_i$ with $\hat{z}_i$, $c^*_{m,1}(\alpha; F_{\mathbf{x}})$ with $c^*_{m, 1}(\alpha/2; F_{\widehat{\mathbf{z}}})$, and $\mu$ with $a_{\min} \cdot \mu$. It is then sufficient to show that $A_1$, $A_2$, and $A_3$ are all of order $o_p(1)$.
	
	Specifically, to show that $A_2=o_p(1)$, the following is applied: 
	\[
	\text{P}\left(s^{-1}\sum\limits_{i\in I_1}\mathbb{I}\{\hat z_i^0> \tau_m\}>a\right)\leq\dfrac{1}{as}\sum\limits_{i\in I_1}\text{P}(\hat z_i^0> \tau_m)\leq a^{-1}\bar{\Phi}(\tau_m-a_{\max}C\sqrt{k/m})=o(1)
	\]
	where the second inequality holds by (\ref{theta1_z_hat.4}), and the last step holds from $a_{\max}\sqrt{k}\leq C\sqrt{m/\log m}$ implied by condition (C1').
	
	Moreover, to prove $A_3=o_p(1)$, we have 
	\begin{eqnarray}
	\text{P}(|A_3|>a)&=&\text{P}\left(1-s^{-1}\sum\limits_{i\in I_1}^{}\mathbb{I}\{\hat z_i> \tau_m\}>a\right)\nonumber\\
	&\leq&a^{-1}\left(1-s^{-1}\sum\limits_{i\in I_1}^{}\text{P}(\hat z_i> \tau_m)\right)=a^{-1}\left(1-s^{-1}\sum\limits_{i\in I_1}^{}\text{P}(z_i> \tau_m+a_i\mathbf{b}_i^T(\widehat{\mathbf{w}}-\mathbf{w}))\right)\nonumber\\
	&\leq& a^{-1}\left(1-\min\limits_{i\in I_1}^{}\text{P}(z_i> \tau_m+a_{\max}|\mathbf{b}_i^T(\widehat{\mathbf{w}}-\mathbf{w})|)\right)\nonumber\\
	&\leq&a^{-1}\left(1-\min\limits_{i\in I_1}^{}\bar{\Phi}( \tau_m+a_{\max}C\sqrt{k/m}-a_i\mu_i)\right)\nonumber\\
	&\leq& a^{-1}\left(1-\bar{\Phi}(\tau_m+a_{\max}C\sqrt{k/m}-a_{\min}\mu)\right)=o(1),\nonumber
	\end{eqnarray}
	where the second to last inequality is due to (\ref{def:model.Z}) and (\ref{keyLemma.1}), and the last step holds from $a_{\max}\sqrt{k}\leq C\sqrt{m/\log m}$ implied by condition (C1').
	
	The remainder of the proof mirrors the arguments in Section \ref{sec:proof_7.5}, which demonstrate the sufficiency of the condition $h > \gamma$. Here, by analogous arguments, we establish the sufficiency of the condition $a^2_{\min} \cdot h > \gamma$. To avoid redundancy, we omit the replicated portions.
	
	\subsection{Proof of Lemma \ref{lemma:suf_05_L}} \label{sec:proof_7.10}
	
	The proof follows similar reasoning as in Lemma \ref{lemma:suf_nec_05}, with the following substitutions: replace $x_i$ with $\hat{z}_i$, $c^*_{m,0.5}(\alpha; F_{\mathbf{x}})$ with $c^*_{m, 0.5}(\alpha/2; F_{\widehat{\mathbf{z}}})$, $\mu$ with $a_{\min} \cdot \mu$, and substitute $\underline{\mu}_2$  with
	\[
	\underline{\mu}'_2=\sqrt{2[2\gamma-\tilde\eta]_+\log m+2\log\log m}.
	\]
	The proof parallels the arguments in Section \ref{sec:proof_7.6}, which demonstrate the sufficiency of the condition $h > [2\gamma -\eta]_+$. Here, using similar arguments, we establish the sufficiency of the condition $a^2_{\min} \cdot h > [2\gamma -\tilde \eta]_+$. To avoid redundancy, we omit the replicated sections.

	\subsection{Proof of Theorem \ref{thm:pi_hat_A}} \label{sec:proof_thm4.2}
	
	Because the same critical sequences, $c^*_{m,1}(\alpha/2; F_{{\mathbf{z}}})$ and $c^*_{m,0.5}(\alpha/2; F_{{\mathbf{z}}})$, as in Theorem \ref{thm:pi_hat_z}, are implemented in the construction of $\hat\pi^*_A(\widehat{\mathbf{z}})$,  it is sufficient, by Theorem \ref{thm:pi_hat_z}, to show that under the conditions of Lemma \ref{lemma:w_hat-w} with (C1) strengthened by (C1''), we have
	\begin{equation} \label{eq:pi_z_hat-pi_z}
	|\hat{\pi}^*_A(\widehat{\mathbf{z}}) - \hat{\pi}^*(\mathbf{z})| / \pi = o_p(1),
	\end{equation}
	which then leads to (\ref{eq:lowerbd_A}) and (\ref{eq:upperbd_A}). Note that in (C1''), the order of $\sqrt{m^{1-\tilde \eta}}$ is less than that of $\sqrt{m^{1-2(\tilde \eta - \gamma)_+}}$, which indicates $2\gamma>\tilde \eta$.

	To show (\ref{eq:pi_z_hat-pi_z}), recall that 
	\begin{equation} \label{def:pi*_hat_z_approx}
	\hat{\pi}^*_A(\widehat{\mathbf{z}})=\max\{\hat{\pi}^*_{A,0.5}(\widehat{\mathbf{z}}),~\hat{\pi}^*_{A,1}(\widehat{\mathbf{z}})\},
	\end{equation}
	and 
	\begin{equation}\label{def:pi*_hat_theta_z_approx}
	\hat{\pi}^*_{A,\theta}(\widehat{\mathbf{z}})=\max\limits_{t\in\mathbb{T}'}\dfrac{m^{-1}\sum\limits_{i=1}^{m}\mathbb{I}\{\hat z_i> t\}-\bar{\Phi}(t)-c^*_{m,\theta}(\alpha/2; F_{\mathbf{z}})[\bar{\Phi}(t)]^\theta}{1-\bar{\Phi}(t)},
	\end{equation}
	where $\mathbb{T}'=[\sqrt{2[\gamma- (\tilde \eta-\gamma)_+]_+\log m},\sqrt{5\log m}]\cap\mathbb{N}$. Note that $2\gamma>\tilde \eta$ implies that $[\gamma- (\tilde \eta-\gamma)_+]_+>0$.
	By the constructions of $\hat{\pi}^*_A(\widehat{\mathbf{z}})$ and $\hat{\pi}^*(\mathbf{z})$, it is sufficient to show 
	\begin{equation} \label{0.0}
	|\hat{\pi}^*_{A,0.5}(\widehat{\mathbf{z}}) - \hat{\pi}^*_{0.5}(\mathbf{z})| / \pi = o_p(1) \qquad \text{and} \qquad |\hat{\pi}^*_{A,1}(\widehat{\mathbf{z}}) - \hat{\pi}^*_{1}(\mathbf{z})| / \pi = o_p(1). 
	\end{equation}
	Define 
	\[
	B_{\theta}(t; \widehat{\mathbf{z}}) = \frac{m^{-1}\sum_{i=1}^{m}\mathbb{I}\{\hat z_i> t\}-\bar{\Phi}(t)-c^*_{m,\theta}(\alpha/2; F_{\mathbf{z}}) [\bar{\Phi}(t)]^\theta}{1-\bar{\Phi}(t)}
	\]
	and 
	\[
	B_{\theta}(t; \mathbf{z}) = \frac{m^{-1}\sum_{i=1}^{m}\mathbb{I}\{z_i> t\}-\bar{\Phi}(t)-c^*_{m,\theta}(\alpha/2; F_{\mathbf{z}}) [\bar{\Phi}(t)]^\theta}{1-\bar{\Phi}(t)}. 
	\]
	Then 
	\begin{equation} \label{2.1}
	\left|\hat{\pi}^*_{A,\theta}(\widehat{\mathbf{z}}) - \hat{\pi}^*_{\theta}(\mathbf{z})\right| = \left|\max_{t \in \mathbb{T}'} B_{\theta}(t; \widehat{\mathbf{z}}) - \max_{t \in \mathbb{T}} B_{\theta}(t; \mathbf{z})\right|.
	\end{equation}
	
	First, we show that 
	\begin{equation} \label{2.2}
	\left|\max_{t \in \mathbb{T}'} B_{\theta}(t; \widehat{\mathbf{z}}) - \max_{t \in \mathbb{T}} B_{\theta}(t; \mathbf{z})\right| \le \max_{t \in \mathbb{T}'} \left|B_{\theta}(t; \widehat{\mathbf{z}}) - B_{\theta}(t; \mathbf{z})\right|, 
	\end{equation}
	which is implied by 
	\begin{equation} \label{0.1}
	\max_{t \in \mathbb{T}'} B_{\theta}(t; \widehat{\mathbf{z}}) - \max_{t \in \mathbb{T}} B_{\theta}(t; \mathbf{z}) \le \max_{t \in \mathbb{T}'} \left|B_{\theta}(t; \widehat{\mathbf{z}}) - B_{\theta}(t; \mathbf{z})\right|
	\end{equation} 
	and 
	\begin{equation} \label{0.2}
	\max_{t \in \mathbb{T}'} B_{\theta}(t; \widehat{\mathbf{z}}) - \max_{t \in \mathbb{T}} B_{\theta}(t; \mathbf{z}) \ge - \max_{t \in \mathbb{T}'} \left|B_{\theta}(t; \widehat{\mathbf{z}}) - B_{\theta}(t; \mathbf{z})\right|.
	\end{equation}
	Equation (\ref{0.1}) holds because 
	\begin{eqnarray*} 
		\max_{t \in \mathbb{T}'} B_{\theta}(t; \widehat{\mathbf{z}}) - \max_{t \in \mathbb{T}} B_{\theta}(t; \mathbf{z}) & = & B_{\theta}(t_1; \widehat{\mathbf{z}}) - \max_{t \in \mathbb{T}} B_{\theta}(t; \mathbf{z}) \le B_{\theta}(t_1; \widehat{\mathbf{z}}) - B_{\theta}(t_1; \mathbf{z}) \\
		& \le & \left|  B_{\theta}(t_1; \widehat{\mathbf{z}}) - B_{\theta}(t_1; \mathbf{z}) \right|  \le  \max_{t \in \mathbb{T}'} \left|B_{\theta}(t; \widehat{\mathbf{z}}) - B_{\theta}(t; \mathbf{z})\right|.
	\end{eqnarray*}
	Equation (\ref{0.2}) also holds because 
	\begin{eqnarray*} 
		\max_{t \in \mathbb{T}'} B_{\theta}(t; \widehat{\mathbf{z}}) - \max_{t \in \mathbb{T}} B_{\theta}(t; \mathbf{z}) & = & \max_{t \in \mathbb{T}'} B_{\theta}(t; \widehat{\mathbf{z}}) - B_{\theta}(t_2; \mathbf{z}) \ge B_{\theta}(t_2; \widehat{\mathbf{z}}) - B_{\theta}(t_2; \mathbf{z}) \\
		& \ge & -\left|  B_{\theta}(t_2; \widehat{\mathbf{z}}) - B_{\theta}(t_2; \mathbf{z}) \right|  \ge  -\max_{t \in \mathbb{T}'} \left|B_{\theta}(t; \widehat{\mathbf{z}}) - B_{\theta}(t; \mathbf{z})\right|.
	\end{eqnarray*}
	Then, by (\ref{2.1}) and (\ref{2.2}),  it is sufficient to show 
	\begin{equation*} 
	\max_{t \in \mathbb{T}'} \left|B_{\theta}(t; \widehat{\mathbf{z}}) - B_{\theta}(t; \mathbf{z})\right|/\pi = o_p(1).
	\end{equation*}
	
	On the other hand, because
	\begin{eqnarray*}
		\text{P}(\max_{t \in \mathbb{T}'} \left|B_{\theta}(t; \widehat{\mathbf{z}}) - B_{\theta}(t; \mathbf{z})\right|/\pi > a) & \le & \text{P}(\max_{t \in \mathbb{T}'} \left|B_{\theta}(t; \widehat{\mathbf{z}}) - B_{\theta}(t; \mathbf{z})\right|/\pi > a, \|\widehat{\mathbf{w}} - \mathbf{w}\|_2\leq C\sqrt{k/m}) \\
		& + & \text{P}(\|\widehat{\mathbf{w}} - \mathbf{w}\|_2 > C\sqrt{k/m}) 
	\end{eqnarray*}
	for any $a>0$, and 
	\[
	\text{P}(\|\widehat{\mathbf{w}} - \mathbf{w}\|_2 > C\sqrt{k/m}) = o(1)
	\]
	under the conditions of Lemma \ref{lemma:w_hat-w} with (C1) strengthened by (C1''), then it is sufficient to show that 
	\begin{equation} \label{0.3} 	
	\text{P}(\max_{t \in \mathbb{T}'} \left|B_{\theta}(t; \widehat{\mathbf{z}}) - B_{\theta}(t; \mathbf{z})\right|/\pi > a, \|\widehat{\mathbf{w}} - \mathbf{w}\|_2\leq C\sqrt{k/m}) =o(1).
	\end{equation}
	
	We provide the following lemma. 
	\begin{lemma} \label{lemma:keyLemma_s}
		Given condition (C1''), for any $t = t_m = \sqrt{2 h' \log m}$ with $h' \geq [\gamma-(\tilde \eta - \gamma)_+]_+$ and any $a>0$, we have
		\begin{equation} \label{0.5}
		\text{P} \left(s^{-1} \left| \sum_{i=1}^{m}\mathbb{I}\{\hat z_i> t\} - \sum_{i=1}^{m}\mathbb{I}\{z_i> t\} \right| > a \mid \|\widehat{\mathbf{w}} - \mathbf{w}\|_2\leq C\sqrt{k/m} \right)=o(1/\sqrt{\log m}).
		\end{equation}
		
	\end{lemma}
	
	Then, for any $a>0$, 
	\begin{eqnarray*} 
		& & \text{P}(\max_{t \in \mathbb{T}'} \left|B_{\theta}(t; \widehat{\mathbf{z}}) - B_{\theta}(t; \mathbf{z})\right|/\pi > a, \|\widehat{\mathbf{w}} - \mathbf{w}\|_2\leq C\sqrt{k/m}) \\
		& \le & \text{P}(\max_{t \in \mathbb{T}'} \left|B_{\theta}(t; \widehat{\mathbf{z}}) - B_{\theta}(t; \mathbf{z})\right|/\pi > a \mid \|\widehat{\mathbf{w}} - \mathbf{w}\|_2\leq C\sqrt{k/m}) \\	
		& \le & \sum_{t \in \mathbb{T}'} \text{P} ( \left|B_{\theta}(t; \widehat{\mathbf{z}}) - B_{\theta}(t; \mathbf{z})\right|/\pi > a \mid \|\widehat{\mathbf{w}} - \mathbf{w}\|_2\leq C\sqrt{k/m}) \\
		& \le & \sum_{t \in \mathbb{T}'} \text{P} \left( C s^{-1}\left| \sum_{i=1}^{m}\mathbb{I}\{\hat z_i> t\} - \sum_{i=1}^{m}\mathbb{I}\{z_i> t\} \right| > a \mid \|\widehat{\mathbf{w}} - \mathbf{w}\|_2\leq C\sqrt{k/m} \right) \\
		& \le & C \cdot \sqrt{\log m} \cdot o(1/\sqrt{\log m}) = o(1),
	\end{eqnarray*}
	where the third inequality is by Lemma \ref{lemma:keyLemma_s}. Therefore (\ref{0.3}) is proved. This concludes the proof of Theorem \ref{thm:pi_hat_A}.
	
	\subsection{Proof of Lemma \ref{lemma:keyLemma_s}} \label{sec:proof_7.11}
	
	Note that
	\begin{eqnarray*}
		\left|\sum_{i=1}^m\mathbb{I}\{\hat z_i> t\} - \sum_{i=1}^m\mathbb{I}\{z_i> t\} \right| & \le & \sum_{i=1}^m \left|\mathbb{I}\{\hat z_i> t\} - \mathbb{I}{\{z_i> t\}}\right|  =  \sum_{i=1}^m\mathbb{I}\{z_i \in (t, t+a_i \mathbf{b}_i^T(\widehat{\mathbf{w}} - \mathbf{w}))\} \\
		& \le & \sum_{i=1}^m\mathbb{I}\{z^0_i \in (t, t+a_i \mathbf{b}_i^T(\widehat{\mathbf{w}} - \mathbf{w}))\} \\
		& + & \sum_{i\in I_1}\mathbb{I}\{z_i \in (t, t+a_i \mathbf{b}_i^T(\widehat{\mathbf{w}} - \mathbf{w}))\}
	\end{eqnarray*}
	almost surely. So it is sufficient to show
	\begin{equation} \label{2.3}
	\text{P}\left(s^{-1}\sum_{i=1}^m\mathbb{I}\{z_i^0 \in (t, t+a_i \mathbf{b}_i^T(\widehat{\mathbf{w}} - \mathbf{w}))\}>a/2  \mid \|\widehat{\mathbf{w}} - \mathbf{w}\|_2\leq C\sqrt{k/m} \right) =   o\left({1\over \sqrt{\log m}}\right), 
	\end{equation}
	and 
	\begin{equation} \label{2.4}
	\text{P}\left(s^{-1}\sum_{i\in I_1}\mathbb{I}\{z_i \in (t, t+a_i \mathbf{b}_i^T(\widehat{\mathbf{w}} - \mathbf{w}))\}>a/2  \mid \|\widehat{\mathbf{w}} - \mathbf{w}\|_2\leq C\sqrt{k/m} \right) =   o\left({1\over \sqrt{\log m}}\right)
	\end{equation}
	
	Consider (\ref{2.3}) first. Markov's inequality gives
	\begin{eqnarray*}\label{0.5.1}
		&&\text{P}\left(s^{-1}\sum_{i=1}^m\mathbb{I}\{z_i^0 \in (t, t+a_i \mathbf{b}_i^T(\widehat{\mathbf{w}} - \mathbf{w}))\}>a/2 \mid \|\widehat{\mathbf{w}} - \mathbf{w}\|_2\leq C\sqrt{k/m} \right)\nonumber\\
		&\leq& \dfrac{2}{as}\sum_{i=1}^{m}\text{P}\left(z_i^0 \in (t, t+a_i \mathbf{b}_i^T(\widehat{\mathbf{w}} - \mathbf{w})) \mid \|\widehat{\mathbf{w}} - \mathbf{w}\|_2\leq C\sqrt{k/m} \right)\nonumber\\
		&\leq&\dfrac{2}{as}\sum_{i=1}^{m}\text{P}\left(z_i^0 \in (t, t-a_{\max}C\sqrt{k/m})\right)\nonumber\\
		&\leq&\dfrac{Cm(a_{\max}\sqrt{k/m})\phi(t-a_{\max}C\sqrt{k/m})}{s}\leq\dfrac{Cm^{1/2}(a_{\max}\sqrt{k})e^{-t^2/2}e^{Cta_{\max}\sqrt{k/m}}}{s} \nonumber\\
		&\leq&\dfrac{C m^{1/2} (a_{\max}\sqrt{k})e^{-t^2/2}}{s} \leq \dfrac{Cm^{-(h'-1/2)} a_{\max}\sqrt{k}}{m^{1-\gamma}} = o\left({1\over \sqrt{\log m}}\right),
	\end{eqnarray*}	
	where the second inequality holds by (\ref{keyLemma.1}), the fifth inequality holds by $a_{\max}\sqrt{k}=o(\sqrt{m/\log m})$, the sixth inequality is by $t=\sqrt{2h'\log m}$, and the last step is by $h'\geq [\gamma-(\tilde \eta-\gamma)_+]_+$ and  $a_{\max}\sqrt{k} = o(\sqrt{m^{1-2(\tilde \eta - \gamma)_+}/\log m})$. 
	Therefore, (\ref{2.3}) is verified. 
	
	Next, consider (\ref{2.4}). Markov's inequality gives
	\begin{eqnarray}\label{0.5.4}
	&&\text{P}\left(s^{-1}\sum_{i\in I_1}\mathbb{I}\{z_i \in (t, t+a_i \mathbf{b}_i^T(\widehat{\mathbf{w}} - \mathbf{w}))\}>a/2 \mid \|\widehat{\mathbf{w}} - \mathbf{w}\|_2\leq C\sqrt{k/m} \right)\nonumber\\
	&\leq& \dfrac{2}{as}\sum_{i\in I_1}\text{P}\left(z_i \in (t, t+a_i \mathbf{b}_i^T(\widehat{\mathbf{w}} - \mathbf{w})) \mid \|\widehat{\mathbf{w}} - \mathbf{w}\|_2\leq C\sqrt{k/m}\right)\nonumber\\
	&=& \dfrac{2}{as}\sum_{i\in I_1}\text{P}\left(z_i^0\in(t-a_i\mu_i, t-a_i\mu_i +a_i \mathbf{b}_i^T(\widehat{\mathbf{w}} - \mathbf{w})) \mid \|\widehat{\mathbf{w}} - \mathbf{w}\|_2\leq C\sqrt{k/m} \right)\nonumber\\
	&\le & \dfrac{2}{as}\sum_{i\in I_1}\text{P}\left(z_i^0\in(t-a_i\mu_i, t-a_i\mu_i - a_{\max} C \sqrt{k/m}) \right)\nonumber\\
	&\leq&\dfrac{2}{as}\sum_{i\in I_1}\left(a_{\max}C\sqrt{k/m}\phi(0)\right) = Ca_{\max}\sqrt{k/m}=o\left({1\over \sqrt{\log m}}\right),
	\end{eqnarray}
	where the third step is by (\ref{keyLemma.1}), and the last step is by condition (C1''). Therefore (\ref{2.4}) holds. This concludes the proof.

	\subsection{Proof of Lemma \ref{lemma:w_hat-w}} \label{sec:proof_lemma4.1}
	
	Without loss of generality, we assume that the true value $\mathbf{w}$ is $\mathbf{0}$ and prove that $\|\widehat{\mathbf{w}}\|_2=O_p(\sqrt{k/m})$.
	
	Define $L:\mathbb{R}^k\rightarrow \mathbb{R}$ such that
	\[
	L(\bm{\beta})=m^{-1}\sum_{i=1}^{m}|\mu_i+K_i-\mathbf{b}_i^T\bm{\beta}|.
	\]
	Define $l:\mathbb{R}^k\rightarrow \mathbb{R}^k$ such that for $j=1,\cdots,k$,
	\[
	l_j(\bm{\beta})=m^{-1}\sum_{i=1}^{m}b_{ij}\text{sgn}(\mu_i+K_i-\mathbf{b}_i^T\bm{\beta}).
	\]
	Under the assumption that $\mathbf{w}=\mathbf{0}$, we have $\mathbf{x}=\bm{\mu}+\mathbf{K}$, so
	\[
	\widehat{\mathbf{w}}=\arg\min\limits_{\bm{\beta}\in\mathbb{R}^k}\sum_{i=1}^{m}|x_i-\mathbf{b}_i^T\bm{\beta}|=\arg\min\limits_{\bm{\beta}\in\mathbb{R}^k}L(\bm{\beta}).
	\] 
	Given that $L$ is convex in $\bm{\beta}$, and $-l(\bm{\beta})\in\nabla L(\bm{\beta})$, by classical convexity argument, it suffices to show that with high probability, $\bm{\beta}^Tl(\bm{\beta})<0$ with $\|\bm{\beta}\|_2=B\sqrt{k/m}$ for a sufficiently large constant $B$.
	
	In fact, since $L$ is convex in $\bm{\beta}$, and
	\[
	-l(\bm{\beta})\in\nabla L(\bm{\beta}),
	\]
	we have for $\forall\theta\in(0,1)$,
	\[
	L(\widehat{\mathbf{w}})\geq L(\theta\widehat{\mathbf{w}})-(\widehat{\mathbf{w}}-\theta\widehat{\mathbf{w}})^Tl(\theta\widehat{\mathbf{w}}).
	\]
	Since $L(\theta\widehat{\mathbf{w}})\geq L(\widehat{\mathbf{w}})$, we have $\theta\widehat{\mathbf{w}}^Tl(\theta\widehat{\mathbf{w}})\geq0$. If $\|\widehat{\mathbf{w}}\|_2>B\sqrt{k/m}$, then $\exists\theta_0\in(0,1)$, s.t., $\|\theta_0\widehat{\mathbf{w}}\|_2=B\sqrt{k/m}$, thus $\theta_0\widehat{\mathbf{w}}^Tl(\theta_0\widehat{\mathbf{w}})<0$, which contradicts the fact that $\theta_0\widehat{\mathbf{w}}^Tl(\theta_0\widehat{\mathbf{w}})\geq0$. Therefore, we have $\|\widehat{\mathbf{w}}\|_2\leq B\sqrt{k/m}$ with high probability.
	
	Here
	\begin{eqnarray}
		\bm{\beta}^Tl(\bm{\beta})&=&m^{-1}\sum_{i=1}^m(\mathbf{b}_i^T\bm{\beta})\text{sgn}(\mu_i+K_i-\mathbf{b}_i^T\bm{\beta})\nonumber\\
		&=&m^{-1}\sum_{i\in I_1}(\mathbf{b}_i^T\bm{\beta})\text{sgn}(\mu_i+K_i-\mathbf{b}_i^T\bm{\beta})+m^{-1}\sum_{i\in I_0}(\mathbf{b}_i^T\bm{\beta})\text{sgn}(K_i-\mathbf{b}_i^T\bm{\beta})\nonumber\\
		&\leq&m^{-1}\sum_{i\in I_1}|\mathbf{b}_i^T\bm{\beta}|+m^{-1}\sum_{i=1}^m(\mathbf{b}_i^T\bm{\beta})\text{sgn}(K_i-\mathbf{b}_i^T\bm{\beta}).\nonumber
	\end{eqnarray}
	Let
	\[
	V=m^{-1}\sum_{i=1}^mV_i,
	\]
	where
	$V_i=(\mathbf{b}_i^T\bm{\beta})\text{sgn}(K_i-\mathbf{b}_i^T\bm{\beta})$. Also let $\bm{\beta}=\omega\mathbf{u}$ with $\|\mathbf{u}\|_2=1$, so $\omega=B\sqrt{k/m}$. Since
	\[
	m^{-1}\sum_{i\in I_1}|\mathbf{b}_i^T\bm{\beta}|\leq m^{-1}\omega\sum_{i\in I_1}|\mathbf{b}_i^T\mathbf{u}|\leq m^{-1}\omega\sum_{i\in I_1}\|\mathbf{b}_i\|_2^2\leq m^{-1}\omega|I_1|=\pi\omega=m^{-\gamma}\omega,
	\]
	we have
	\[
	\bm{\beta}^Tl(\bm{\beta})\leq m^{-\gamma}\omega+V.
	\]
	For $\forall h>0$, by Chebyshev's inequality,
	\[
	\text{P}\left(V<\text{E}(V)+h\text{SD}(V)\right)>1-h^{-2}.
	\]
	To show that $\bm{\beta}^Tl(\bm{\beta})<0$ with high probability, it is sufficient to show that $\exists B$ and $\exists M$, s.t., $\forall m>M$, $\text{P}(m^{-\gamma}\omega+V<0)>1-h^{-2}$, and it suffices to show that when $B$ is large enough, $m^{-\gamma}\omega+\text{E}(V)+h\text{SD}(V)<0$, i.e., $-\text{E}(V)>m^{-\gamma}\omega+h\text{SD}(V)$, as $m\rightarrow\infty$. 
	
	In Section \ref{section_variance}, we show that
	\[
	\text{SD}(V)=O\left(\omega \sqrt{m^{-\tilde \eta}}\exp\left\{-\dfrac{1}{\pi}(a_{\min}d\omega)^2\right\}\right),
	\]
	and in Section \ref{section_expectation}, we show that
	\[
	-\text{E}(V)\geq\sqrt{\dfrac{2}{9\pi}}a_{\min}d^2\omega^2\exp\left\{-\dfrac{1}{2\pi}(a_{\min}d\omega)^2\right\}.
	\]
	Then to show that $-\text{E}(V)>m^{-\gamma}\omega+h\text{SD}(V)$, it is sufficient to show that
	\[
	\sqrt{\dfrac{2}{9\pi}}a_{\min}d^2\omega^2\exp\left\{-\dfrac{1}{2\pi}(a_{\min}d\omega)^2\right\}>m^{-\gamma}\omega+Ch\omega \sqrt{m^{-\tilde \eta}}\exp\left\{-\dfrac{1}{\pi}(a_{\min}d\omega)^2\right\},
	\]
	i.e.,
	\[
	a_{\min}\omega>C\left[m^{-\gamma}\exp\left\{\dfrac{1}{2\pi}(a_{\min}d\omega)^2\right\}+h\sqrt{m^{-\tilde \eta}}\exp\left\{-\dfrac{1}{2\pi}(a_{\min}d\omega)^2\right\}\right].
	\]
	Here $\exp\left\{-\dfrac{1}{2\pi}(a_{\min}d\omega)^2\right\}\leq1$. Recall that $\omega=B\sqrt{k/m}$. Under condition (C1), we have $a_{\max}\sqrt{k/m}=o(1)$, thus  $1\leq\exp\left\{\dfrac{1}{2\pi}(a_{\min}d\omega)^2\right\}\leq\exp\left\{\dfrac{1}{2\pi}(a_{\max}d\omega)^2\right\}\to1$. Therefore, the conclusion holds when
	\[
	B>C(a_{\min}\sqrt{k/m})^{-1}(m^{-\gamma}+h\sqrt{m^{-\tilde \eta}}).
	\]
	This can be achieved when $a_{\min}\sqrt{k/m}\geq C\sqrt{m^{-\min\{\tilde \eta,~2\gamma\}}}$, which is given by condition (C1).
	
	The remaining of the proof is similar to the proof of Theorem 3 in \cite{fan2012estimating}, where Polya's approximation is used, i.e., for $x>0$, with $\sup\limits_{x>0}|\delta(x)|<0.004$, we have
	\[
	\Phi(x)=\dfrac{1}{2}\left\{1+\sqrt{1-\exp\left(-\dfrac{2}{\pi}x^2\right)}\right\}(1+\delta(x)).
	\]

	\subsubsection{Bounding the variance}\label{section_variance}
	
	The variance of $V$ is
	\begin{eqnarray}\label{appendix.var_V}
	\text{Var}(V)&=&m^{-2}\sum_{i=1}^{m}\text{Var}(V_i)+m^{-2}\sum_{i\neq j}\text{Cov}(V_i,V_j)\nonumber\\
	&\leq& m^{-2}\sum_{i=1}^{m}\text{Var}(V_i)+m^{-2}\sum_{i\neq j}|\text{Cov}(V_i,V_j)|,
	\end{eqnarray}
	where
	\begin{equation}\label{appendix.var_V.1}
	\sum_{i=1}^{m}\text{Var}(V_i)=\sum_{i=1}^{m}\mathbb{I}\{|\mathbf{b}_i^T\mathbf{u}|\leq d\}\text{Var}(V_i)+\sum_{i=1}^{m}\mathbb{I}\{|\mathbf{b}_i^T\mathbf{u}|>d\}\text{Var}(V_i)
	\end{equation}
	and
	\begin{eqnarray}\label{appendix.var_V.2}
	\sum_{i\neq j}|\text{Cov}(V_i,V_j)|&=&\sum_{i\neq j}\mathbb{I}\{|\mathbf{b}_i^T\mathbf{u}|\leq d\}\mathbb{I}\{|\mathbf{b}_j^T\mathbf{u}|\leq d\}|\text{Cov}(V_i,V_j)|\nonumber\\
	&+&2\sum_{i\neq j}\mathbb{I}\{|\mathbf{b}_i^T\mathbf{u}|\leq d\}\mathbb{I}\{|\mathbf{b}_j^T\mathbf{u}|> d\}|\text{Cov}(V_i,V_j)|\nonumber\\
	&+&\sum_{i\neq j}\mathbb{I}\{|\mathbf{b}_i^T\mathbf{u}|> d\}\mathbb{I}\{|\mathbf{b}_j^T\mathbf{u}|> d\}|\text{Cov}(V_i,V_j)|.
	\end{eqnarray}
	
	As for $\text{Var}(V_i)$, by Polya's approximation, we have
	\begin{eqnarray}
	\text{Var}(V_i)&=&4\text{Var}\left(\dfrac{1-V_i}{2}\right)=4(\mathbf{b}_i^T\bm{\beta})^2\text{Var}(\text{Bern}(\text{P}(K_i<\mathbf{b}_i^T\bm{\beta})))\nonumber\\
	&=&4(\mathbf{b}_i^T\bm{\beta})^2\Phi(a_i\mathbf{b}_i^T\bm{\beta})(1-\Phi(a_i\mathbf{b}_i^T\bm{\beta}))\nonumber\\
	&=&(\mathbf{b}_i^T\bm{\beta})^2\exp\left\{-\dfrac{2}{\pi}(a_i\mathbf{b}_i^T\bm{\beta})^2\right\}(1+\delta_i)^2\nonumber
	\end{eqnarray}
	with $|\delta_i|<0.004$. On the one hand,
	\begin{eqnarray}\label{appendix.var_Vi.1}
	\sum_{i=1}^{m}\mathbb{I}\{|\mathbf{b}_i^T\mathbf{u}|>d\}\text{Var}(V_i)&=&\sum_{i=1}^{m}\mathbb{I}\{|\mathbf{b}_i^T\mathbf{u}|>d\}(\mathbf{b}_i^T\mathbf{u})^2\omega^2\exp\left\{-\dfrac{2}{\pi}(a_i\mathbf{b}_i^T\mathbf{u}\omega)^2\right\}(1+\delta_i)^2\nonumber\\
	&\leq&\sum_{i=1}^{m}\omega^2\exp\left\{-\dfrac{2}{\pi}(a_id\omega)^2\right\}(1+\delta_i)^2\nonumber\\
	&\leq&2m\omega^2\exp\left\{-\dfrac{2}{\pi}(a_{\min}d\omega)^2\right\},
	\end{eqnarray}
	where the second step holds since $\mathbf{b}_i^T\mathbf{u}\leq\|\mathbf{b}_i\|_2^2\leq1$. On the other hand,
	\begin{eqnarray}
	\sum_{i=1}^{m}\mathbb{I}\{|\mathbf{b}_i^T\mathbf{u}|\leq d\}\text{Var}(V_i)&=&\sum_{i=1}^{m}\mathbb{I}\{|\mathbf{b}_i^T\mathbf{u}|\leq d\}(\mathbf{b}_i^T\mathbf{u})^2\omega^2\exp\left\{-\dfrac{2}{\pi}(a_i\mathbf{b}_i^T\mathbf{u}\omega)^2\right\}(1+\delta_i)^2\nonumber\\
	&\leq&2\sum_{i=1}^{m}\mathbb{I}\{|\mathbf{b}_i^T\mathbf{u}|\leq d\}d^2\omega^2\nonumber
	\end{eqnarray}
	and
	\begin{eqnarray}
	\sum_{i=1}^{m}\mathbb{I}\{|\mathbf{b}_i^T\mathbf{u}|>d\}\text{Var}(V_i)&=&\sum_{i=1}^{m}\mathbb{I}\{|\mathbf{b}_i^T\mathbf{u}|>d\}(\mathbf{b}_i^T\mathbf{u})^2\omega^2\exp\left\{-\dfrac{2}{\pi}(a_i\mathbf{b}_i^T\mathbf{u}\omega)^2\right\}(1+\delta_i)^2\nonumber\\
	&\geq&\dfrac{1}{2}\sum_{i=1}^{m}\mathbb{I}\{|\mathbf{b}_i^T\mathbf{u}|>d\}d^2\omega^2\exp\left\{-\dfrac{2}{\pi}(a_{\max}\omega)^2\right\},\nonumber
	\end{eqnarray}
	where the last step holds since $\mathbf{b}_i^T\mathbf{u}\leq\|\mathbf{b}_i\|_2^2\leq1$. This implies that
	\begin{equation}
	\dfrac{\sum\limits_{i=1}^{m}\mathbb{I}\{|\mathbf{b}_i^T\mathbf{u}|\leq d\}\text{Var}(V_i)}{\sum\limits_{i=1}^{m}\mathbb{I}\{|\mathbf{b}_i^T\mathbf{u}|>d\}\text{Var}(V_i)}\leq\dfrac{4m^{-1}\sum\limits_{i=1}^{m}\mathbb{I}\{|\mathbf{b}_i^T\mathbf{u}|\leq d\}}{1-m^{-1}\sum\limits_{i=1}^{m}\mathbb{I}\{|\mathbf{b}_i^T\mathbf{u}|\leq d\}}\exp\left\{\dfrac{2}{\pi}(a_{\max}\omega)^2\right\}.\nonumber
	\end{equation}
	Under condition (C1), we have $a_{\max}\sqrt{k/m}=o(1)$, thus  $\exp\left\{\dfrac{2}{\pi}(a_{\max}\omega)^2\right\}=O(1)$. Then condition (C2) gives
	\begin{equation}\label{appendix.var_Vi.2}
	\sum\limits_{i=1}^{m}\mathbb{I}\{|\mathbf{b}_i^T\mathbf{u}|\leq d\}\text{Var}(V_i)=o\left(\sum\limits_{i=1}^{m}\mathbb{I}\{|\mathbf{b}_i^T\mathbf{u}|>d\}\text{Var}(V_i)\right).
	\end{equation}
	Combining (\ref{appendix.var_V.1}), (\ref{appendix.var_Vi.1}), and (\ref{appendix.var_Vi.2}) gives
	\begin{equation}\label{appendix.var_Vi}
	\sum_{i=1}^{m}\text{Var}(V_i)=O\left(m\omega^2\exp\left\{-\dfrac{2}{\pi}(a_{\min}d\omega)^2\right\}\right).
	\end{equation}
	
	As for $\text{Cov}(V_i,V_j)$, we have
	\begin{eqnarray}
	\text{Cov}(V_i,V_j)&=&4\text{Cov}\left(\dfrac{1-V_i}{2},\dfrac{1-V_j}{2}\right)\nonumber\\
	&=&4(\mathbf{b}_i^T\bm{\beta})(\mathbf{b}_j^T\bm{\beta})\{\text{P}(K_i<\mathbf{b}_i^T\bm{\beta},K_j<\mathbf{b}_j^T\bm{\beta})-\text{P}(K_i<\mathbf{b}_i^T\bm{\beta})\text{P}(K_j<\mathbf{b}_j^T\bm{\beta})\}\nonumber\\
	&=&4(\mathbf{b}_i^T\bm{\beta})(\mathbf{b}_j^T\bm{\beta})\{\Phi(a_i\mathbf{b}_i^T\bm{\beta},a_j\mathbf{b}_j^T\bm{\beta};\tilde\Sigma_{ij})-\Phi(a_i\mathbf{b}_i^T\bm{\beta})\Phi(a_j\mathbf{b}_j^T\bm{\beta})\}.\nonumber\\
	&=&4(\mathbf{b}_i^T\bm{\beta})(\mathbf{b}_j^T\bm{\beta})\phi(a_i\mathbf{b}_i^T\bm{\beta})\phi(a_j\mathbf{b}_j^T\bm{\beta})\tilde\Sigma_{ij}(1+o(1)),\nonumber
	\end{eqnarray}
	where $\tilde\Sigma_{ij}$ is the correlation between $K_i$ and $K_j$.  On the one hand,
	\begin{eqnarray}\label{appendix.cov_Vi.1}
	&&\sum_{i\neq j}\mathbb{I}\{|\mathbf{b}_i^T\mathbf{u}|> d\}\mathbb{I}\{|\mathbf{b}_j^T\mathbf{u}|> d\}|\text{Cov}(V_i,V_j)|\nonumber\\
	&=&4\sum_{i\neq j}\mathbb{I}\{|\mathbf{b}_i^T\mathbf{u}|> d\}\mathbb{I}\{|\mathbf{b}_j^T\mathbf{u}|> d\}\left|(\mathbf{b}_i^T\mathbf{u})(\mathbf{b}_j^T\mathbf{u})\omega^2\phi(a_i\mathbf{b}_i^T\mathbf{u}\omega)\phi(a_j\mathbf{b}_j^T\mathbf{u}\omega)\tilde\Sigma_{ij}(1+o(1))\right|\nonumber\\
	&\leq&4\sum_{i\neq j}\omega^2\phi(a_id\omega)\phi(a_jd\omega)|\tilde\Sigma_{ij}|(1+o(1))\nonumber\\
	&\leq&C\sum_{i\neq j}\omega^2\exp\left\{-(a_{\min}d\omega)^2\right\}|\tilde\Sigma_{ij}|,\nonumber\\
	\end{eqnarray}
	where the second step holds since $\mathbf{b}_i^T\mathbf{u}\leq\|\mathbf{b}_i\|_2^2\leq1$. On the other hand, similar arguments to those used in proving (\ref{appendix.var_Vi.2}) give
	\begin{equation}\label{appendix.cov_Vi.2}
	\sum_{i\neq j}\mathbb{I}\{|\mathbf{b}_i^T\mathbf{u}|\leq d\}\mathbb{I}\{|\mathbf{b}_j^T\mathbf{u}|> d\}|\text{Cov}(V_i,V_j)|=o\left(\sum_{i\neq j}\mathbb{I}\{|\mathbf{b}_i^T\mathbf{u}|> d\}\mathbb{I}\{|\mathbf{b}_j^T\mathbf{u}|> d\}|\text{Cov}(V_i,V_j)|\right),
	\end{equation}
	and
	\begin{equation}\label{appendix.cov_Vi.3}
	\sum_{i\neq j}\mathbb{I}\{|\mathbf{b}_i^T\mathbf{u}|\leq d\}\mathbb{I}\{|\mathbf{b}_j^T\mathbf{u}|\leq d\}|\text{Cov}(V_i,V_j)|=o\left(\sum_{i\neq j}\mathbb{I}\{|\mathbf{b}_i^T\mathbf{u}|> d\}\mathbb{I}\{|\mathbf{b}_j^T\mathbf{u}|> d\}|\text{Cov}(V_i,V_j)|\right).
	\end{equation}
	Combining (\ref{appendix.var_V.2}), (\ref{appendix.cov_Vi.1}), (\ref{appendix.cov_Vi.2}), and (\ref{appendix.cov_Vi.3}) gives
	\begin{equation}\label{appendix.cov_Vi}
	\sum_{i\neq j}|\text{Cov}(V_i,V_j)|=O\left(\omega^2\exp\left\{-(a_{\min}d\omega)^2\right\}\sum_{i\neq j}|\tilde\Sigma_{ij}|\right).
	\end{equation}
	
	Combining (\ref{appendix.var_V}), (\ref{appendix.var_Vi}), and (\ref{appendix.cov_Vi}), we have
	\[
	\text{Var}(V)=O\left(\omega^2m^{-\tilde \eta}\exp\left\{-\dfrac{2}{\pi}(a_{\min}d\omega)^2\right\}\right),
	\]
	so the standard deviation of $V$ is bounded by
	\begin{equation}\label{appendix.SD_V}
	\text{SD}(V)=O\left(\omega \sqrt{m^{-\tilde \eta}}\exp\left\{-\dfrac{1}{\pi}(a_{\min}d\omega)^2\right\}\right).
	\end{equation}

	\subsubsection{Bounding the expectation}\label{section_expectation}
	
	As for the expectation of $V$, we have
	\begin{eqnarray}
	-\text{E}(V)&=&m^{-1}\sum_{i=1}^m\mathbf{b}_i^T\bm{\beta}\text{E}\left(-\text{sgn}(K_i-\mathbf{b}_i^T\bm{\beta})\right)=m^{-1}\sum_{i=1}^m\mathbf{b}_i^T\bm{\beta}\left(2\text{P}(K_i<\mathbf{b}_i^T\bm{\beta})-1\right)\nonumber\nonumber\\
	&=&\dfrac{2}{m}\sum_{i=1}^m\mathbf{b}_i^T\bm{\beta}\left(\Phi(a_i\mathbf{b}_i^T\bm{\beta})-\dfrac{1}{2}\right)\geq\dfrac{2}{m}\sum_{i=1}^m\mathbb{I}\{|\mathbf{b}_i^T\mathbf{u}|> d\}\mathbf{b}_i^T\mathbf{u}\omega\left(\Phi(a_i\mathbf{b}_i^T\mathbf{u}\omega)-\dfrac{1}{2}\right)\nonumber\\
	&\geq&\dfrac{2}{m}\sum_{i=1}^m\mathbb{I}\{|\mathbf{b}_i^T\mathbf{u}|> d\}d\omega\left(\Phi(a_id\omega)-\dfrac{1}{2}\right)\nonumber\\
	&\geq&\dfrac{2d\omega}{m}\sum_{i=1}^m\left(\Phi(a_id\omega)-\dfrac{1}{2}\right)-\dfrac{2d\omega}{m}\sum_{i=1}^m\mathbb{I}\{|\mathbf{b}_i^T\mathbf{u}|\leq d\}\left(\Phi(a_id\omega)-\dfrac{1}{2}\right),\nonumber
	\end{eqnarray}
	where the first inequality holds due to the fact that $x(\Phi(x)-1/2)\geq0$ for $x\in\mathbb{R}$. Under condition (C2), we have
	\[
	(0\leq)\dfrac{2d\omega}{m}\sum_{i=1}^m\mathbb{I}\{|\mathbf{b}_i^T\mathbf{u}|\leq d\}\left(\Phi(a_id\omega)-\dfrac{1}{2}\right)\leq\dfrac{d\omega}{m}\sum_{i=1}^m\mathbb{I}\{|\mathbf{b}_i^T\mathbf{u}|\leq d\}\rightarrow0
	\]
	as $m\rightarrow0$. Thus for sufficiently large $m$, we have
	\[
	-\text{E}(V)\geq\dfrac{d\omega}{m}\sum_{i=1}^m\left(\Phi(a_id\omega)-\dfrac{1}{2}\right)\geq d\omega\left(\Phi(a_{\min}d\omega)-\dfrac{1}{2}\right)(\geq0).\nonumber
	\]
	Using Polya's approximation, we have
	\[
	-\text{E}(V)\geq\dfrac{d\omega}{2}\sqrt{1-\exp\left\{-\dfrac{2}{\pi}(a_{\min}d\omega)^2\right\}}(1+\delta_0)\geq\dfrac{d\omega}{3}\sqrt{1-\exp\left\{-\dfrac{2}{\pi}(a_{\min}d\omega)^2\right\}}.
	\]
	Consider $f(\nu)=e^{-2\nu}+2\nu e^{-\nu}$. Taking derivatives we have $f'(\nu)=2e^{-\nu}(1-2\nu-e^{-\nu})$. Then consider $g(\nu)=1-2\nu-e^{-\nu}$. Taking derivatives we have $g'(\nu)=-2+e^{-\nu}<0$ when $\nu\geq0$, so $g(\nu)<g(0)=0$ when $\nu>0$. Thus $f'(\nu)\leq 0$ when $\nu\geq0$. So $f(\nu)\leq f(0)=1$ when $\nu>0$. Since $(a_{\min}d\omega)^2/\pi>0$, we have $f((a_{\min}d\omega)^2/\pi)\leq1$, i.e.,
	\[
	1-\exp\left\{-\dfrac{2}{\pi}(a_{\min}d\omega)^2\right\}\geq\dfrac{2}{\pi}(a_{\min}d\omega)^2\exp\left\{-\dfrac{1}{\pi}(a_{\min}d\omega)^2\right\}.
	\]
	Thus the expectation of $V$ is bounded from above by a negative constant, i.e.,
	\begin{equation}\label{appendix.E_V}
	-\text{E}(V)\geq\sqrt{\dfrac{2}{9\pi}}a_{\min}d^2\omega^2\exp\left\{-\dfrac{1}{2\pi}(a_{\min}d\omega)^2\right\}.
	\end{equation}

\bibliographystyle{chicago}
\bibliography{proportion_ref}

\end{document}